\documentclass[12pt, a4article]{amsart}   	
\usepackage[top=1.17in, bottom=1.17in, left=1.19in, right=1.19in]{geometry}

\usepackage[colorlinks=true, pdfstartview=FitV, linkcolor=blue, citecolor=blue, urlcolor=blue, breaklinks=true]{hyperref}
\usepackage{amssymb,amscd,graphics, amsfonts, multicol}
\usepackage{graphics}
\usepackage{stmaryrd}
\usepackage{DotArrow}
\usepackage{amsmath, amsthm,wasysym}
\usepackage{epsf}
\usepackage{xypic}
\usepackage{tikz-cd}
\usepackage{color}
\definecolor{R}{rgb}{1, 1, 1}
\definecolor{P}{rgb}{1, 1, 1}
\definecolor{P2}{rgb}{1, 1, 1}
\definecolor{P3}{rgb}{1, 1, 1}
\definecolor{X}{rgb}{1, 0, 1}

\newcommand{\lam}{{\lambda}}

\newcommand{\charc}{\hbox{\rm Char}}

\newcommand{\Z}{{\mathbb Z}}

\def\v{{\mathbf v}}

\def\a{{\alpha}}

\def\0{{\theta}}
\def\tR{{\tilde R}}

\DeclareMathOperator{\supp}{supp}

\title[Seshadri stratification for Schubert varieties]{Seshadri stratification for Schubert varieties and Standard Monomial Theory}
\author{Rocco Chiriv\`i}
\address{Dipartimento di Matematica e Fisica ``Ennio De Giorgi'', Universit\`a del Salento, Lecce, Italy}
\email{rocco.chirivi@unisalento.it}

\author{Xin Fang} 
\address{Department Mathematik/Informatik, Universit\"at zu K\"oln, 50931, Cologne, Germany}
\email{xinfang.math@gmail.com}

\author{Peter Littelmann}
\address{Department Mathematik/Informatik, Universit\"at zu K\"oln, 50931, Cologne, Germany}
\email{peter.littelmann@math.uni-koeln.de}

\theoremstyle{plain}
\newtheorem{theorem}{Theorem}[section]
\newtheorem{coro}[theorem]{Corollary}
\newtheorem{proposition}[theorem]{Proposition}
\newtheorem{lemma}[theorem]{Lemma}

\newtheorem{thmintro}{Theorem}

\theoremstyle{definition}

\newtheorem{remark}[theorem]{Remark}
\newtheorem{example}[theorem]{Example}
\newtheorem{definition}[theorem]{Definition}
\makeindex



\begin{document}
\maketitle
\begin{center}
\textit{To the memory of C.\ S.\ Seshadri, a friend and a wonderful teacher}
\end{center}

\begin{abstract}
{The theory of Seshadri stratifications has been developed by the authors with the intention to build up a new geometric approach towards a standard monomial theory for embedded projective varieties 
with certain nice properties. In this article, we investigate the 
Seshadri stratification on a Schubert variety arising from its 
Schubert subvarieties. We show that the standard monomial theory 
developed in \cite{L1} is compatible with this new strategy.}
\end{abstract}

\section{Introduction}
According to Seshadri, \emph{standard monomial theory} 
(for short \emph{SMT}) 
deals with the construction 
of nice bases of finite-dimensional representations of a semi-simple algebraic group $G$, and more 
generally of what are called Demazure modules in a symmetrizable Kac-Moody setup, see \cite{S2}.
The term \emph{standard monomial theory} might be a bit misleading  because there is no 
definition of what a \emph{standard monomial theory} really is. Despite this ambiguity,
the idea of a standard monomial theory stimulated a wealth of research projects  
connecting the geometry of flag varieties with many aspects in representation theory, 
K-theory and algebraic combinatorics.

Already in \cite{LSIV}, see also  \cite{S2},
Seshadri, Musili and Lakshmibai pointed out the close connection between a nice filtration
of the vanishing ideal $\mathcal I(\partial X(\tau))$ in $\mathcal O_{X(\tau)}$, standard monomial theory,
and their indexing system by admissible pairs. Here $X(\tau)$ is a Schubert variety in a (partial) flag variety 
$G/Q$ for a { maximal parabolic subgroup of classical type $Q$},
and $\partial X(\tau)$ is the reduced union of all Schubert varieties 
properly contained in $X(\tau)$.  This indication of a connection between standard monomial 
theory and valuation theory  { was the starting point for \cite{FL}. 
In \emph{ibidem}, the second and third author 
tested out in the simplest (but nontrivial) case, the
\emph{Grassmann variety},
how to join ideas from standard monomial theory and associated
semi-toric degenerations \cite{Ch,DEP,S2} together with the theory of Newton-Okounkov
bodies \cite{KK} and its associated toric degenerations \cite{And}.
The next aim was to get a Newton-Okounkov type interpretation of 
standard monomial theory and its associated combinatorics, the
Lakshmibai-Seshadri path model \cite{L2}.}

The task to find a setup like in \cite{FL} for a larger class of embedded projective varieties $X\hookrightarrow \mathbb P(V)$ was 
accomplished in \cite{CFL} in a surprisingly general framework: \emph{Seshadri stratifications}.

{ 
A Seshadri stratification of an irreducible embedded projective variety $X\hookrightarrow \mathbb P(V)$ comes with a collection of projective subvarieties 
$X_p$ in $X$, indexed by a finite set $A$.
The set $A$ inherits naturally a partial order 
$\le$ defined by: for $p, q \in A$, $p\le q$ 
if and only if $X_p \subseteq X_q$. In addition, one 
needs { for} each $p\in A$ a homogeneous function 
$f_p\in {\textrm{Sym}}(V^*)$ of degree larger or
equal to one. These subvarieties and functions {have} 
to satisfy certain compatibility conditions, 
see Definition~\ref{Defn:SS}. The example 
studied in this article is the case of a 
(generalized) Schubert variety embedded in 
the projective space over a Demazure module. 
The subvarieties are the Schubert 
subvarieties which are contained in the given 
one, and the functions are the extremal weight 
vectors in the dual of the Demazure module 
(see Section~\ref{section:Schubert_varieties}).

Fix a total order ``$\ge_t$'' on $A$ 
(compatible with the given partial order, see 
Section~\ref{higher_valued_valuation}). 
The existence of a Seshadri stratification on 
$X$ implies the existence of a quasi-valuation 
(which depends on the choice of ``$\ge_t$'') 
with values in $\mathbb Q^A$.  
The quasi-valuation in turn defines a 
filtration of the homogeneous coordinate {ring}
$\mathbb K[X]$, which can be used to 
construct a semi-toric degeneration of $X$ 
(see Section~\ref{maxchainAndquasivalu}). 
The degenerate variety is a reduced union 
of equidimensional projective toric varieties, 
one irreducible component for each maximal 
chain in $A$. The geometry of the degenerate variety
is completely governed by the image $\Gamma\subseteq \mathbb Q^A$ of the quasi-valuation, which is a finite union of finitely generated monoids, one for each maximal chain in $A$. This is the reason why $\Gamma$ is called the associated \emph{fan of monoids}.

{ 
Now having this tool, there are several natural tasks to accomplish. In this sense the present paper is one in a series of articles. 
Particularly nice Seshadri stratifications are those which are balanced and normal \cite{CFL}, it is described in \textit{ibidem} how to 
construct a standard monomial theory in this case.
Other aspects of normal stratifications are discussed in \cite{CFL3}: for example, a Gr\"obner basis of the defining ideal of 
the semi-toric degeneration of the variety can be lifted to define the embedded projective variety.}

Another task is to investigate the
compatibility with already established structures 
like Hodge algebras \cite{DEP}, LS-algebras 
\cite{Ch}, or the \emph{SMT} developed in 
\cite{L1}. The case of the Grassmann variety and 
its classical standard monomial theory was 
studied in \cite{FL}, which was the starting 
point for us to develop the notion of a Seshadri 
stratification. The close connection of Seshadri stratifications with LS-algebras is studied in \cite{CFL2}.  

{ 
The connection between
the \emph{SMT} in \cite{L1} and Seshadri stratifications is one topic of 
this article. The approach we take here is a mixture of the methods developed
in \cite{CFL} and \cite{L1}. We start by defining a basis $\mathbb B(V(\lambda)_\tau)$ of $V(\lambda)^*_\tau$,
the dual of the Demazure module in which we have embedded 
$X(\tau)\hookrightarrow \mathbb P(V(\lambda)_\tau)$ the Schubert variety.
This basis has the following nice property: if $\mathcal V$ is a quasi-valuation associated to the Seshadri
stratification (which depends on the choice of a total order $\ge_t$), then the value of $\mathcal V$ on 
the elements in $\mathbb B(V(\lambda)_\tau)$ is independent {of the choice of the} total order!

This is the starting point for the SMT in Section~\ref{standard-monomial-theory}.
But note that the basis $\mathbb B(V(\lambda)_\tau)$ is far
from being unique. Roughly speaking:
we define a filtration indexed by a partially ordered set such that the associated graded components, which are called leaves, are one-dimensional.
Choosing a nice basis $\mathbb B(V(\lambda)_\tau)$
corresponds to picking representatives of these leaves. It would be interesting to get 
a representation theoretic interpretation of this combinatorially defined collection of subspaces. 

Now the connection with  \cite{L1} is the fact that the basis of $V(\lambda)^*_\tau$ constructed in
\textit{ibidem} is an example of a such a nice basis $\mathbb B(V(\lambda)_\tau)$. This is proved 
in Appendix III, Section~\ref{The_path_vectors}.
To see that the \emph{SMT} for generalized flag varieties  developed in \cite{L1} can be viewed 
as a special example  for the theory developed in \cite{CFL},
it remains to see that the notion of a \textit{standard monomial} is the same. The proof
of this purely combinatorial part can be found in Appendix I, Section~\ref{appendix1}.
Here we show that the fan of monoids coming from our Newton-Okounkov type construction is a reincarnation 
of the path model given by the LS-paths \cite{L2}. So one gets in addition the desired geometric interpretation of the 
LS-path model  as sequences of  renormalized vanishing multiplicities 
of certain functions.}

Another task is to use the new {algebro-geometric} tool to construct in a purely geometric way a standard monomial theory for generalized flag varieties, avoiding the representation theoretic tools used in \cite{L1}.
{ This aspect is discussed in \cite{CFL4}, where a standard monomial theory 
for Schubert varieties is constructed exploiting: the geometry of the Seshadri stratifications of Schubert
varieties by their Schubert subvarieties and the combinatorial LS-path character formula for Demazure modules.} {In the same article, the general theory of Seshadri stratifications is somehow improved by: (1) allowing arbitrary linearizations in the definition of the quasi-valuation (see Section~2.6 in \cite{CFL4}) and (2) weakening the definition of balanced Seshadri stratification (see Section 2.9 in \cite{CFL4}).}
}

{
We give now a short summary of the results of this article.}
The wording of the following theorems is rather informal, the precise formulations can be found later in this article:

\begin{thmintro}
The Seshadri stratification of an embedded Schubert variety
$X(\tau)\hookrightarrow \mathbb P(V(\lambda))$ arising from its Schubert subvarieties
is balanced and normal (see Section~\ref{standard-monomial-theory}).
\end{thmintro}

The main part of this article is dedicated to the proof of the following theorem. As in the case of Newton-Okounkov theory,
one of the main problems is to get an explicit description of the semigroups associated to the (quasi-)valuation.
\begin{thmintro}\label{Thm:B}
\begin{enumerate}
\item[\textit{1)}]  The associated quasi-valuation $\mathcal V$ on the homogeneous coordinate ring $\mathbb K[\hat X(\tau)]$
has as fan of monoids a lattice type realization of the LS-path model (Section~\ref{standard-monomial-theory}). 
\item[\textit{2)}] The standard monomial theory for $X(\tau)\hookrightarrow \mathbb P(V(\lambda))$ constructed in \cite{L1}
is an example of a standard monomial theory for a balanced and normal Seshadri stratification in 
\cite{CFL} (Section~\ref{standard-monomial-theory2}).
\end{enumerate}
\end{thmintro}

{Part 2) of the above theorem is proved also in \cite{CFL2} in the algebraic context of LS-algebras.}

By the general theory of Seshadri stratifications one recovers known results on semi-toric degenerations
of Schubert varieties and their properties \cite{Ch}, and polytopal realizations of the LS-path models \cite{De}
turn up as Newton-Okounkov simplicial complexes, see Sections~\ref{standard-monomial-theory2}, \ref{section:projective_normality} and \ref{section:Newton_Okounkov}: 

\begin{thmintro}
\begin{enumerate}
\item[\textit{1)}] There exists a flat degeneration of $X(\tau)$ to $X_0$, a union of projective toric varieties. Moreover, the special fibre 
$X_0$ is equidimensional, it is isomorphic to $\mathrm{Proj}(\mathrm{gr}_{\mathcal V}\mathbb K[\hat X(\tau)])$,
and its irreducible components are normal toric varieties. 

\item[\textit{2)}] Schubert varieties are projectively normal, so one has a standard monomial theory also for the ring of sections
$\bigoplus_{s\ge 0} \mathrm{H}^0(X(\tau),\mathcal L_{s\lambda})$. 
\item[\textit{3)}] The Newton-Okounkov simplicial complex associated to the quasi-valuation coincides with the 
polytope with an integral structure defined in \cite{De}.
\end{enumerate}
\end{thmintro}


The article is structured as follows: we start in Section~\ref{A_partially_ordered_collection} with a quick review
on the concept of a Seshadri stratification. In Section~\ref{section:Schubert_varieties} we present the setup:
Schubert varieties embedded into the projective space over a Demazure module for a symmetrizable Kac-Moody group over an 
algebraically closed field $\mathbb K$. In particular, we show that the stratification by Schubert subvarieties
is a Seshadri stratification.

In Section~\ref{maxchainAndquasivalu} we recall the construction of the associated quasi-valuation in the
case of Schubert varieties. To determine in such a general setting an explicit description of 
the image of the quasi-valuation (the fan of monoids), one needs to make a { good guess. The favorite candidate}
is presented in Section~\ref{CandidatForMonoid}. The connection with the path model and a conjecture by Lakshmibai
is presented in Appendix I (Section~\ref{appendix1}).

{ We present in Section~\ref{special-functions}
a new approach towards a construction of a basis of the dual $V(\lambda)^*_\tau$  of the Demazure module $V(\lambda)_\tau$,
these elements are called 
path vectors. Note that the notion of a path vector in this article is more general than {that} in \cite{L1}.
So for the convenience of the reader we give in Appendix II and III (Sections~\ref{A_filtration_on_V} and \ref{The_path_vectors}) a 
detailed and adapted review of the constructions in \cite{L1}. We have rewritten the formulation and the proofs in 
a way more adjusted to the point  of view of this article. 

As a consequence we describe in Section~\ref{standard-monomial-theory} a standard monomial theory.
In addition, we show that the Seshadri stratification is balanced and normal.
In Section~\ref{standard-monomial-theory2} we {discuss} some applications: Koszul property, Khovanskii basis,  
compatibilty of the standard monomial theory with the strata, compatibilty with the SMT in \cite{L1}, straightening relations.

In Section~\ref{section:projective_normality} we give (yet another) 
proof of the projective normality of embedded Schubert varieties.}
{The associated 
Newton-Okounkov simplicial complex is described in Section~\ref{section:Newton_Okounkov}. Here we recover the polytopes with integral structure constructed by Dehy in \cite{De} in connection with the LS-path model.

The authors would like to the thank the referee, whose suggestions helped to improve
the presentation of the article.
}

\vskip 0.2cm

We learned so much of what we know about flag and Schubert varieties from Seshadri.
The stimulating, humorous and encouraging discussions with him have always been  
a special event.

\section{Seshadri stratifications}\label{A_partially_ordered_collection}

Throughout the article we fix $\mathbb{K}$ to be an algebraically closed field. 

For a partially ordered set $(A,\leq)$ (in the following called a \textit{poset}) and $p\in A$, we denote by 
${A_p}$ the subset $\{q\in A\mid q\leq p\}$: $(A_p,\leq)$ is again a poset.

Let $V$ be a finite dimensional vector space over $\mathbb{K}$. The hypersurface defined as the vanishing set of a homogeneous polynomial function $f\in\mathrm{Sym}(V^*)$ will be denoted by $\mathcal{H}_f:=\{[v]\in\mathbb{P}(V)\mid f(v)=0\}$.
\index{$\mathcal{H}_f$, vanishing set of a homogeneous polynomial function $f$}

{A variety $X$ is always assumed to be irreducible.}
The affine cone in $V$ over a projective variety $X\subseteq\mathbb{P}(V)$ will be denoted by $\hat{X}$.

\subsection{Seshadri stratifications}


We start with a brief recollection on Seshadri stratifications, as introduced in \cite{CFL}.

Let $X\subseteq \mathbb P(V)$ be an embedded projective variety with graded
homogeneous coordinate ring $R:=\mathbb K[\hat{X}]$.  Let $X_p$, $p\in A$, be a collection of projective subvarieties $X_p$ in $X$, indexed by a finite set ${A}$. The set $A$ inherits naturally a partial order $\leq$ defined by: for $p,q\in A$, $p\leq q$ if and only if $X_p\subseteq X_q$. We assume that there exists a unique maximal element $p_{\max}$ in $A$ with $X_{p_{\max}}=X$.

For each $p\in A$, we fix a homogeneous function $f_{p}\in\mathrm{Sym}(V^*)$ of degree larger or equal to one.

\begin{definition}[\cite{CFL}]\label{Defn:SS}
The collection of subvarieties ${X_p}$ and homogeneous functions ${f_p}$ for $p\in A$ is termed a \emph{Seshadri stratification}, if the following conditions are satisfied:\index{Seshadri stratification}
\begin{enumerate}
\item[(S1)] the projective varieties $X_p$, $p\in A$, are smooth in codimension one; if $q<p$ is a covering relation in $A$, then $X_q\subseteq X_p$ is a codimension one subvariety; 
\item[(S2)] for any $p\in A$ and any $q\not\leq p$, the function $f_q$ vanishes on $X_p$;
\item[(S3)] for $p\in A$, the set-theoretical intersection satisfies 
$$\mathcal{H}_{f_p}\cap X_p=\bigcup_{q\text{ covered by }p} X_q.$$
\end{enumerate}
In a Seshadri stratification, the functions ${f_p}$ are called \emph{extremal functions}.
\end{definition}

\begin{remark}
In the axiom (S3), if $p\in A$ is a minimal element, there is no element covered by $p$. 
Hence the right hand side is the empty set.
\end{remark}

Throughout the paper, the following lemma will be used often without mention.

\begin{lemma}[{\cite{CFL}}]\label{Lem:SS}
If the collection of subvarieties ${X_p}$ of $X$ and homogeneous functions ${f_p}$ for $p\in A$ defines 
a Seshadri stratification, then
\begin{enumerate}
    \item[\it i)] the function $f_p$ does not identically vanish on $X_p$,
    \item[\it ii)] all maximal chains in $A$ have the same length, which coincides with $\dim X$. In particular, the poset $A$ is graded.
\end{enumerate}
\end{lemma}


\begin{definition}\label{lenghtfunction}
Let $p\in A$. The \emph{length} ${\ell}(p)$ of $p$ 
is the length of a (hence any) maximal chain joining $p$ with a minimal element in $A$. 
\end{definition}
According to the above lemma, the length is well-defined and satisfies $\ell(p)=\dim X_p$.

\begin{example}[{\cite{CFL}}]\label{induction}
For a fixed $p\in A$, the collection of varieties $X_q$ and the extremal functions 
$f_q$ for $q\in A_p$ satisfies the conditions (S1)-(S3), and hence defines a Seshadri stratification for $X_p\hookrightarrow  \mathbb P(V)$.
\end{example}

\begin{remark}\label{Rmk:Extend}
Later in the article we will consider the affine cones of the subvarieties in a Seshadri stratification. It is useful to extend the notation one step further. For a minimal element $p\in A$, the affine cone $\hat{X}_p\cong \mathbb{A}^1$. We set $\hat{A}:=A\cup\{p_{-1}\}$
\index{$\hat{A}$, extended poset}
with $\hat{X}_{p_{-1}}:=\{0\}\in V$. Since the variety $\hat{X}_{p_{-1}}$\index{$\hat{X}_{p_{-1}}$} is contained in the affine cone $\hat{X}_p$ for any minimal element $p\in A$, the set $\hat{A}$ inherits a poset structure by requiring $p_{-1}$ to be the unique minimal element.
\end{remark}

\subsection{A Hasse diagram with bonds}\label{Hasse}

We associate an edge-coloured directed graph to a Seshadri stratification of a projective variety $X$ consisting of subvarieties $X_p$ and extremal functions $f_p$ for $p\in A$.

The Hasse diagram ${\mathcal G}_A$ of the poset $A$ is a directed graph on $A$ whose edges are covering relations, pointing to the larger element.

For a covering relation $p>q$ in $A$, $\hat X_q$ is a prime divisor in $ \hat X_p$. According to (S1), the local ring $\mathcal O_{\hat X_p,\hat X_q}$ is a discrete valuation ring. Let ${\nu_{p,q}}:\mathcal O_{\hat X_p,\hat X_q}\setminus\{0\}\to\mathbb{Z}$ be the associated valuation (see also Section \ref{Sec:Valuation}). Let $R_p:=\mathbb{K}[\hat{X}_p]$ denote the homogeneous coordinate ring of $X_p$. For $f\in R_p\setminus\{0\}$, the value $\nu_{p,q}(f)$ is the \textit{vanishing multiplicity} of $f$ in the divisor $\hat X_q$. The integer $b_{p,q}:=\nu_{p,q}(f_p)$ will be called the \textit{bond} between $p$ and $q$. By (S3), we have $b_{p,q}\geq 1$.

The Hasse diagram with bonds\index{Hasse diagram with bonds} is the diagram with edges coloured with the corresponding bonds: 
$q\stackrel{b_{p,q}}{\longrightarrow}p$.

\begin{remark}\label{Rmk:Extend2}
We extend the construction to the poset $\hat{A}$ (Remark \ref{Rmk:Extend}) and the associated extended Hasse diagram $\mathcal{G}_{\hat{A}}$. For a minimal element $p\in A$, the bond $b_{p,p_{-1}}$ is defined to be the vanishing multiplicity of $f_{p}$ at $\hat{X}_{p_{-1}}=\{0\}$, which coincides with the degree of $f_{p}$.
\end{remark}

\subsection{The generic hyperplane stratification}
The definition of a Seshadri stratification looks restrictive, but it has been shown in \cite{CFL}:
\begin{proposition} 
Every embedded projective variety $X\subseteq\mathbb P(V)$, smooth in codimension
one, admits a Seshadri stratification.
\end{proposition}
The proof uses Bertini's Theorem \cite{Jou}. If $X$ is of dimension $r$, then one can find 
linear functions $f_r,\ldots,f_1\in V^*$, so that the inductively defined intersections: $X_r=X$,
and for $j=1,\ldots,r-1$: $X_{j}=X_{j+1}\cap \mathcal H_{f_{j+1}}$, are reduced varieties which are moreover smooth in codimension one.
Now $X_1$ is a smooth curve, and one can find $f_1\in V^*$ such that the intersection $X_{0}=X_{1}\cap \mathcal H_{f_{1}}$ is a union 
of $s=\deg X$ many points. By choosing appropriate functions $f_{0,1}, \ldots, f_{0,s}$ 
for the ($0$-dimensional) irreducible components $X_{0,1},\ldots,X_{0,s}$, one gets in this way a Seshadri stratification for $X$.

\section{Seshadri stratification on Schubert varieties}\label{section:Schubert_varieties}

The Seshadri stratification on an embedded flag variety  $G/B\hookrightarrow \mathbb P(V(\lambda))$
consisting of Schubert subvarieties and the extremal weight functions is briefly discussed in \cite{CFL}. 
In this section we introduce these stratifications  on Schubert varieties in full generality.


\subsection{}
Let $A$ be a symmetrizable generalized  Cartan matrix and
let $\{\alpha_1,\ldots,\alpha_n\}$ be the set of simple roots of the root system. 
We denote by $s_i\in W$ the simple reflection in the Weyl group associated to the simple root $\alpha_i$.
Denote by $G$ the associated 
maximal Kac-Moody group \cite{Marquis, Ma1}. This is in general not an algebraic group but an ind-scheme. 
The group $G$ comes equipped with a (positive) Borel subgroup $B\subseteq G$ and a maximal torus 
$T\subseteq B$. The unipotent radicals of $B$ and its opposite will be denoted by $U$ and $U^-$.
The action of $B$ on $G$ from the {right is free}, which makes it possible to talk about the flag variety $G/B$ as in the 
finite type case \cite{Ma1}. We assume the root datum to be simply connected.

Denote by $\Lambda=\Lambda(T)$ the weight lattice and, according to the choice of $B$, let $\Lambda^+$ be the monoid of integral dominant weights. 
For $\lambda\in \Lambda^+$ let  $V_{\mathbb C}(\lambda)$ be
the irreducible highest weight  representation for the complex version $\mathfrak g_{\mathbb C}$
of the associated Kac-Moody algebra. Fixing a highest weight vector $v_\lambda\in V_{\mathbb C}(\lambda)$, 
we get an admissible lattice $V_{\mathbb Z}(\lambda)\subseteq V_{\mathbb C}(\lambda)$ by applying the 
Kostant integral $\mathbb Z$-form $U(\mathfrak g_{\mathbb C})_{\mathbb Z}$ of the enveloping algebra 
of $\mathfrak g_{\mathbb C}$ to the fixed highest weight vector  $v_\lambda\in V_{\mathbb C}(\lambda)$. 
The tensor product $V_{\mathbb Z}(\lambda)\otimes_{\mathbb Z}\mathbb K$ with the field $\mathbb K$ 
gives the desired $G$-module $V(\lambda)$ of highest weight $\lambda$. { The module $V(\lambda)$ is called the Weyl module associated to the dominant weight $\lambda$. 
See \cite{AndHH,Hum,Ma1} for more information about Weyl modules,
their universal property and a realization as the 
dual representation of the module of global sections 
$\mathrm{H}^0(G/B,\mathcal L_{\lambda^*})$,
where the dominant weight $\lambda^*$ is appropriately chosen.}

Let $Q\subseteq G$ be the parabolic subgroup generated by $B$ and the root subgroups $U_{-\alpha}\subseteq G$
associated to simple roots $\alpha$ such that $\langle \lambda,\alpha^\vee\rangle=0$. 
Denote by $W_Q$ the Weyl group of $Q$.
We associate to $\tau\in W/W_Q$ an extremal weight vector $v_\tau\in V_{\mathbb Z}(\lambda)$ of weight 
$\tau(\lambda)$, see \eqref{extremalweightvector} for an explicit construction.

According to the choice of $B\subseteq G$ let $\mathfrak b_{\mathbb C}\subseteq \mathfrak g_{\mathbb C}$ be a 
Borel subalgebra, and let $U(\mathfrak b)_{\mathbb Z}\subseteq U(\mathfrak g)_{\mathbb Z}$ 
be a corresponding Kostant integral $\mathbb Z$-form of the enveloping algebra. Let $U(\mathfrak b)_{\mathbb K}$ be the algebra obtained by tensoring with $\mathbb K$. 

We denote by $V(\lambda)_\tau$ the corresponding Demazure module: it is the linear span 
$V(\lambda)_\tau=\langle B\cdot v_\tau\rangle$ of the $B$-orbit through $v_{\tau}$. The Demazure module is 
finite dimensional, and it is equipped with an action 
of $B$. More precisely, there exists a normal subgroup $B'\subseteq B$, such that $B'$ acts trivially 
on $V(\lambda)_\tau$ and $B/B'$  is a finite dimensional affine algebraic group acting on $V(\lambda)_\tau$. The orbit closure 
$X(\tau)=\overline{B\cdot[v_{\tau(\lambda)}]}\subseteq \mathbb P(V(\lambda)_\tau)$ is the Schubert variety associated to $\tau$.
This construction of the Schubert variety is independent of the choice of $\lambda$ \cite{Ma1,Ma2}, i.e. for any
dominant weight $\lambda'$ having the same associated parabolic subgroup $Q$, the above construction
leads to an isomorphic variety. In particular, the Schubert variety contains only finitely many $B$-orbits.

\subsection{The homogeneous coordinate ring of Schubert varieties}
{ 
The following connection between Demazure modules 
and the homogeneous coordinate ring of Schubert varieties is well known:
\begin{lemma}\label{CoordDemazure}
We have a graded epimorphism of $B$-modules:
$$
\bigoplus_{s\ge 0} V(s\lambda)_\tau^*\rightarrow \mathbb K[\hat{X}(\tau)]=\bigoplus_{s\ge 0}\mathbb K[\hat{X}(\tau)]_s .
$$
\end{lemma}
\begin{proof}
For $s\ge 1$, we consider the diagonal embedding of $X(\tau)$ into its $s$-fold product with itself,
followed by the Segre embedding:
\begin{equation}\label{Demazure1}
\begin{tikzcd}
X(\tau) \arrow[r, hook] \arrow[d, hook]     & X(\tau)\times\ldots\times X(\tau) \arrow[d, hook] & \\
\mathbb P(V(\lambda)_\tau) \arrow[r, hook] & \mathbb P(V(\lambda)_\tau)\times \ldots\times \mathbb P(V(\lambda)_\tau) \arrow[r, hook] & \mathbb P(V(\lambda)_\tau^{\otimes s}).
\end{tikzcd}
\end{equation}
Every homogeneous function of degree $s$ on $V(\lambda)_\tau$ can be written as a composition
of the map $v\mapsto v^{\otimes s}$ and a linear function on $V(\lambda)_\tau^{\otimes s}$. So we get by restriction
a surjective morphism $(V(\lambda)_\tau^{\otimes s})^*\rightarrow \mathbb K[V(\lambda)_\tau]_s$.

Let $v_\lambda\in V(\lambda)$ and  $v_{s\lambda}\in V(s\lambda)$ be highest weight vectors.
One has a natural $G$-equivariant morphism 
$\phi:V(s\lambda)\rightarrow V(\lambda)^{\otimes s}$, sending $v_{s\lambda}$ to $v_\lambda^{\otimes s}$
and $v_{\tau(s\lambda)}$ to $v_{\tau(\lambda)}^{\otimes s}$ (up to a nonzero scalar multiple).
The morphism $\phi$  induces hence a natural $B$-equivariant morphism
$\phi_\tau:V(s\lambda)_\tau\rightarrow V(\lambda)^{\otimes s}_\tau$. 

The inclusion $X(\tau)  \rightarrow \mathbb P(V(\lambda)_\tau^{\otimes s})$ described in \eqref{Demazure1} is $B$-equivariant, 
and it sends $[v_{\tau(\lambda)}]\in \mathbb P(V(\lambda)_{\tau})$  to the element $[v_{\tau(\lambda)}^{\otimes s}]$  in $\mathbb P(V(\lambda)^{\otimes s}_{\tau})$. 
It follows by the $B$-equivariance that the image of the Schubert variety 
in \eqref{Demazure1} lies in $\mathbb P(\phi_\tau(V(s\lambda)_\tau))\subseteq \mathbb P(V(\lambda)_\tau^{\otimes s})$. 
Since $V(s\lambda)_\tau$ is by definition the linear span of the $B$-orbit $B. v_{\tau(s\lambda)}$, the $B$-equivariance 
implies that $\phi_\tau(V(s\lambda)_\tau)$ is indeed the linear span of the cone over the image of the Schubert variety.

The surjective morphism $(V(\lambda)_\tau^{\otimes s})^*\rightarrow \mathbb K[V(\lambda)_\tau]_s$ induces hence a surjective $B$-equivariant morphism
$V(s\lambda)^*_\tau\rightarrow \phi_\tau(V(s\lambda)^*_\tau)\simeq \mathbb K[\hat X(\tau)]_s$. 
\end{proof}
If $\mathbb K$ is of characteristic zero, then it is well known that the morphism $\phi_\tau$ is in fact  injective. In positive characteristic one can either  use 
Frobenius splitting  (see \cite{AndHH2,Ma1,MR,R,RR}) or the quantum Frobenius splitting and its applications (see \cite{KL1,KL2,L1})
to show that the map $\phi_\tau$ is always injective. 
We will use later, usually without mention, the following representation theoretic interpretation of the 
homogeneous coordinate ring $\mathbb K[\hat{X}(\tau)]$ of $X(\tau)\subseteq \mathbb P(V(\lambda)_\tau)$:
\begin{coro}\label{DemazureIso}
We have a graded isomorphism of $B$-modules:
$$
 \mathbb K[\hat{X}(\tau)]=\bigoplus_{s\ge 0}\mathbb K[\hat{X}(\tau)]_s\simeq  \bigoplus_{s\ge 0} V(s\lambda)_\tau^*.
$$
\end{coro}
}

\subsection{Partial orders}\label{Schubert:Partial_orders}

Let $\ell$ be the length function on the Weyl group $W$. The value $\ell(w)$ is the minimal
length of a reduced decomposition of $w$ as a product of simple reflections. 

We often identify $W/W_Q$ with the subset $W^Q\subseteq W$
of representatives of elements of $W/W_Q$ in $W$ of minimal length. The Weyl group is naturally
endowed with a partial order by viewing the pair consisting of $W$ and the simple reflections as a Coxeter system. 
This partial order is called the Bruhat order on $W$. We get an induced Bruhat order on   $W/W_Q$ via the identification
with the subset $W^Q\subseteq W$. We view the length function $\ell(\cdot)$ as a function on $W/W_Q$ as follows: we define
$\ell(\tau)$ for $\tau\in W/W_Q$ to be $\ell(\hat\tau)$, where $\hat\tau\in W^Q$ is the unique minimal representative of $\tau$. This is the same as the length in the graded poset $W/W_Q$.

The partial order and the length function have the following geometric interpretation: for $\tau\in W/W_Q$, the length $\ell(\tau)$ is the dimension
$\dim X(\tau)$ of the corresponding Schubert variety, and if $\kappa\in  W/W_Q$ is a second element, then
$ X(\kappa)\subseteq X(\tau)$ if and only if $\kappa\le \tau$ in the Bruhat order. 

\subsection{Roots and subgroups}
Let $\Phi$ be the set of real roots of $G$. Having fixed the Borel subgroup $B$, we can divide $\Phi$  into the set of positive real roots $\Phi^+$
and negative real roots $\Phi^-$. Let $\Delta$ denote the set of simple roots. 

For a real root $\beta$ denote by $U_\beta\subseteq G$
the one-dimensional root subgroup corresponding to $\beta$.
The  parabolic subgroup $Q$ is determined by a subset $\Upsilon$ of the simple roots $\Delta$. 
Let $\Phi_{\Upsilon}\subseteq \Phi$ be the subset of real roots spanned by $\Upsilon$. The Weyl group $W_Q$ is the subgroup of $W$  
generated by the simple reflections $s_\alpha$, $\alpha\in \Upsilon$.

Since $\Phi_\Upsilon$ is stable under $W_Q$,  the following
conventions make sense: for $\tau\in W/W_Q$ and
$\gamma\in\Phi^+$ we write $\tau^{-1}(\gamma)\in\Phi_Q$ if
$w^{-1}(\gamma)\in\Phi_Q$ for one, and hence every representative
$w\in W$ of $\tau$. Similarly, write $\tau^{-1}(\gamma)\not\in\Phi_Q$
and $\tau^{-1}(\gamma)\prec 0$ if $w^{-1}(\gamma)\not\in\Phi_Q$
and $w^{-1}(\gamma)\prec 0$ for one, and hence every representative $w$ of $\tau$ in $W$.
For an element $\sigma\in W/W_Q$ set 
$$
\Phi^+_\sigma:=\{\gamma\in\Phi^{+}\mid \sigma^{-1}(\gamma)\not\in\Phi_Q,\ \sigma^{-1}(\gamma)\prec 0 \}.
$$
The set $\Phi^+_\sigma$ is finite and closed under addition, so the product of the root subgroups is a subgroup $U_\sigma$ of $U$.  The decomposition
$U_\sigma=\prod_{\gamma\in \Phi^+_\sigma} U_\gamma$
as product of root subgroups holds for any chosen ordering of the elements in $\Phi^+_\sigma$. 
The orbit map $U_\sigma \rightarrow  \mathbb P(V(\lambda)_\sigma)$, 
$u\mapsto u\cdot [v_\sigma]$, is an isomorphism onto the  Schubert cell $C(\sigma)$, {  viewed as a subset of { $\mathbb P(V(\lambda)_\sigma)$. }}

\subsection{Smooth in codimension one}
Our goal is to endow the embedded Schubert variety $X(\tau) \hookrightarrow  \mathbb P(V(\lambda)_\tau)$
with a Seshadri stratification and use the associated valuations and quasi-valuations
to study the corresponding homogeneous coordinate rings. In this section we start from the requirement of smoothness in codimension one. Such a property is well known, we provide a short proof for completeness.

We fix the Schubert variety $X(\tau)$ in $\mathbb P(V(\lambda)_\tau)$ and let $X(\kappa)\subset X(\tau)$ be a 
Schubert variety of codimension one. There exists a positive real root $\beta$ such that $s_\beta\kappa=\tau$.

\begin{lemma}\label{Lemmaone}
Let $\tau>\kappa$ be two elements in $W/W_Q$ such that $\ell(\tau)=\ell(\kappa)+1$. Let $\beta$ be the positive real root such that $s_\beta\kappa=\tau$.
The product map 
$$
U_\kappa\times U_{-\beta}\rightarrow  \mathbb P(V(\lambda)_\tau), \quad (u,v)\mapsto uv\cdot[v_\kappa],
$$
induces an isomorphism of $U_\kappa\times U_{-\beta}$ with a $T$-stable affine dense subset of $X(\tau)$, containing $\kappa$ as the only $T$-fixed point.
\end{lemma}

\begin{proof}
Denote by $U^-_Q$ the subgroup of $G$ generated by the one-parameter
subgroups $U_\gamma$  corresponding to  negative real roots $\gamma$ which are not in $\Phi_Q$, it is an affine ind-group.
The orbit map $\kappa U_Q^-\kappa^{-1}\rightarrow G/Q$, $\kappa u\kappa^{-1}\mapsto \kappa u \kappa^{-1}\cdot\kappa$, defines an isomorphism
of the subgroup $\kappa U^-\kappa^{-1}\subset G$ onto an open $T$-stable affine ind-subvariety in $G/Q$, containing $\kappa$ as the only $T$-fixed point. 
By definition, $U_\kappa$
is contained in $\kappa U^-\kappa^{-1}$.
Since $U_\beta\not\subset U_\kappa$, it follows $U_{-\beta}\subset \kappa U^-\kappa^{-1}$. So one can regard $U_\kappa\times U_{-\beta}$
as a $T$-stable ($T$-acting by conjugation) affine subspace in $\kappa U^-\kappa^{-1}$, and
 the restriction of the map $\kappa U^-\kappa^{-1}\rightarrow G/Q$ to $U_\kappa\times U_{-\beta}$ induces a  $T$-equivariant isomorphism of the latter to an affine subspace in $G/Q$. 

Let $\mathrm{SL}_2(\beta)$ be the $\mathrm{SL}_2$-copy in $G$ corresponding to the positive real root $\beta$. 
Since
$$
\overline{U_{-\beta}\cdot\kappa}= \mathrm{SL}_2(\beta)\cdot\kappa=\mathrm{SL}_2(\beta)\cdot\tau=\overline{U_{\beta}\cdot\tau},
$$
it follows $U_{-\beta}\cdot\kappa \subset X(\tau)$. The Schubert variety $X(\tau)$ is $U$-stable and hence $U_\kappa U_{-\beta}\cdot\kappa$
is an affine subspace in $X(\tau)$ of the same dimension as $X(\tau)$.
\end{proof} 
{ 
Denote the image of $U_\kappa\times U_{-\beta}$ by $C(\tau,\kappa)$, it is a cell. Let $C(\kappa):=U_\kappa\cdot\kappa\subset X(\kappa)$
and $C(\tau):=U_{\tau}\cdot\tau\subset X(\tau)$ be the usual Schubert cells. By construction one has
$$
C(\kappa)\subset C(\tau,\kappa),\hbox{\rm \ and $C(\tau)\cap C(\tau,\kappa)$ is dense in $X(\tau)$}.
$$
\begin{remark}
The product $U_\kappa\times U_{-\beta}$ is isomorphic to $\mathbb A^{\ell(\tau)}$. By choosing a parameterization of the root subgroups of $U_\kappa$ by parameters $t_\gamma$, 
$\gamma\in \Phi_\kappa^+$, and $U_{-\beta}$ by $t_\beta$, we get a parameterization of an open and dense subset of $X(\tau)$, which contains the open cell $C(\sigma)$ in $X(\sigma)$. The lemma above implies
that a generic element in the image of $U_\kappa\times U_{-\beta}$ in $\mathbb P(V(\lambda)_\tau)$ is of the form 
$$[ct_{\beta}^{\vert\langle \tau(\lambda),\beta^\vee\rangle\vert} v_{\tau(\lambda)}+\text{sum of weight vectors of weight}> \tau(\lambda)].
$$
Here $c$ is a nonzero constant and $v_{\tau(\lambda)}$ is a weight vector of weight $\tau(\lambda)$. Denote by $f_\tau\in (V(\lambda)_\tau)^*$ the corresponding dual vector, i.e. $f_\tau$ is a $T$-eigenvector
of weight $-\tau(\lambda)$, and $f_\tau(v_{\tau(\lambda)})=1$. 
It follows that $f_\tau$ vanishes on $X(\sigma)$ with multiplicity
$\vert \langle \tau(\lambda),\beta^\vee\rangle\vert$.
\end{remark}
}
\begin{coro}\label{codimensionOne}
Schubert varieties are smooth in codimension one.
\end{coro}

\begin{proof}
The Borel subgroup $B$ acts on $X(\tau)$ with a finite number of orbits. The set of singular points is stable under this action,
so we have only to check that the codimension one orbits contain no singular points of $X(\tau)$. But such an orbit 
contains a $T$-fixed point $\kappa$ so that the associated Schubert variety 
$X(\kappa)$ is of codimension one in $X(\tau)$. By the above construction, there exists an open and dense cell $C(\tau,\kappa)$
in $X(\tau)$ containing the orbit $B\cdot\kappa$ and meeting the open orbit $B\cdot\tau$ in a dense subset. It follows that 
the $B$-orbit through $\kappa$ contains only smooth points in $X(\tau)$.
\end{proof}

Let $\kappa,\tau$ und $\beta$ be as in Lemma~\ref{Lemmaone}.
\begin{coro}\label{onemultiplicity}
The following hold:
\begin{itemize}
\item[{\it i)}]  $\beta$ is not in the non-negative real span of $\Phi^+_\kappa$.
\item[{\it ii)}] The $T$-weight spaces in $V(\lambda)_\tau$ of weight $\kappa(\lambda)+\ell\beta$, $\ell\ge 0$, are at most
one-dimensional.
\end{itemize}
\end{coro}
\begin{proof}
An element $\delta$
in the non-negative real span of $\Phi^+_\kappa$ has the property that 
$\kappa^{-1}(\delta)$ is a non-negative real span of negative roots, but  $\kappa^{-1}(\beta)$ is a positive root. 
Since  $V(\lambda)_\tau=\langle U_\kappa U_{-\beta}\cdot v_\kappa\rangle$
by Lemma~\ref{Lemmaone}, it follows that all $T$-weight spaces in $V(\lambda)_\tau$ of weight $\kappa(\lambda)+\ell\beta$, $\ell\ge 0$, are at most
one-dimensional.
\end{proof}

\subsection{A Seshadri stratification for $X(\tau)\hookrightarrow \mathbb P(V(\lambda)_\tau)$}\label{Seshadri_stratification}

As before, let $Q\supseteq B$ be the standard parabolic subgroup of $G$ associated to $\lambda$, i.e. $Q$ is generated by
$B$ and the root subgroups $U_{-\alpha}$ for all simple roots $\alpha$ such that $\langle\lambda,\alpha^\vee\rangle=0$. 
Let $\tau\in W/W_Q$, in the following we consider the Schubert variety $X(\tau)\hookrightarrow \mathbb P(V(\lambda)_\tau)$ embedded in the projective space $\mathbb P(V(\lambda)_\tau)$ associated to the Demazure module $V(\lambda)_\tau$.

The set $A_{\tau}:=\{\sigma\in W/W_Q\mid \sigma\le \tau\}$, endowed with the Bruhat order, is a poset. To each $\sigma\in A_{\tau}$ we associate the Schubert variety $X(\sigma)$,
which is a closed subvariety of $X(\tau)$. So we have a collection of subvarieties $X(\sigma)$ of $X(\tau)$, indexed by 
the partially ordered $A_\tau$ such that $\kappa\le\sigma$ if an only if $X(\kappa)\subseteq X(\sigma)$.
In addition, all the subvarieties are smooth in codimension one by Corollary~\ref{codimensionOne}, and it is well-known
that the covering relations correspond to codimension one subvarieties. The length function defined in Definition~\ref{lenghtfunction} 
corresponds to the usual length function on $W$ as a Coxeter group (see Section~\ref{Schubert:Partial_orders}). 

To get a Seshadri stratification, we need in addition a collection of homogeneous functions $f_\sigma\in \mathbb K[V(\lambda)_\tau]$. 
We have fixed 
for all $\sigma\in A_\tau$ a $T$-eigenvector $v_{\sigma}\in V(\lambda)_\tau$ of weight $\sigma(\lambda)$.
The corresponding weight space is one dimensional, so the vector is unique up to scalar multiple.
Denote by $f_\sigma\in (V(\lambda)_\tau)^*$ the corresponding dual vector, i.e. $f_\sigma$ is a $T$-eigenvector
of weight $-\sigma(\lambda)$, and $f_\sigma(v_{\sigma})=1$.

\begin{proposition}\label{SeshadriStraficication_Schubert}
The collection of subvarieties $X(\sigma)$ and linear functions $f_\sigma$, $\sigma\in A_\tau$,
defines a Seshadri stratification for $X(\tau)$.
\end{proposition}

\begin{proof}
By Corollary \ref{codimensionOne}, it remains to prove that the functions satisfy the conditions (S2) and (S3). To show (S2), we prove in fact a bit more: 
$f_\kappa\vert_{X(\tau)}$ is identically $0$ if and only if $\kappa\not\le\tau$.

If $f_\kappa\vert_{X(\tau)}\not\equiv 0$, then there exists an element 
$z\in  C(\tau)$ such that $f_\kappa(z)\not=0$. Let $\hat z\in \hat X(\tau)\setminus\{0\}$ be an element in the affine cone over $X(\tau)$ 
such that $[\hat z]=z$, and let $\hat z=\sum_{\mu\in\Lambda} w_\mu $ be a decomposition into weight vectors. 
Since $f_\kappa(\hat z)\not=0$, we know $w_{\kappa(\lambda)}\not=0$. Now $\kappa(\lambda)$ is an extremal weight
(it is a vertex of the polytope obtained as the convex hull of all weights occuring in $V(\lambda)_\tau$), so there
exists a one-parameter subgroup $\eta:\mathbb{K}^*\rightarrow T$ such that $\lim_{t\rightarrow 0}\eta(t)z=[v_\kappa]$. 
We find $[v_\kappa]\in \overline{C(\tau)}=X(\tau)$, which implies $X(\kappa)\subseteq X(\tau)$ and hence $\kappa\le \tau$.
Vice versa, suppose $\kappa\le \tau$. But then we know $f_\kappa\vert_{C(\kappa)}\not\equiv 0$ by construction.
Since $X(\kappa)\subseteq X(\tau)$, it follows $f_\kappa\vert_{X(\tau)}\not\equiv 0$.

Note that $f_\tau\in (V(\lambda)_\tau)^*$ spans a $B$-stable line, so the intersection $\mathcal{H}_{f_\tau}\cap X(\tau)$
is closed and $B$-stable. Indeed, $f_\tau(z)\not=0$ for $z\in X(\tau)$ if and only  if $[z]\in C(\tau)$. 
The intersection is hence the union of all Schubert varieties contained in, but not equal to $X(\tau)$, i.e. it is the union of all 
Schubert varieties of codimension one in $X(\tau)$, which proves (S3).
\end{proof}

The corresponding Hasse diagram with bonds is just the usual Bruhat graph for the set $A_{\tau}$: one has an edge between 
$\sigma,\kappa\in A_{\tau}$, $\sigma>\kappa$, if and only if there exists a positive real root $\beta$ such that $\sigma=s_{\beta}\kappa$ 
and $\ell(\sigma)=\ell(\kappa)+1$. Lemma~\ref{Lemmaone} implies that the bond is equal to 
$b_{\sigma,\kappa}=\langle\kappa(\lambda),\beta^\vee\rangle$, where $\beta^\vee$ is the coroot associated to $\beta$.
The Hasse diagram with bonds for the adjoint representation of $\mathrm{SL}_3$ is depicted below, with the Weyl group $W$ identified with the symmetric group $\mathfrak{S}_3$:
$${\scriptsize
\xymatrix{
& (132)  \ar@{<-}[r]^2\ar@{<-}[rdd]^1&(12)\ar@{<-}[dr]^1&\\
(13)\ar@{<-}[ur]^1\ar@{<-}[dr]_1&&&\mathrm{id}\\
& (123) \ar@{<-}[r]_2\ar@{<-}[ruu]^1&(23)\ar@{<-}[ur]_1&\\
}
}
$$

\begin{remark}\label{extrabond}
As in Remark \ref{Rmk:Extend} and \ref{Rmk:Extend2}, we let $\hat{X}(\tau_{-1})$ denote the origin of $\hat X(\mathrm{id})$: it is contained in  all affine cones $\hat{X}(\tau)$. 
The bond $b_{\mathrm{id},\tau_{-1}}$ is the vanishing multiplicity of $f_{\mathrm{id}}$ at the origin, which is equal to $1$, the degree of 
$f_{\mathrm{id}}$. We get in this way the extended set $\hat A_\tau:=A_\tau\cup\{\tau_{-1}\}$ and the extended Hasse diagram $\mathcal G_{\hat A_\tau}$.
\end{remark}

\section{Maximal chains, valuations and quasi-valuations}\label{maxchainAndquasivalu}
Let $X(\tau)\subseteq G/Q$ be a Schubert variety of dimension $\ell(\tau)=r$. 
Fix a {\it maximal chain} $\mathfrak C$\index{$\mathfrak C$, maximal chain} in $A_\tau$, i.e., we have a sequence of Weyl group elements
\begin{equation}\label{chain}
\mathfrak C:\quad
\tau=\tau_r\mathop{\longleftarrow}^{\beta_r} \tau_{r-1}\mathop{\longleftarrow}^{\beta_{r-1}}
 \tau_{r-2}\mathop{\longleftarrow}^{\beta_{r-2}} 
 \quad\cdots \quad
 \mathop{\longleftarrow}^{\beta_{3}}\tau_2
 \mathop{\longleftarrow}^{\beta_{2}}\tau_1
 \mathop{\longleftarrow}^{\beta_{1}}\tau_0=\mathrm{id},
\end{equation}
such that $\tau_j>\tau_{j-1}$, $\ell(\tau_j)=\ell(\tau_{j-1})+1$ and
the $\beta_j$ are positive roots such that $s_{\beta_j}\tau_j=\tau_{j-1}$ for all $j=1,\ldots,r$. 
Or, in geometric terms, we have a sequence of Schubert varieties
\begin{equation}\label{Schubertsequence}
X(\tau)=X(\tau_r)\supsetneq X(\tau_{r-1})\supsetneq X(\tau_{r-2})\supsetneq \cdots \supsetneq X(\tau_2)\supsetneq X(\tau_1)\supsetneq X(\tau_0)=\mathrm{id},
\end{equation}
successively contained in each other of codimension one, and $s_{\beta_j}\tau_j=\tau_{j-1}$ for $j=1,\ldots,r$.

To such a chain of subvarieties we associate a higher rank valuation (see \cite{CFL}). 


\subsection{Basics on valuations and quasi-valuations}\label{Sec:Valuation}

\begin{definition}\label{quasidef}
Let $\mathcal R$ be a $\mathbb K$-algebra. A {\it quasi-valuation} on $\mathcal R$ with 
values in a totally ordered abelian group $\mathbb{G}$ 
is a map $\nu: \mathcal R\setminus \{0\}\rightarrow \mathbb{G}$ satisfying the following conditions:
\begin{itemize}
\item[{(a)}] $\nu(x + y) \ge \min\{\nu(x), \nu(y)\}$ for all $x, y\in \mathcal R \setminus \{0\}$ with $x+y\not=0$;
\item[{(b)}] $\nu(\lambda x) = \nu(x)$ for all $x\in \mathcal R  \setminus \{0\}$ and $\lambda \in \mathbb K^*$;
\item[{(c)}] $\nu(xy) \ge \nu(x) + \nu(y)$ for all $x,y \in \mathcal R \setminus \{0\}$ with $xy\not=0$.
\end{itemize}
The map $\nu$ is called a \emph{valuation} if  
the inequality in (c) can be replaced by an equality:
\begin{itemize}
\item[{(c')}] $\nu(xy)= \nu(x) + \nu(y)$ for all $x,y \in \mathcal R \setminus \{0\}$ with $xy\not=0$.
\end{itemize}
\end{definition}

Quasi-valuations on $\mathcal{R}$ can be thought of as algebra filtrations on $\mathcal{R}$ (see \cite[Section 2.4]{KM}).
The proof of the following lemma is straightforward, see \cite{CFL}.

\begin{lemma}\label{simpleproperties}
Let $\nu, \nu_1,\ldots,\nu_k:\mathcal{R}\setminus\{0\}\to \mathbb{G}$ be quasi-valuations and let $x, y\in \mathcal R \setminus \{0\}$. 
\begin{enumerate}
\item [{\it i)}] If $\nu(x)\not= \nu(y)$, then $\nu(x + y) = \min\{\nu(x), \nu(y)\}$.
\item[{\it ii)}] If $x+y\not=0$ and $\nu(x + y)>\nu(x)$, then  $\nu(x)=\nu(y)$.
\item[{\it iii)}] The map $\mathcal{R}\setminus\{0\}\to\mathbb{G}$, $x\mapsto\min\{\nu_j(x)\mid j=1,\ldots,k\}$ defines a quasi-valuation on $\mathcal{R}$.
\end{enumerate}
\end{lemma}

Natural examples of valuations arise from vanishing orders of functions.
By Corollary~\ref{codimensionOne}, Schubert varieties are smooth in codimension one. The local ring 
$\mathcal{O}_{\hat X(\tau_i),\hat X(\tau_{i-1})}$ is hence a discrete valuation ring for all $i=0,\ldots,r$ (see also Section \ref{Hasse}). Since $\mathbb K(\hat X(\tau_i))$ is the quotient field of $\mathcal{O}_{\hat X(\tau_i),\hat X(\tau_{i-1})}$, we obtain a $\mathbb Z$-valued valuation
$$
\nu_i:\mathbb K(\hat X(\tau_i))\setminus\{0\}\rightarrow \mathbb Z.
$$ 
For $g\in \mathbb K(\hat X(\tau_i))\setminus\{0\}$, the value $\nu_i(g)$ is called the vanishing order of $g$ along the divisor $\hat X(\tau_{i-1})$.

\subsection{Higher rank valuations}\label{Hrv}
We come back to the case of a chain of Schubert varieties as in \eqref{Schubertsequence}. In the following we  
consider the affine cones $\hat X(\tau_{j})$, $j=0,\ldots,r$, over the embedded Schubert varieties. 
%

Since $\mathfrak C$ is a chain, the elements in $\mathfrak C$ are totally ordered. The vector space $\mathbb Q^{\mathfrak C}$
has as standard basis the elements $e_{\tau_j}$, $j=0,\ldots,r$. The lexicographic order on $\mathbb{Q}^{\mathfrak{C}}$ is defined as follows: 
{ 
for $\underline{c}=\sum_{j=0}^r c_je_{\tau_j}$ and $\underline{c'}=\sum_{j=0}^r c'_je_{\tau_j}\in\mathbb{Q}^{\mathfrak{C}}$, $\underline{c'}>\underline{c}$ if there exists $0\leq j\leq r$ 
such that $c'_r=c_r$, $\ldots$, $c'_{j+1}=c_{j+1}$ and $c'_j>c_j$.}
From now on we take on $\mathbb Q^{\mathfrak C}$ this lexicographic order.

Let $N$ be the l.c.m. of all bonds appearing in the Hasse diagram with bonds $\mathcal G_{A_\tau}$.

Let $g_r:=g\in\mathbb{K}(\hat{X}(\tau_r))$ be a non-zero rational function with vanishing order { $o_r\in \mathbb Z$ along 
the prime divisor $\hat{X}(\tau_{r-1})$ in $\hat{X}(\tau_r)$. We define a rational function 
$$
h:=\frac{g_r^N}{f_{\tau_r}^{N\frac{o_r}{b_r}}}\in\mathbb{K}(\hat{X}(\tau_r)),
$$
where} $b_r$ is the vanishing order of $f_{\tau_r}$ along $\hat{X}(\tau_{r-1})$. In other words, $b_r$ is the bond in the Hasse diagram
$\mathcal G_{A_\tau}$ associated to the edge joining $\tau_{r}$ and $\tau_{r-1}$. 

The restriction of $h$ to $\hat{X}(\tau_{r-1})$, denoted by $g_{r-1}$, gives rise to a well-defined non-zero rational function in 
$\mathbb{K}(\hat{X}(\tau_{r-1}))$ (\cite[Lemma 4.1]{CFL}). This procedure can be henceforth iterated, yielding a sequence of rational functions $g_{\mathfrak{C}}:=(g_r,g_{r-1},\ldots,g_0)$ with $g_k\in\mathbb{K}(\hat{X}(\tau_k))\setminus\{0\}$. 
{  The vanishing order of $f_{\tau_k}$ on 
$\hat{X}(\tau_{k-1})$ will be denoted by $b_k$.}

In view of the $N$-th powers appearing in the sequence of rational functions, we define a map 
$\mathcal{V}_{\mathfrak{C}}:\mathbb{K}[\hat{X}(\tau)]\setminus\{0\}\to\mathbb{Q}^\mathfrak{C}$ in the following way:
$$
g\mapsto \frac{\nu_r(g_r)}{b_r}e_{\tau_r}+\frac{1}{N}\frac{\nu_{r-1}(g_{r-1})}{b_{r-1}}e_{\tau_{r-1}}+\ldots+\frac{1}{N^r}\frac{\nu_{0}(g_0)}{b_{0}}e_{\tau_0}.
$$

It is proved in \cite[Section 6]{CFL}:
\begin{proposition}
$\mathcal{V}_{\mathfrak{C}}$ is a $\mathbb Q^{\mathfrak C}$-valued valuation on $\mathbb{K}[\hat{X}(\tau)]$.
\end{proposition}

Denote by 
$\mathbb V_{\mathfrak C}(X(\tau))= \{\mathcal V_{\mathfrak C}(f)\mid f\in \mathbb{K}[\hat{X}(\tau)]\setminus\{0\}\}\subseteq  \mathbb Q^{\mathfrak C}$ 
the image of the valuation map $\mathcal V_{\mathfrak C}$, it is a monoid called the \emph{valuation monoid}. 
The valuation $\mathcal{V}_{\mathfrak{C}}$ defines a filtration of $\mathbb{K}[\hat{X}(\tau)]$ by subspaces,
for $\underline a\in \mathbb Q^{\mathfrak{C}}$ we set:
\begin{equation}\label{filtration:1}
{\mathbb{K}[\hat{X}(\tau)]}_{\ge \underline a}^{\mathfrak{C}}:=\{ g\in   \mathbb{K}[\hat{X}(\tau)]\setminus\{0\}\mid  
\mathcal V_{\mathfrak C}(g)\ge \underline a\}\cup\{0\},
\end{equation}
respectively
\begin{equation}\label{filtration:2}
{\mathbb{K}[\hat{X}(\tau)]}_{> \underline a}^{\mathfrak{C}}:=\{ g\in   \mathbb{K}[\hat{X}(\tau)]\setminus\{0\}\mid  
\mathcal V_{\mathfrak C}(g)> \underline a\}\cup\{0\}.
\end{equation}
The subquotient ${\mathbb{K}[\hat{X}(\tau)]}_{\ge \underline a}^{\mathfrak{C}}/{\mathbb{K}[\hat{X}(\tau)]}_{> \underline a}^{\mathfrak{C}}$
is called a \emph{leaf} of the valuation. It has been shown in \cite[Section 6]{CFL}:

\begin{proposition}
The valuation $\mathcal{V}_{\mathfrak{C}}$ has at most one dimensional leaves.
\end{proposition}

\subsection{A higher rank quasi-valuation}\label{higher_valued_valuation}

There is no obvious reason why one maximal chain should better reflect the geometry of the Schubert variety than another. Moreover, it is not at all clear why in general the associated valuation monoid $\mathbb V_{\mathfrak C}(X(\tau))$ should be finitely generated. 

We fix on $A_\tau$ a total order $\le_t$ refining the partial order which respects
the length, i.e. if $\ell(\sigma)>\ell(\kappa)$, then $\sigma>_t\kappa$.\index{$\ge_t$, chosen total order}
The set $\mathbb Q^{A_\tau}$ is thus endowed with a total order by taking the lexicographic order as in Section \ref{Hrv}. Such a total order is compatible with the addition in $\mathbb Q^{A_\tau}$.

In the following we consider for a maximal chain $\mathfrak C$
the vector space $\mathbb Q^{\mathfrak C}$ as a subspace of $\mathbb Q^{A_\tau}$ spanned by the coordinate functions $e_\sigma$, $\sigma\in\mathfrak C$. The total order on $\mathbb Q^{A_\tau}$ induces a total order on $\mathbb Q^{\mathfrak C}$. Note
that this order coincides with the total order on $\mathbb Q^{\mathfrak C}$ fixed in Section~\ref{Hrv}.

Since we regard $\mathbb Q^{\mathfrak C}$ as a subspace of $\mathbb Q^{A_\tau}$, it makes sense to write $\mathcal V_{\mathfrak C}(g)\in \mathbb Q^{A_\tau}$ for a regular function 
$g\in \mathbb{K}[\hat{X}(\tau)]\setminus\{0\}$. 

Denote by $\mathcal C$ the set of all maximal chains in $A_{\tau}$.
By Lemma~\ref{simpleproperties}, the minimum
over a finite list of valuations is a quasi-valuation. 

\begin{definition}\label{defnquasivaluation}
\begin{enumerate}
\item[i)] We define the \emph{quasi-valuation} $\mathcal V:  \mathbb{K}[\hat{X}(\tau)]\setminus\{0\}\rightarrow \mathbb Q^{A_\tau}$ by
$$
\mathcal{V}(g):=\min\{\mathcal V_{\mathfrak C}(g) \mid \mathfrak C\in\mathcal{C}\}.
$$
\item[ii)] For $\underline{a}=(a_\sigma)_{\sigma\in A_\tau}\in\mathbb{Q}^{A_\tau}$, 
the \textit{support} of $\underline{a}$ is defined by
$$
\supp \underline{a}:=\{\sigma\in A_\tau\mid a_\sigma\not=0\}.
$$ 
\end{enumerate}
\end{definition}

\begin{remark}
Let $g\in \mathbb{K}[\hat{X}(\tau)]\setminus\{0\}$. Unless $\supp\mathcal V(g)$ is a maximal chain, there might be several
maximal chains $\mathfrak{C}$ such that $\mathcal V(g)=\mathcal V_{\mathfrak C}(g)$.
\end{remark}

As an example let us consider an extremal function $f_\kappa$ for $\kappa\in A_\tau$.
\begin{lemma}\label{extremalvaluationII}
For any $\kappa\in A_\tau$, $\mathcal V(f_\kappa)=e_\kappa$.
\end{lemma}
\begin{proof}
Let $\mathfrak C:\tau=:\tau_r>\tau_{r-1}>\ldots>\tau_0=\mathrm{id}$ be a maximal chain in $A_\tau$, and denote by $b_i$ the bond between $\tau_i$ and $\tau_{i-1}$, $i=0,\ldots,r$ (recall that, as in Remark \ref{extrabond}, $\tau_{-1}$ is the additional element in the extended set $\hat A_{\tau}$.)
Suppose first that 
$\kappa=\tau_j$ for some $0\le j\le r$. By Lemma~\ref{Lem:SS}, $f_{\kappa}$ does not vanish identically on $\hat X(\tau_i)$ for $i\ge j$.
The inductive procedure to determine $(f_\kappa)_{\mathfrak C}=(g_r=f_\kappa,g_{r-1},\ldots,g_0)$ gives:
$$
g_r=f_{\kappa},\ g_{r-1}=f^{N}_{\kappa},\ \ldots, \textrm{\ and\ } g_{j}=f^{N^{r-j}}_{\kappa}.
$$
The function $g_j=f_{\kappa}^{N^{r-j}}$ vanishes on the divisor $\hat X(\tau_{j-1})$ with the multiplicity $a_j=N^{r-j} b_j$,
where $b_j$ is the bond between $\kappa=\tau_j$ and $\tau_{j-1}$. So we have
$$
g_{j-1}=\frac{g_j^N}{(f_{\tau_j})^{N\frac{a_j}{b_j}}}=\frac{f_{\tau_j}^{N^{r-i+1}}}{f_{\tau_j}^{N^{r-j+1}}} =1.
$$
This function $g_{j-1}$ does not vanish on any of the Schubert varieties, so the procedure implies $g_i=1$ for all $i<j$. 
We have proved:
$$
(f_{\kappa})_{\mathfrak C}=(f_{\kappa},f_{\kappa}^N,f_{\kappa}^{N^2},\ldots,f_{\kappa}^{N^{r-j}},1,\ldots,1).
$$
For the valuation $\mathcal V_{\mathfrak C}$ this implies:
$$
\mathcal V_{\mathfrak C}(f_\kappa)=\sum_{i=0}^r \frac{1}{N^{r-i}}\frac{\nu_i(g_i)}{b_i} e_{\tau_i}
=\frac{1}{N^{r-j}}\frac{\nu_j(f_\kappa^{N^{r-j}})}{b_j} e_{\kappa}=e_{\kappa}.
$$
If $\kappa\not\in \mathfrak C$, then let $\tau_k\in \mathfrak C$ be the unique element such that $\ell(\tau_k)=\ell(\kappa)$. 
Since $\tau_k$ and $\kappa$ are not comparable with respect to the Bruhat order on $A_\tau$, $f_\kappa$ 
vanishes on $\hat{X}(\tau_k)$. But $f_\kappa$ is not the zero function, so there exist elements 
$\tau_{j}>\tau_{j-1}\ge \tau_k$ in $\mathfrak C$,
such that $\tau_{j}$ covers $\tau_{j-1}$, $f_\kappa \vert_{X(\tau_{j})}\not\equiv 0$, but $f_\kappa$ vanishes on $\hat{X}(\tau_{j-1})$.
It follows with the same arguments as above:
$$
(f_\kappa)_{\mathfrak C}=(f_\kappa,f_\kappa^N,\ldots,f_\kappa^{N^{r-j}},\ldots)\Longrightarrow \mathcal V_{\mathfrak C}(f_\kappa)
=\frac{\nu_{{j}}(f_\kappa^{N^{r-j}})} {N^{r-j}b_j} e_{p_j} +\sum_{i<j} c_ie_{p_i}
$$
for some rational numbers $c_i\in \mathbb Q$, $0\le i\le j-1$. By assumption we have
$\nu_{j}(f_\kappa)>0$ and $\ell(\tau_j)>\ell(\kappa)$, so $\mathcal V_{\mathfrak C}(f_\kappa)>e_\kappa$ with respect to the total order
on $\mathbb Q^{A_\tau}$. Since  $\mathcal V(f_\kappa)$ is the minimum of the $V_{\mathfrak C}(f_\kappa)$,  $\mathfrak C$
running over all maximal chains in $A_\tau$, we get $\mathcal V(f_\kappa)=e_\kappa$.
\end{proof}

\begin{remark}
The proof shows in fact more precisely: the property $\mathcal V(f_\kappa)=e_\kappa$ holds for every choice 
of a total order $\ge_t$ that respects the length.
\end{remark}

Let $\Gamma:=\{\mathcal V(g)\mid g\in\mathbb K[\hat X(\tau)]\setminus\{0\} \}\subseteq \mathbb Q^{A_\tau}$ 
be the image of the quasi-valuation. Since $\mathcal V$ is only a quasi-valuation, 
$\Gamma$ is in general not anymore a monoid as in the case
of valuations. 

But we have the following compensation:
for a maximal chain $\mathfrak{C}\in\mathcal{C}$, we set
\begin{equation}\label{GammaC}
\Gamma_{\mathfrak C}:=\{\underline{a}\in\Gamma\mid \supp\,\underline{a}\subseteq \mathfrak C\}.
\end{equation}
 By Proposition 8.7 in \cite{CFL} we have $\mathcal V(g)=\mathcal V_{\mathfrak C}(g)$ if and only if 
$\supp\,\mathcal V(g)\subseteq \mathfrak C$. So if $g,h\in \mathbb{K}[\hat{X}(\tau)]\setminus\{0\}$ are such that 
$\supp\,\mathcal V(g),\supp\,\mathcal V(h)\subseteq \mathfrak C$, then 
$$
\mathcal V_{\mathfrak C}(g)+\mathcal V_{\mathfrak C}(h)=\mathcal V_{\mathfrak C}(gh)\ge \mathcal V(gh)
\ge \mathcal V(g)+\mathcal V(h)=\mathcal V_{\mathfrak C}(g)+\mathcal V_{\mathfrak C}(h),
$$
which implies equality everywhere and hence:
\begin{lemma} 
$\Gamma_{\mathfrak C}$ is a monoid. 
\end{lemma}
Being a union of monoids, $\Gamma=\bigcup_{\mathfrak C\in\mathcal{C}} \Gamma_{\mathfrak C}$ gets a name:
\begin{definition}\label{sectionfour_fan_of_monoids}
$\Gamma$ is called the \emph{fan of monoids} associated to the quasi-valuation $\mathcal V$. 
\end{definition}
The following is a summary of some of the results in \cite{CFL}, reformulated in the
setting of Schubert varieties.  
Recall that a quasi-valuation defines a filtration of $\mathbb{K}[\hat{X}(\tau)]$ by subspaces in the same way as in \eqref{filtration:1} and \eqref{filtration:2}. 
{ 
For $\underline a\in \mathbb Q^{A_\tau}$ we set:
\begin{equation}\label{filtration:3}
{\mathbb{K}[\hat{X}(\tau)]}_{\ge_t \underline a}:=\{ g\in   \mathbb{K}[\hat{X}(\tau)]\setminus\{0\}\mid  
\mathcal V(g)\ge_t \underline a\}\cup\{0\},
\end{equation}
respectively
\begin{equation}\label{filtration:4}
{\mathbb{K}[\hat{X}(\tau)]}_{>_t \underline a}:=\{ g\in   \mathbb{K}[\hat{X}(\tau)]\setminus\{0\}\mid  
\mathcal V(g)>_t \underline a\}\cup\{0\}.
\end{equation}
}
The subquotients { of these subspaces} are called leaves. 
\begin{theorem}[\cite{CFL}]\label{theorem4}
Let $\mathcal V$ be the quasi-valuation on $\mathbb{K}[\hat X(\tau)]\setminus\{0\}$ associated to the Seshadri stratification
defined in Section~\ref{Seshadri_stratification} and the total order ``$\le_t$'' fixed Section~\ref{higher_valued_valuation}. The following statements hold:
\begin{itemize}
\item[{\it i)}] $\mathcal{V}$ has at most one-dimensional leaves.
\item[{\it ii)}] The fan of monoids $\Gamma$ is contained in $\mathbb{Q}_{\ge 0}^{A_\tau}$.  
\item[{\it iii)}] The fan of monoids $\Gamma$ is the union of the finitely generated monoids $\Gamma_{\mathfrak C}$, where the union is running over all maximal chains
$\mathfrak C$.
\item[{\it iv)}] The quasi-valuation is additive if the supports of both functions are contained in a maximal chain: for $g,h\in \mathbb K[\hat{X}(\tau)]\setminus\{0\}$, $\mathcal V(gh) = \mathcal V(g)+\mathcal V(h)$ if and only if 
there exists a maximal chain $\mathfrak C$ such that $\supp \mathcal V(g),\supp \mathcal V(h) \subseteq \mathfrak C$.
\end{itemize}
\end{theorem}

{ We can  
recover the degree of a homogeneous function from its quasi-valuation.} Notice that 
$\mathrm{deg}(f_\kappa)=1$ for all $\kappa\in A_\tau$.

\begin{proposition}[\cite{CFL}]\label{powerrelation2}
Let $g\in \mathbb K[\hat X(\tau)]\setminus\{0\}$ and suppose $\mathcal V(g)=\sum_{\kappa\in A_\tau}a_\kappa e_\kappa$. Let $m$ be such that $m a_\kappa\in \mathbb N$ for all $\kappa\in A_\tau$. Then there exist $\lambda\in\mathbb{K}^*$ and $g'\in \mathbb K[\hat X(\tau)]$ such that 
$$
g^m=\lambda \prod_{\kappa\in A_\tau} f^{ma_\kappa}_\kappa +g'
$$
with $\mathcal V(g')>\mathcal V(g^m)$ when $g'\neq 0$. If $g$ is homogeneous and $g'\not=0$, 
then $\prod_{\kappa\in A_\tau} f^{ma_\kappa}_\kappa$ and $g'$ are homogeneous of the same degree as $g$.
In particular, if $g$ is homogeneous, then: 
$$
\deg g=\sum_{\kappa\in A_\tau}a_\kappa.
$$
\end{proposition}
Let $C$ be a (not necessarily maximal) chain in $A_\tau$. We define a cone $K_C\subseteq\mathbb{R}^{A_\tau}$ by:
$$
K_C:=\sum_{\kappa\in C}\mathbb{R}_{\geq 0}e_\kappa.
$$
The positivity property of the quasi-valuation implies that the subspaces defined by the filtration are ideals, so we can consider
the associated graded algebra $\mathrm{gr}_{\mathcal V}\mathbb K[\hat X(\tau)]$. To describe this degenerate algebra, we introduce the notion of a 
fan algebra:

\begin{definition}[\cite{CFL}]\label{Defn:FanAlgebra}
The \textit{fan algebra} $\mathbb K[\Gamma]$ associated to the fan of monoids $\Gamma$ is defined as
$$
\mathbb K[\Gamma]:=\mathbb K[x_{\underline{a}}\mid \underline{a}\in \Gamma] / I(\Gamma)
$$
where $I(\Gamma)$ is the ideal generated by the following elements:
$$
\left\{
\begin{array}{rl}
x_{\underline{a}}\cdot x_{\underline{b}}-x_{\underline{a}+\underline{b}},&\textrm{if there exists a chain $C\subseteq A_\tau$ such that }\underline{a},\underline{b}\in K_C;\\
x_{\underline{a}}\cdot x_{\underline{b}},&\textrm{if there exists no such a chain.}\\
\end{array}
\right.
$$
\end{definition}

To simplify the notation, we shall write $x_{\underline{a}}$ also for its class in $\mathbb{K}[\Gamma]$ when there is no ambiguity.
For a maximal chain $\mathfrak C$ denote by $\mathbb K[\Gamma_{\mathfrak C}]$ the following subalgebra:
$$
\mathbb K[\Gamma_{\mathfrak C}]:=\bigoplus_{\underline{a}\in\Gamma_{\mathfrak C}} \mathbb Kx_{\underline{a}}\subseteq \mathbb K[\Gamma],
$$
then $\mathbb K[\Gamma_{\mathfrak C}]$ is naturally isomorphic to the usual semigroup algebra associated to the monoid $\Gamma_{\mathfrak C}$.
We endow the algebra $\mathbb{K}[\Gamma]$ with a grading inspired by Proposition~\ref{powerrelation2}: for $\underline{a}\in\Gamma\subseteq \mathbb Q_{\ge 0}^{A_\tau}$, the degree of $x_{\underline{a}}$ is defined by
$$
\deg x_{\underline{a}}=\sum_{p\in A_\tau} a_p.
$$

One of the main results in \cite{CFL}, when applied to this Seshadri stratification on $X(\tau)$, gives:

\begin{theorem}[\cite{CFL}]\label{degenerationTheorem}
There exists a flat degeneration $\Psi:\mathcal A\rightarrow \mathbb A^1$ such that the generic fibre is isomorphic to $X(\tau)$ and the special fibre $X_0$ is isomorphic
to $\mathrm{Proj\,}(\mathrm{gr}_{\mathcal V} \mathbb{K}[\hat X(\tau)])$. The degenerate variety  $X_0$ is a reduced union of equidimensional projective toric varieties, one irreducible component for each maximal chain in $A_\tau$. The irreducible component associated to a maximal chain $\mathfrak C$ is isomorphic to $\mathrm{Proj\,}(\mathbb K[\Gamma_{\mathfrak C}])$.
\end{theorem}

\section{
{ 
The fan of monoids $\Gamma$ 
and the Lakshmibai-Seshadri lattice}}\label{CandidatForMonoid}
{ 
The fan of monoids $\Gamma$ is a central object 
in the theory of Seshadri stratifications. For example,
by Theorem \ref{degenerationTheorem} the fan algebra completely describes the semi-toric degeneration of the variety $X$.   
So it is important to get a concrete description { of} the monoids $\Gamma_{\mathfrak C}$ for all maximal chains $\mathfrak C$.

We consider the Seshadri stratification of a Schubert variety $X(\tau)$ defined in section~\ref{Seshadri_stratification}.
We will see later (Theorem~\ref{LScandidat_FanOfMonoids}), that the fan of monoids $\Gamma$ associated to the 
quasi-valuation $\mathcal V$ coincides with the fan of monoids $\mathrm{LS}^+_{\lambda}$ we present now. 

The fan of monoids $\mathrm{LS}^+_{\lambda}$ can be described as a reformulation of a conjecture made by Lakshmibai. In
\cite{LS2}, see also Appendix C of \cite{SMT}, she made a conjecture on a possible index system for a basis of the space of global sections $\mathrm{H}^0(X(\tau),\mathcal{L} _\lambda)$. The conjecture of Lakshmibai was reformulated later in terms of the Lakshmibai-Seshadri path model for representations and was proved in \cite{L2}. 

More details about the connection between Lakshmibai's conjecture, the path model, as well as the lattices and monoids
defined below are explained in Appendix I, Section~\ref{appendix1}.

\subsection{The LS-lattice}
We fix a maximal chain }
$\mathfrak C:\ \tau=\tau_r>\ldots>\tau_0=\mathrm{id}$ 
in $A_\tau$.  
For such a fixed maximal chain, we simplify the notation by writing $b_i$ instead of $b_{\tau_i,\tau_{i-1}}$ for the bonds. Note that 
all extremal functions $f_{\tau_k}$ have degree $1$ and hence $b_0=1$.

\begin{definition}
The lattice $\mathrm{LS}_{\mathfrak C,\lambda}\subseteq \mathbb Q^{\mathfrak C}$ defined underneath 
is called the \textit{Lakshmibai-Seshadri lattice} (for short \textit{LS-lattice}) 
associated to $\lambda$ and $\mathfrak C$:\index{LS-lattice}\index{$\mathrm{LS}_{\mathfrak C,\lambda}$, LS-lattice}
\begin{equation}\label{EqLSLattice}
\mathrm{LS}_{\mathfrak C,\lambda}:=\left\{\underline{a}=\left(\begin{array}{c}a_r \\ \vdots \\ a_0\end{array}\right)\in \mathbb Q^{\mathfrak C}\left\vert 
{\scriptsize
\begin{array}{r}
b_{r}a_r\in\mathbb Z\\
b_{r-1}(a_r+a_{r-1})\in\mathbb Z\\
\ldots\\
b_{1}(a_r+a_{r-1}+\ldots+a_1)\in\mathbb Z\\ 
a_0 +a_1+\ldots+a_r \in \mathbb Z\\
\end{array} }
 \right.\right\}.
\end{equation}
\end{definition}
Motivated by Theorem~\ref{theorem4} and Proposition \ref{powerrelation2} we set:
\begin{definition}\label{monoiddefinitionintersection}
The  monoid obtained as the intersection $\mathrm{LS}_{\mathfrak C,\lambda}\cap  \mathbb Q_{\ge 0}^{\mathfrak C}$
of the LS-lattice with the positive quadrant is denoted by $\mathrm{LS}_{\mathfrak C,\lambda}^+$.

For $\underline{a}\in \mathrm{LS}_{\mathfrak C,\lambda}$ we call the sum $a_0 +a_1+\ldots+a_r$ 
the \textit{degree of $\underline{a}$}. For $m\ge 0$ let $\mathrm{LS}_{\mathfrak C,\lambda}^+(m)$ be the subset of elements in $\mathrm{LS}^+_{\mathfrak C,\lambda}$ of degree $m$.
\end{definition}{ 
\subsection{The LS-fan of monoids}
As in the case of the fan of monoids associated to a 
quasi-valuation (Definition~\ref{sectionfour_fan_of_monoids}) we set:}
\begin{definition}
Denote by $\mathrm{LS}_{\lambda}\subseteq \mathbb Q^{A_\tau}$ 
the set theoretic union of the lattices $\mathrm{LS}_{\mathfrak C,\lambda}$ over all maximal chains 
$\mathfrak{C}\in\mathcal{C}$. The intersection $\mathrm{LS}_{\lambda}\cap  \mathbb Q_{\ge 0}^{A_\tau}$
is denoted by $\mathrm{LS}_{\lambda}^+$, which is called the \textit{Lakshmibai-Seshadri fan of monoids}
(for short: the LS-fan of monoids) associated to $\lambda$ and $\tau$.
\end{definition}
The notion of the degree of an element and its support is extended to $\mathrm{LS}_{\lambda}$ in the obvious way. 
Let $\mathrm{LS}_{\lambda}^+(m)$ for $m\ge 0$ be the subset of elements of degree $m$.  
For $\underline{a}=(a_\kappa)_{\kappa\in A_\tau}\in \mathrm{LS}_{\lambda}^+(m)$ we define \index{$\textrm{weight}(\underline{a})$, weight of $\underline{a}$}
$$
\textrm{weight}(\underline{a}):=\sum_{\kappa\in A_\tau} a_\kappa\kappa(\lambda).
$$
To avoid too many indices, we do not add an index $\tau$ to the fan of monoids $\mathrm{LS}_{\lambda}^+$ as long
as it is clear from the context that the fan is associated to the partially ordered set $A_\tau$ with
the corresponding bonds. 

For example, right now we add the index $\mathrm{LS}_{\lambda,\tau}^+$, just to point out: if $\tau'>\tau$, then the inclusion
$A_\tau\hookrightarrow A_{\tau'}$ induces an inclusion $\mathbb Q^{A_\tau}\hookrightarrow \mathbb Q^{A_{\tau'}}$
sending $e_\kappa\in \mathbb Q^{A_\tau}$ for $\kappa\in A_{\tau}$  to $e_\kappa\in\mathbb Q^{A_{\tau'}}$. It follows directly from the definition of the fan of monoids
that this inclusion induces an inclusion of fans of monoids:
\begin{equation}\label{extending_the_fan}
\begin{tikzcd}
\mathbb Q^{A_\tau} \arrow[r, hook]  & \mathbb Q^{A_{\tau'}} \\
\mathrm{LS}_{\lambda,\tau}^+ \arrow[r, hook] \arrow[u, hook] & \mathrm{LS}_{\lambda,\tau'}^+ \arrow[u, hook].
\end{tikzcd}
\end{equation}
We have more precisely: $\mathrm{LS}_{\lambda,\tau}^+=\mathrm{LS}_{\lambda,\tau'}^+\cap \mathbb Q^{A_\tau}$.
This ``extension'' $\mathrm{LS}_{\lambda,\tau}^+\hookrightarrow \mathrm{LS}_{\lambda,\tau'}^+$
 of the fan of monoids makes it possible to use induction arguments.
 \vskip 5pt
\subsection{Panorama on the next sections} Consider the Seshadri stratification defined in Section~\ref{Seshadri_stratification} for the embedded Schubert variety 
$X(\tau)\hookrightarrow \mathbb P(V(\lambda)_\tau)$ and the associated quasi-valuation defined in 
Section~\ref{higher_valued_valuation}. We have associated to a maximal chain in \eqref{GammaC} a monoid 
$\Gamma_{\mathfrak C}\subseteq \mathbb Q^{\mathfrak{C}}$, let  ${\mathcal L}^{\mathfrak{C}}\subseteq \mathbb Q^{\mathfrak{C}}$ 
be the lattice generated by the monoid $\Gamma_{\mathfrak C}$.
\vskip 5pt
The aim of the next sections is to show:{ 
\begin{enumerate}
\item[i)] the lattice ${\mathcal L}^{\mathfrak{C}}$
 coincides with the lattice $\mathrm{LS}_{\mathfrak C,\lambda}$;
\item[ii)] the monoid $\Gamma_{\mathfrak{C}}$ coincides with 
the monoid $\mathrm{LS}_{\mathfrak C,\lambda}^+$;
\item[iii)] the path vectors, which will be described in the next section, define representatives for all the leaves of the 
quasi-valuation $\mathcal V$.
\end{enumerate}
}

\section{Some special vectors and functions}\label{special-functions}
We come back to the situation as before: a Schubert variety $X(\tau)\hookrightarrow\mathbb P(V(\lambda)_\tau)$
embedded in the projective space over a Demazure submodule $V(\lambda)_\tau\subseteq V(\lambda)$.
To describe a filtration on the Demazure modules $V(\lambda)_\tau$, we need the following partial order on $\mathbb Q^{A_\tau}$.

For $\underline{a}\in \mathbb Q^{A_\tau}$ denote by $\{\supp\underline{a}\}_j$ the subset of elements in $\supp\underline{a}$ of length $j$.
We declare $\underline{a}$ is \textit{thin} if $\# \{\supp\underline{a}\}_j\le 1$ for all $j\ge 0$. For example, all the elements
in $\mathrm{LS}_{\lambda}^+$ are thin.

\begin{definition}\label{newpartialorder}\rm
Let $\underline{a}, \underline{b}\in \mathbb Q^{A_\tau}_{\ge 0}$ be both thin and of the same degree  $m$.
Let $\supp\underline{a}=\{\tau_1,\tau_2,\ldots\}$ and $\supp\underline{b}=\{\kappa_1,\kappa_2,\ldots\}$
be enumerated in such a way that $\ell(\tau_1)>\ell(\tau_2)>\ldots$ and $\ell(\kappa_1)>\ell(\kappa_2)>\ldots$.

We say $\underline{a}\rhd \underline{b}$ \index{$\lhd$, partial order on $\mathbb Q^{A_\tau}$}
if $\tau_1>\kappa_1$, or $\tau_1=\kappa_1$ and $a_{\tau_1}>b_{\kappa_1}$,
or $\tau_1=\kappa_1$ and $a_{\tau_1}=b_{\kappa_1}$ and $\tau_2>\kappa_2$, or 
$\tau_1=\kappa_1$ and $a_{\tau_1}=b_{\kappa_1}$ and $\tau_2=\kappa_2$ and $a_{\tau_2}>b_{\kappa_2}$
or $\ldots$.
\end{definition} 
An immediate consequence for the fixed total order $\ge_t$ on $\mathbb Q^{A_\tau}$ is the following:
\begin{lemma}\label{rhd_and-total_order}
If $\underline{a},\underline{b}\in {\textrm{\rm LS}}^+_{\lambda}(s)$, then $\underline{a}\rhd \underline{b}$ implies $\underline{a}\ge_t \underline{b}$
for all possible choices of $\ge_t $. 
\end{lemma}
{ 
Given $\underline a\in \mathrm{LS}_{\lambda}^+(1)$, let $\underline \sigma$ be a reduced decomposition
of the maximal element $\sigma$ in $\supp\underline a$. (Note that for $\underline a\in \mathrm{LS}_{\lambda}^+$ 
the support is linearly ordered with respect to ``$\ge$'', so the maximal element is independent of the choice of $\ge_t$.)

In Section~\ref{algorithm} we describe a purely combinatorial algorithm how to associate to such a pair
$(\underline a,\underline \sigma)$ a sequence of integers $(n_1,\ldots,n_t)$ and a sequence of simple roots
$(\alpha_{i_1},\ldots,\alpha_{i_t})$. The vector $v_{\underline{a},\underline{\sigma}}\in V(\lambda)_\tau$ 
\textit{associated} to $\underline{a}\in \mathrm{LS}_{\lambda}^+(1)$
and the reduced decomposition $\underline{\sigma}$ is defined by:
$$
v_{\underline{a},\underline{\sigma}}=X_{-i_1}^{(n_1)}\cdots X_{-i_t}^{(n_t)} v_{\lambda}.
\index{$v_{\underline{a},\underline{\sigma}}$, vector associated to $\underline{a}$}
$$
For every $\underline{a}\in \mathrm{LS}_{\lambda}^+(1)$ fix a reduced decomposition $\underline{\sigma}^{\underline a}$
of its maximal element in $\sigma$ in $\supp\underline a$. It was shown in \cite{L1} that the collection
of vectors $\{v_{\underline{a},\underline{\sigma}^{\underline a}}\mid \underline{a}\in  \mathrm{LS}_{\lambda}^+(1)\}$
is a basis of the Demazure module $V(\lambda)_\tau$.
The basis was actually constructed over $\mathbb Z$, so that it specializes to a basis
for any algebraically closed field. For more details see \textit{ibidem}, respectively Appendix II, Section~\ref{A_filtration_on_V}. 

The basis is far from being canonical. 
What is canonical about this construction is the { following collection of subspaces} defined by these vectors.
For more details see the Appendix II, Section~\ref{A_filtration_on_V}. }%
For $\underline{a}\in \mathrm{LS}_{\lambda}^+(1)$ set
\begin{equation}\label{section_six_filt_one}
V(\lambda)_{\tau,\unlhd\underline a}=\left\langle 
v_{\underline a',\underline \sigma'}\,\left\vert\, \underline a\unrhd \underline a';
\begin{array}{l}
\underline a'\in \mathrm{LS}_{\lambda}^+(1), \underline\sigma' \textit{\ reduced decomposition}\\
\textit{of maximal element in\ }\supp \underline a'
\end{array}\right.
\right\rangle_{\mathbb K},
\end{equation}
and let 
\begin{equation}\label{section_six_filt_two}
V_{}(\lambda)_{\tau,\lhd\underline{a}}=\left\langle 
v_{\underline a',\underline{\sigma}'}\,\left\vert\,  \underline a\rhd \underline a';
\begin{array}{l}
\underline a'\in \mathrm{LS}_{\lambda}^+(1), \underline{\sigma}' \textit{\ reduced decomposition}\\
\textit{of maximal element in\ }\supp \underline a'
\end{array}\right.
\right\rangle_{\mathbb K}.
\end{equation}
The following is known about these { subspaces}, see Appendix II, Section~\ref{A_filtration_on_V} for details:

\begin{proposition}\label{summary1}
\begin{enumerate}
\item[\it i)] For $\underline{a}\in \mathrm{LS}_{\lambda}^+(1)$, the subspaces $V(\lambda)_{\tau,\unlhd\underline{a}}$ of $V(\lambda)_\tau$ are
$U(\mathfrak g)_{\mathbb K}^+$-stable.
\item[\it ii)] For $\underline{a}\in \mathrm{LS}_{\lambda}^+(1)$, the leaves  $V(\lambda)_{\tau,\unlhd\underline{a}}/V(\lambda)_{\tau,\lhd\underline{a}}$ are one dimensional.
\item[\it iii)] For any choice of a reduced decomposition $\underline{\sigma}$ of the maximal element 
$\sigma$ in the support $\supp\underline{a}$ of $\underline{a}\in \mathrm{LS}_{\lambda}^+(1)$, 
$v_{\underline{a},\underline{\sigma}}$ is a representative
of the leaf $V(\lambda)_{\tau,\unlhd\underline{a}}/V(\lambda)_{\tau,\lhd\underline{a}}$, and its
class $\bar v_{\underline{a},\underline{\sigma}}$ is a generator of the leaf.
\item[\it iv)] If $\underline{\sigma}'$ and $\underline{\sigma}$ are different reduced decompositions 
of the maximal element 
$\sigma$ in the support $\supp\underline{a}$ of $\underline{a}\in \mathrm{LS}_{\lambda}^+(1)$, then 
$v_{\underline{a},\underline{\sigma}'}=v_{\underline{a},\underline{\sigma}}+
\sum b_{\underline{a}''}v_{\underline{a}'',\underline{\sigma}''}$, where the sum is running over 
elements $\underline{a}''\in \mathrm{LS}_{\lambda}^+(1)$ such that $\underline{a}\rhd \underline{a}''$.
\end{enumerate}
\end{proposition}

So it makes sense to write just $v_{\underline{a}}$ instead
of $v_{\underline{a},\underline{\sigma}}$ if no confusion is possible. 
{  For later purpose it is important to note that neither the definition of the fan of monoids $\mathrm{LS}_{\lambda}^+$
nor the definition of the vectors $v_{\underline{a},\underline{\sigma}}$ respectively of the subspaces $V(\lambda)_{\tau,\unlhd\underline a}$
and $V(\lambda)_{\tau,\lhd\underline a}$ involves or depends on the choice of the total order $\ge_t$.}
We will now define a kind of dual basis which again induces a filtration, but this time on the dual space. {  Note that the
following definition of a path vector is more general than the one in \cite{L1} and does not use the quantum Frobenius splitting.}
\begin{definition}\label{define_path_vector}
Let $\underline{a}\in \mathrm{LS}_{\lambda}^+(1)$.  A \textit{path vector associated to $\underline a$} is 
a linear function $p_{\underline a}\in  V(\lambda)^*_{\tau}$ which is a 
$T$-eigenvector of weight $(-\textrm{weight}(\underline a))$, and such that 
\begin{enumerate}
\item[\it i)] {  there exists a reduced decomposition $\underline{\sigma}$ of the maximal element $\sigma$ in $\mathrm{supp}\,\underline{a}$ such that 
$p_{\underline a}(v_{\underline a,\underline \sigma})=1$};
\item[\it ii)] 
for $\underline{a'}\in\mathrm{LS}_\lambda^+(1)$ and for some reduced
decomposition $\underline{\sigma}'$ of the maximal element $\sigma'$ 
in $\mathrm{supp}\,\underline{a}'$,  
$p_{\underline a} (v_{\underline a',\underline \sigma'})\not=0$
implies $\underline a'\unrhd \underline a$.
\end{enumerate}
\end{definition}
{ 
\begin{lemma}\label{indipendentproofpathvector}
The definition of a path vector is independent of the choice of the reduced decompositions of 
$\sigma$ and  $\sigma'$. 
\end{lemma} 
\begin{proof}
We have to show that if the properties in the definition hold for one choice of a reduced decomposition,
then it holds for any choice of a reduced decomposition.

We start with property \textit{ii)} and prove it by induction over the partial order $\unrhd$.
The set $\mathrm{LS}_{\lambda}^+(1)$ has a unique minimal element with respect to $\unrhd$, it is the element
$\underline{a}''=(0,\ldots,0,1)$. It can also be characterized as the unique element of weight $\lambda$, and hence
$p_{\underline a} (v_{\underline{a}'',\underline{id}})=0$ unless $ \underline a=\underline{a}''$.

Suppose now $\underline{a}''\in \mathrm{LS}_{\lambda}^+(1)$ is given and $\underline a''\not\hskip -3pt\unrhd \underline a$.
Note if $\underline{a}'''\in  \mathrm{LS}_{\lambda}^+(1)$ is such that $\underline{a}''\rhd \underline{a}'''$, then 
$\underline{a}'''\not\hskip -3pt\unrhd \underline a$ and we may assume by induction: 
$p_{\underline a}(v_{\underline{a}''',\underline{\sigma}'''})=0$ for any reduced decomposition $\underline{\sigma}'''$
of the maximal element ${\sigma}'''$ in the support of  $\underline{a}'''$. By the definition of a path vector, we know there
exists at least one reduced decomposition $\underline{\sigma}''$ of the maximal element in the support of $\underline{a}''$
such that $p_{\underline a}(v_{\underline{a}'',\underline{\sigma}''})=0$. Now Proposition~\ref{summary1}, part \textit{iv)},
and induction implies: $p_{\underline a}(v_{\underline{a}'',\underline{\underline{\sigma}}''})=0$ for any reduced decomposition
$\underline{\underline{\sigma}}''$ of $\sigma''$. 

To prove part \textit{i)}, let $p_{\underline a}$ be a path vector associated to $\underline{a}\in \mathrm{LS}_{\lambda}^+(1)$.
By Proposition~\ref{summary1}, part {\textit{iv)}},  and part \textit{ii)} of Lemma~\ref{indipendentproofpathvector}, 
$p_{\underline a} (v_{\underline a,\underline \sigma})=1$
for any reduced decomposition $\underline \sigma$ of the maximal element $\sigma$ in $\mathrm{supp}\,\underline{a}$.
\end{proof}
}

The following lemma is helpful when using induction arguments. Notations in \eqref{extending_the_fan} will be adopted without mention.
{  Since now and in the following lemma we have two Schubert varieties: $X(\tau)$ and $X(\sigma)$, it makes sense to add the index $\tau$ respectively $\sigma$  to the monoids: $\mathrm{LS}_{\lambda,\tau}^+$ respectively $\mathrm{LS}_{\lambda,\sigma}^+$, to avoid confusion.
Note that $\supp\underline{a}$ for $\underline{a}\in \mathrm{LS}^+_{\lambda,\tau}$
is contained in a maximal chain $\mathfrak C\subseteq A_\tau$, 
hence it is linearly ordered and the maximal element is independent 
of the choice of ``$\ge_t$''.

Recall that if $\underline{a}$ is such that the maximal element in $\mathrm{supp}\,\underline{a}$ is smaller or equal to $\sigma$,
then by \eqref{extending_the_fan}, $\underline{a}$ can be naturally viewed as an element in $\mathrm{LS}^+_{\lambda,\sigma}(1)$.

\begin{lemma}\label{inductionOneLemma}
Suppose $\sigma<\tau$, and for $\underline{a}\in\mathrm{LS}^+_{\lambda,\tau}(1)$ and let $p_{\underline{a}}\in  V(\lambda)^*_{\tau}$
be a path vector associated to $\underline{a}$.
If $\underline{a}$ is such that the maximal element in $\mathrm{supp}\,\underline{a}$ is smaller or equal to $\sigma$,
then the restriction $p_{\underline{a}}\vert_{V(\lambda)_{\sigma}}$ is a path vector in $V(\lambda)^*_{\sigma}$ associated to $\underline{a}\in\mathrm{LS}^+_{\lambda,\sigma}(1)$.
\end{lemma}
}
\begin{proof}
{  For all 
$\underline{a}'\in \mathrm{LS}^+_{\lambda,\sigma}(1)$ 
fix a reduced decomposition $\underline{\sigma}'$
of the maximal element $\sigma'$ in $\supp{\underline{a}'}$}.
By Theorem~\ref{AuxBasisTheorem}, the
set $\{v_{\underline{a}',\underline{\sigma}'}\mid\underline{a}'\in \mathrm{LS}^+_{\lambda,\sigma}(1)\}$ forms a basis of the Demazure module
$V(\lambda)_{\sigma}$. It follows by the definition of a path vector that the restriction
$p_{\underline{a}}\vert_{V(\lambda)_{\sigma}}$ is still a weight vector of weight $(-\textrm{weight}(\underline a))$. This restricted function satisfies $p_{\underline{a}}\vert_{V(\lambda)_{\sigma}}(v_{\underline{a},\underline{\sigma}})=1$, and for all $\underline{a}'\in \mathrm{LS}^+_{\lambda,\sigma}(1)$
we have: $p_{\underline{a}}\vert_{V(\lambda)_{\sigma}}(v_{\underline a',\underline \sigma'})\not=0$
only if $\underline a'\rhd \underline a$. Hence  $p_{\underline{a}}$, restricted to $V(\lambda)_\sigma$, is a path vector in $V(\lambda)_\sigma^*$
associated to $\underline{a}\in \mathrm{LS}^+_{\lambda,\sigma}(1)$.
\end{proof}

The path vectors have another special property which helps to apply induction arguments. 
Again suppose $\sigma<\tau$. We get an induced Seshadri stratification for $X(\sigma)$, which has as associated 
partially ordered set the subset $A_\sigma\subseteq A_\tau$ (Example~\ref{induction}). To define the quasi-valuation $\mathcal V$, we have
fixed a total order $\ge_t$ on $A_\tau$ (satisfying the conditions in Section~\ref{higher_valued_valuation}). We denote by 
the same symbol $\ge_t$ the induced 
total order on  $A_\sigma$. Let $\mathcal V_\sigma$ be the corresponding quasi-valuation on the homogeneous
coordinate ring $\mathbb K[\hat X(\sigma)]$
given by the embedding $X(\sigma)\hookrightarrow X(\tau)\hookrightarrow \mathbb P(V(\lambda)_\tau)$.
Note that the image of $X(\sigma)$ is contained in $\mathbb P(V(\lambda)_\sigma)$.
{ 
\begin{lemma}\label{inductionTwoLemma}
Let $\underline{a}\in\mathrm{LS}^+_{\lambda,\tau}(1)$ be such that $\sigma$ is the largest  elements in $\supp\underline{a}$,
so $\underline{a}$ can be naturally viewed as an element in $\mathrm{LS}^+_{\lambda,\sigma}(1)$. For a path vector $p_{\underline{a}}$
associated to $\underline{a}$ the quasi-valuation $\mathcal V(p_{\underline{a}})\in\mathbb Q^{A_\tau}$ is equal to the extension by 
zeros of $\mathcal V_\sigma(p_{\underline{a}}\vert_{\hat{X}(\sigma)})\in\mathbb Q^{A_\sigma}$.
\end{lemma}}
%
\begin{proof}
The Demazure module $V(\lambda)_\sigma$ is the linear span of $\hat X(\sigma)$. Let $\underline{\sigma}$ 
be a reduced decomposition of $\sigma$. By assumption, $v_{\underline{a},\underline{\sigma}}\in V(\lambda)_\sigma$, so 
$p_{\underline{a}}\vert_{V(\lambda)_\sigma}$ does not vanish identically, which implies 
$p_{\underline{a}}\vert_{\hat X(\sigma)}$ does not vanish identically.

A maximal chain $\frak C'$ in $A_\sigma$ can always be extended to a maximal chain $\mathfrak C$ in $A_\tau$. From the definition of the valuation in Section \ref{Hrv}, $\mathcal V_{\frak C'}(p_{\underline{a}}\vert_{\hat X(\sigma)})=\mathcal V_{\frak C}(p_{\underline{a}})$
for such a maximal chain.

Suppose now $\frak C$ in $A_\tau$ is a maximal chain that does not contain $\sigma$. Let $\sigma'\in  \frak C$ be the unique element
in $\frak C$ having the same length as $\sigma$. The linear span of $\hat X(\sigma')$ is the Demazure module $V(\lambda)_{\sigma'}$,
which is spanned by vectors $v_{\underline{a}'',\underline{\sigma}''}$, where $\underline{a}''\in \mathrm{LS}^+_{\lambda,\sigma'}(1)\subseteq { \mathrm{LS}^+_{\lambda,\tau}(1)}$.
Again, the inclusion is given by extension by zeros.

By the definition of a path vector we know $p_{\underline a}(v_{\underline a'',\underline \sigma''})\not=0$
only if $\underline a''\unrhd \underline a$. But this would imply that the maximal element in the support of $\underline a''$
(which is smaller or equal to $\sigma'$) would be larger or equal to $\sigma$, which is impossible. It follows: 
$p_{\underline{a}}\vert_{\hat X(\sigma')}\equiv 0$. Since $p_{\underline{a}}$ is not the zero function on $\hat{X}(\tau)$, there exists
elements $\eta>\zeta\ge \sigma'$ in the maximal chain $\frak C$ such that $p_{\underline{a}}$ does not vanish identically on
$\hat X(\eta)$, but vanishes on the divisor $\hat X(\zeta)$. It follows that the coefficient of $e_\eta$ in $\mathcal V_{\frak C}(p_{\underline{a}})$
is strictly positive.  Hence for any choice of the total order $\ge_t$ (satisfying the rules in Section~\ref{higher_valued_valuation}),
$\mathcal V_{\frak C}(p_{\underline{a}})\ge_t \mathcal V_{\frak C''}(p_{\underline{a}})$, where $\frak C''$
is any maximal chain in $A_\tau$ containing $\sigma$. In the following we consider the inclusion $\mathbb Q^{A_\sigma}\hookrightarrow \mathbb Q^{A_\tau} $
as given by extension by zeros.  It follows:
$$
\begin{array}{rcl}
\mathcal V_{\sigma}(p_{\underline{a}}\vert_{X(\sigma)})
&=&\min\left\{ \mathcal V_{\frak C'}(p_{\underline{a}}\vert_{X(\sigma)})\mid \frak C'\textrm{\ maximal chain in $A_\sigma$}\right\}\\
&=&\min\left\{ \mathcal V_{\frak C}(p_{\underline{a}})\mid \frak C\textrm{\ maximal chain in $A_\tau$}, \sigma\in\frak C\right\}\\
&=&\min\left\{ \mathcal V_{\frak C}(p_{\underline{a}})\mid \frak C\textrm{\ maximal chain in $A_\tau$}\right\}
=\mathcal V(p_{\underline{a}}).
\end{array}
$$
\end{proof}
{  We drop now the index $\tau$ again and use  the same  notation as before: ${\mathrm{LS}^+_{\lambda}}$ instead of ${\mathrm{LS}^+_{\lambda,\tau}}$.}
In connection with the Seshadri stratification, the  following ``dual filtration''  
is important. For $\underline{a}\in {\mathrm{LS}^+_{\lambda}}(1)$  
denote by $V(\lambda)_{\tau,\unrhd\underline a}^*$ the  subspace defined as follows
{ 
\begin{equation}\label{section_six_dual_filtr_one}
V(\lambda)_{\tau,\unrhd\underline a}^*
=\left\langle 
p_{\underline{a}'} \bigg\vert p_{\underline{a}'}\textrm{\ a path vector,\ } \underline{a}'\unrhd\underline{a}, 
\underline{a'}\in \mathrm{LS}_{\lambda}^+(1)
\right\rangle_{\mathbb K},
\end{equation}
}
and set 
{ 
\begin{equation}\label{section_six_dual_filtr_two}
V(\lambda)_{\tau,\rhd\underline a}^*=\left\langle 
p_{\underline{a}'} \bigg\vert p_{\underline{a}'}\textrm{\ a path vector,\ } \underline{a}'\rhd\underline{a}, 
\underline{a'}\in \mathrm{LS}_{\lambda}^+(1)\right\rangle_{\mathbb K}.
\end{equation}}
The following is known about these subspaces, see Appendix III, Section~\ref{The_path_vectors}:
\begin{proposition}\label{dualfiltration}
\begin{enumerate}
\item[\it i)] For every $\underline{a}\in {\mathrm{LS}^+_{\lambda}}(1)$ fix a path vector $p_{\underline{a}}$ associated to $\underline{a}$.
The set $\mathbb B(V(\lambda)_\tau)=\{p_{\underline{a}}\mid \underline{a}\in {\mathrm{LS}^+_{\lambda}}(1)\}$ is a basis 
for  $V(\lambda)_{\tau}^*$.
\item[\it ii)] For $\underline{a}\in \mathrm{LS}_{\lambda}^+(1)$,
the leaves  $V(\lambda)_{\tau,\unrhd\underline a}^*/V(\lambda)_{\tau,\rhd\underline a}^*$ are one dimensional,
and any path vector $p_{\underline{a}}$ is a representative of such a leaf.
\end{enumerate}
\end{proposition}
The following is the most important result needed to start a standard monomial theory. 
{ 
\begin{theorem}\label{some_vanishingA}
If $p_{\underline{a}}\in V(\lambda)_{\tau}^*
\simeq \mathbb K[\hat X(\tau)]_1$ is a path vector associated to $\underline{a}\in \mathrm{LS}^+_{\lambda}(1)$,
then $\mathcal V(p_{\underline{a}})=\underline{a}$, independent of the choice of the total order $\ge_t$.
\end{theorem}
As an immediate consequence we get the following relationship
between the collection of subspaces defined in \eqref{section_six_dual_filtr_one}
and \eqref{section_six_dual_filtr_two} and the filtration of $\mathbb{K}[\hat{X}(\tau)]$ induced by the quasi-valuation, see \eqref{filtration:3} and \eqref{filtration:4}.
Recall that the latter has at most one-dimensional leaves (Theorem~\ref{theorem4}).
\begin{coro}
For all $\underline{a}\in \mathrm{LS}^+_{\lambda}(1)$,
the subspaces $V(\lambda)_{\tau,\unrhd\underline a}^*$ respectively $V(\lambda)_{\tau,\rhd\underline a}^*$
are compatible with the filtration induced by the quasi-valuation. { That is to say}, independently of the choice of the total order $\ge_t$, we have
$$
V(\lambda)_{\tau,\unrhd\underline a}^*\subseteq {\mathbb{K}[\hat{X}(\tau)]}_{\ge_t \underline a},\quad
V(\lambda)_{\tau,\rhd\underline a}^*\subseteq {\mathbb{K}[\hat{X}(\tau)]}_{>_t \underline a},
$$ 
and 
$V(\lambda)_{\tau,\unrhd\underline a}^*/V(\lambda)_{\tau,\rhd\underline a}^* \simeq 
{\mathbb{K}[\hat{X}(\tau)]}_{\ge_t \underline a}/{\mathbb{K}[\hat{X}(\tau)]}_{>_t \underline a}$.
\end{coro}
\begin{proof}
The proof is by induction of $\ell(\tau)$,  the case $\ell(\tau)=0$ being obvious. So suppose now $\ell(\tau)\ge 1$,
and the theorem holds for all smaller Schubert varieties. 

To be more precise, if $\sigma<\tau$, then we have an induced Seshadri stratification for the Schubert variety $X(\sigma)$, with associated 
partially ordered set $A_\sigma\subseteq A_\tau$, the fixed total order $\ge_t$ on  $A_\tau$
induces a total order $\ge_t$ on  $A_\sigma$, we get an induced quasi-valuation $\mathcal V_\sigma$ on $\mathbb K[\hat X(\sigma)]$,
and an associated LS-fan of monoids $\mathrm{LS}^+_{\lambda,\sigma}$.

We view this LS-fan of monoids $\mathrm{LS}^+_{\lambda,\sigma}$ as being embedded into $\mathrm{LS}^+_{\lambda}$
via extension by zeros, see \eqref{extending_the_fan}. For $\underline{a}\in \mathrm{LS}^+_{\lambda}(1)$
let $p_{\underline{a}}$ be a corresponding path vector. If the largest element in $\supp\underline{a}$ is $\sigma$ and
$\sigma<\tau$, then $\underline{a}\in \mathrm{LS}^+_{\lambda,\sigma}(1)$. We know by Lemma~\ref{inductionOneLemma} that the
restriction  $p_{\underline{a}}\vert_{\hat X(\sigma)}$ is a path vector for the Seshadri stratification on $X(\sigma)$.

And, by Lemma~\ref{inductionTwoLemma}, we know $\mathcal V(p_{\underline{a}})=\mathcal V_\sigma(p_{\underline{a}}\vert_{\hat X(\sigma)})$, where the
equality has to be read as: $\mathcal V_\sigma(p_{\underline{a}}\vert_{\hat X(\sigma)})$ is treated as an element in $\mathbb Q^{A_\tau}$ via extension by zeros.
So our induction assumption gives: $\mathcal V_\sigma(p_{\underline{a}}\vert_{\hat X(\sigma)})=\underline{a}$, independent of the choice
of the total order on $A_\sigma$.

Now in the proof of Lemma~\ref{inductionTwoLemma} we have seen if $\mathfrak C$ is a maximal chain in $A_\tau$ such that $\sigma\in\mathfrak C$,
then $\mathcal V_{\frak C}(p_{\underline{a}})=\mathcal V_{\frak C'}(p_{\underline{a}}\vert_{\hat{X}(\sigma)})$, where $\frak C'$ is obtained
from $\frak C$ by omitting all elements larger than $\sigma$. And if $\frak C$ is a maximal chain in $A_\tau$ not containing 
$\sigma$, then we have seen in the proof of Lemma~\ref{inductionTwoLemma} that, independent of the choice of the total order $\ge_t$ (satisfying the rules in Section~\ref{higher_valued_valuation}),
$\mathcal V_{\frak C}(p_{\underline{a}})\ge_t \mathcal V_{\frak C''}(p_{\underline{a}})$, where $\frak C''$
is any maximal chain in $A_\tau$ containing $\sigma$.
It follows: $\mathcal V(p_{\underline{a}})=\underline{a}$, independent of the choice of the total order $\ge_t$.

It remains to consider the case: the largest element in $\supp\underline{a}$ is $\tau$. In this case we use 
Corollary~\ref{factoring}. Let $m\ge 1$ be such that $m a_\tau\in\mathbb N$. Then we know that $p_{\underline{a}}^m$
is, up to multiplication by a root of unity, equal to $p_\tau^{ma_\tau} p_{\underline{b}}$. Here
$\underline b=m\underline{a}-ma_\tau e_\tau\in \mathrm{LS}_{\mathfrak C,\lam}^+(m-ma_\tau)$ and $p_{\underline b}\in V((m-ma_\tau)\lambda)^*_\tau$ is a path vector
associated to the leaf $\underline b$. 

By Theorem~\ref{theorem4}, we know the quasi-valuation is additive if the support of the functions is contained in a common maximal chain.
Now by induction we know: $\mathcal V(p_{\underline b})=\underline b$, and this holds independent of the choice of the total order.
By Lemma~\ref{extremalvaluationII} we know: $\mathcal V(p_\tau^{ma_\tau} )=ma_\tau e_\tau$, and this holds independent of the choice of the total order
on $A_{\tau}$. Since $\supp \underline{a}=\supp \underline{b}\cup \{e_\tau\}$ is contained in a maximal chain, we see:
$$
\mathcal V(p_{\underline{a}})=\frac{1}{m}\mathcal V(p^m_{\underline{a}})=\frac{1}{m}\mathcal V(p_\tau^{ma_\tau} p_{\underline{b}})
=\frac{1}{m}\big(\mathcal V(p_\tau^{ma_\tau}) +\mathcal V(p_{\underline{b}})\big)
=a_\tau e_\tau + (\underline{a}-a_\tau e_\tau)=\underline{a},
$$
and this holds independent of the choice of the total order.
\end{proof}
}

{ 
\section{Standard monomial theory}\label{standard-monomial-theory}

\subsection{Standard monomial basis}\label{sectionname_standard_monomial}
Recall that the elements in the support of $\underline{a}\in \mathrm{LS}^+_{\lambda}(s)$
are  linearly ordered with respect to the partial order ``$\ge$'' on $A_\tau$, so 
there exists always a unique maximal element in the support: $\max\supp \underline{a}$,
and a unique minimal element in the support: $\min\supp \underline{a}$.

\begin{definition}\label{def:standard:monomial}
A monomial $p_{\underline{a}^1}\cdots p_{\underline{a}^m}\in \mathbb K[\hat X(\tau)]$ of path vectors with $\underline{a}^1,\ldots,\underline{a}^m\in
\mathrm{LS}_\lambda^+(1)$
is called \emph{standard} if for each $1\leq j\leq m-1$ we have $\min\supp \underline{a}^j \geq \max\supp \underline{a}^{j+1}$.
\end{definition}
An important property of the elements of the LS-fan of monoids is stated in the Appendix I in Lemma~\ref{latticedecomp} (see also \cite{Ch}):
Every element $\underline a\in  \mathrm{LS}^+_{\lambda}(m)$
has a unique decomposition $\underline a=\underline a^1+\ldots+\underline a^m$ into $m$ elements
$\underline a^i\in \mathrm{LS}^+_{\lambda}(1)$, $i=1,\ldots,m$, such that 
$\supp \underline a^1\ge \supp \underline a^2\ge \ldots\ge \supp \underline a^m$.
 
For every $\underline a\in  \mathrm{LS}^+_{\lambda}(1)$ fix a path vector $p_{\underline{a}}$, so this collection 
$\mathbb B(V(\lambda)_\tau)$ of path vectors
forms a basis for $V(\lambda)_\tau^*$ (see Proposition~\ref{dualfiltration}).
 
So we start with the fixed basis $\mathbb B(V(\lambda)_\tau)$ consisting of path vectors.
Given $\underline a\in  \mathrm{LS}^+_{\lambda}(m)$, $m\ge 1$, we have a unique
decomposition $\underline a=\underline a^1+\ldots+\underline a^m$ with
$\underline a^i\in \mathrm{LS}^+_{\lambda}(1)$, $i=1,\ldots,m$.
We associate to $\underline a$ the standard monomial
$$
p_{\underline a}=p_{\underline a^1}\cdots p_{\underline a^m}.
$$
By Theorem~\ref{theorem4}, we know the quasi-valuation is additive if the support of the functions is contained in a common maximal chain.
By Theorem~\ref{some_vanishingA} we know $\mathcal V(p_{\underline{a}^j})=\underline{a}^j$, independent of the choice of the total order $\ge_t$.
As a consequence we see:
$$
\mathcal V(p_{\underline{a}})=\mathcal V(p_{\underline{a}^1})+\ldots+\mathcal V(p_{\underline{a}^m})=\underline a^1+\ldots+\underline a^m=\underline a,
$$
independent of the choice of the total order $\ge_t$,
and hence:
\begin{coro}\label{SMTvaluationIndependent}
For the standard monomial $p_{\underline a}$, we have $\mathcal V(p_{\underline a})=\underline a$, independent of the choice of the total order $\ge_t$.
\end{coro}
This leads us directly to a vector space basis of $\mathbb K[\hat X(\tau)]$ by standard monomials:
\begin{theorem}\label{proposition_standard_monomial_basis1}
Let $\mathbb B(V(\lambda)_\tau)$ be a fixed basis of $V(\lambda)_\tau^*=\mathbb K[\hat X(\tau)]_1$ consisting of path vectors.
The set of standard monomials in the elements of $\mathbb B(V(\lambda)_\tau)$ forms a vector space basis for $\mathbb K[\hat X(\tau)]$. 
\end{theorem}
\begin{proof}
For ${\underline a}\in \mathrm{LS}^+_{\lambda}(m)$ we have $\mathcal V(p_{\underline{a}})=\underline a$ by Corollary~\ref{SMTvaluationIndependent},
which implies that the set of standard monomials $\{p_{\underline{a}}\mid \underline{a}\in \mathrm{LS}^+_{\lambda}(m)\}$ is a set of linearly independent vectors.
On the other hand, by Lemma~\ref{CoordDemazure}, the dimension of $\mathbb K[\hat X(\tau)]_m$ is bounded above by the dimension of $V(m\lambda)_\tau$,
and by Corollary~\ref{DemazureIso} we have $\dim V(m\lambda)_\tau=\mathbb K[\hat X(\tau)]_m$. And one can either use Demzure's character formula and its 
combinatorial version in \cite{L2} (see also Theorem~\ref{path-character} in Appendix I), or one can use \cite{L1} (see also Appendix II, summary before Theorem~\ref{AuxBasisTheorem}) to see that 
$\dim V(m\lambda)_\tau$ is equal to the cardinality of  $\mathrm{LS}^+_{\lambda}(m)$. 
This implies that $\{p_{\underline{a}}\mid \underline{a}\in \mathrm{LS}^+_{\lambda}(m)\}$
is a basis for $\mathbb K[\hat X(\tau)]_m$.
\end{proof}

In \cite{CFL} we have introduced two special properties of a Seshadri stratification. 
Recall that in the process of associating a quasi-valuation to a Seshadri stratification, the only choice we made is a total order $\leq^t$ on $A$ refining the given partial order and preserving the length function. 

\begin{definition}\label{definition_balanced_statification}
A Seshadri stratification of the embedded projective variety $X\hookrightarrow \mathbb P(V)$, with homogeneous coordinate ring $R$, 
is called \emph{balanced} if the following two properties hold:
\begin{enumerate}
\item[\it i)] the fan of monoids $\Gamma$ associated to the quasi-valuation ${\mathcal V}$ is independent of the choice of the total order $\geq_t$;
\item[\it ii)] for each $\underline{a}\in\Gamma$ there exists a regular function $x_{\underline{a}}\in R$ such that 
${\mathcal V}(x_{\underline{a}}) = \underline{a}$ for all possible choices of a total order $\geq_t$.
\end{enumerate}
\end{definition}

The second special property introduced in \cite{CFL} is related to the question whether the monoids $\Gamma_{\mathfrak{C}}$ are saturated, that is to say: if $K_{\mathfrak C}$ is the real cone generated by $\Gamma_{\mathfrak{C}}$
and $\mathcal L^{\mathfrak C}$ is the lattice generated by $\Gamma_{\mathfrak{C}}$, then 
$\mathcal L^{\mathfrak C} \cap  K_{\mathfrak C} =\Gamma_{\mathfrak C}$ for all maximal chains $\mathfrak{C}$. 

Algebraically this is equivalent to say that the algebras $\mathbb{K}[\Gamma_{\mathfrak{C}}]$ (Section~\ref{higher_valued_valuation}) are normal for all maximal chains $\mathfrak{C}$.
Geometrically this condition is equivalent to the normality of all the irreducible components of the degenerate variety
$\mathrm{Spec}(\mathrm{gr}_{\mathcal V}\mathbb{K}[\hat{X}(\tau)])$ (see Theorem~\ref{degenerationTheorem}). We give a name to such a situation:

\begin{definition}\label{def:stratification:normal}
A Seshadri stratification is called \emph{normal} if for all maximal chains $\mathfrak C$ the monoid
$\Gamma_\mathfrak{C}$ is saturated.
\end{definition}

\begin{remark}
In the above definition, we have chosen implicitly a total order $\geq^t$ refining the given partial order on $A$. Different such choices would produce different monoids $\Gamma_\mathfrak{C}$. When the Seshadri stratification is balanced, the monoids $\Gamma_\mathfrak{C}$ are independent of the choice of the refinement of the partial order. In this case, being normal is indeed a property associated to the stratification.
\end{remark}

\begin{theorem}\label{LScandidat_FanOfMonoids}
The Seshadri stratification of the Schubert variety $X(\tau)\hookrightarrow \mathbb P(V(\lambda))$ defined in Section~\ref{Seshadri_stratification} is normal and balanced.
The fan of monoids associated to the quasi-valuation $\mathcal V$ coincides
with the fan of monoids $\mathrm{LS}^+_{\lambda}$. 
\end{theorem}

\begin{proof}
Fix a basis $\mathbb B(V(\lambda)_\tau)$ of $V(\lambda)_\tau ^*=\mathbb K[\hat X(\tau)]_1$ consisting of path vectors.
By Theorem~\ref{proposition_standard_monomial_basis1} we know that 
$\mathbb K[\hat X(\tau)]$ has a vector space basis consisting of the standard monomials in the elements
of $\mathbb B(V(\lambda)_\tau)$. This vector space basis of $\mathbb K[\hat X(\tau)]$ is, by construction (and Lemma~\ref{latticedecomp}),
indexed by the elements of $\mathrm{LS}^+_{\lambda}$, and Theorem~\ref{proposition_standard_monomial_basis1} 
states: for $\underline a\in \mathrm{LS}^+_{\lambda}$ we have  $\mathcal V(p_{\underline a})=\underline a$, independent of the choice of the total order $\ge_t$.

This property has several consequences. Since the quasi-valuation differs for pairwise different basis elements, 
Lemma~\ref{simpleproperties} implies $\mathcal V(p)\in \mathrm{LS}^+_{\lambda}$ for any nonzero element in $\mathbb K[\hat X(\tau)]$.
In particular, the fan of monoids associated to the quasi-valuation $\mathcal V$ coincides
with the fan of monoids $\mathrm{LS}^+_{\lambda}$, and this holds independent of the choice of the total order ``$\ge_t$''.
Moreover, the basis given by the standard monomials has the properties desired in Definition~\ref{definition_balanced_statification},
so the Seshadri stratification is balanced.

For each maximal chain, the monoid $\mathrm{LS}^+_{\mathfrak C,\lambda}$ is defined (see Definition~\ref{monoiddefinitionintersection}) as the intersection of the lattice 
$\mathrm{LS}_{\mathfrak C,\lambda}$ with the positive quadrant $\mathbb Q^{\mathfrak C}_{\ge 0}$, in particular,
the monoid is saturated. The stratification is hence normal.
\end{proof}
}

{  
\section{Some applications}\label{standard-monomial-theory2}
The fact that the fan of monoids coincides with $\mathrm{LS}_\lambda^+$ implies that the stratification is of LS-type
in the terminology of  \cite{CFL3}. Hence Theorem 3.2 and Theorem 3.4 in \cite{CFL3} imply:
\begin{coro}
The homogeneous coordinate ring of $X(\tau)$ admits a quadratic Gr\"obner basis and is a Koszul algebra.
\end{coro}

The positivity property of the quasi-valuation implies that the subspaces defined by the filtration associated to $\mathcal V$ (see Section~\ref{maxchainAndquasivalu})
are ideals, so we can consider the associated graded algebra $\mathrm{gr}_{\mathcal V}\mathbb K[\hat X(\tau)]$. Theorems~11.1 and 12.2 in \cite{CFL} 
imply:
\begin{theorem}  
\begin{enumerate}
\item[\it i)] The degenerate algebra $\mathrm{gr}_{\mathcal V}\mathbb K[\hat X(\tau)]$ is isomorphic to the fan algebra $\mathbb K[\mathrm{LS}^+_{\lambda}]$
(see Definition~\ref{Defn:FanAlgebra}).
\item[\it ii)] There exists a flat degeneration of $X(\tau)$ into $X_0$, a union of projective toric varieties. Moreover, the special fibre 
$X_0$ is equidimensional, it is isomorphic to $\mathrm{Proj}(\mathrm{gr}_{\mathcal V}\mathbb K[\hat X(\tau)])$,
and its irreducible components are normal varieties and  in bijection with maximal chains in $A_\tau$.
\end{enumerate}
\end{theorem}

\subsection{Straightening laws}
Let $\mathbb B(V(\lambda)_\tau)$ be a fixed basis of $V(\lambda)_\tau^*=\mathbb K[\hat X(\tau)]_1$ consisting of path vectors.
For the special class of path vectors constructed in \cite{L1}
one can find relations in \cite{L1} and \cite{LLM}. But it turns out that these relations actually fit perfectly into the framework of normal and balanced Seshadri stratifications. The following is just a reformulation of part of Proposition 2.20 in \cite{CFL4} in terms of the Seshadri stratification of Schubert varieties and shows that the above mentioned relations can be viewed as a special case:

\begin{proposition}\label{proposition_standard_monomial_basis}
\begin{enumerate}
\item[(i)]
If a monomial $p_{\underline{a}_1}\cdots p_{\underline{a}_n}$ of path vectors is not standard, then there exists a straightening relation expressing it as a linear combination of standard monomials
\[
p_{\underline{a}_1}\cdots p_{\underline{a}_n} = \sum_h u_h p_{\underline{a}_{h,1}}\cdots p_{\underline{a}_{h,{ n}}},
\]
where $u_h\not=0$ only if 
$\underline{a}_1+\ldots+\underline{a}_n \trianglelefteq \underline{a}_{h,1}+\ldots+\underline{a}_{h,{ n}}$.
\item[(ii)] If in \emph{(i)} there exists a chain $\mathfrak{C}$ such that $\supp\underline{a}_i\subseteq\mathfrak{C}$ for all $i=1,\ldots,n$, and $\underline{a}'_1 + \cdots + \underline{a}'_{ n} $ is the decomposition of $\underline{a}_1 + \cdots +\underline{a}_n \in \mathrm{LS}^+_{\lambda}$, then the standard monomial $p_{\underline{a}'_1}\cdots p_{\underline{a}'_{ n}}$ appears in the right side of the equation in \emph{(i)} with a non-zero coefficient.
\end{enumerate}
\end{proposition}
To rewrite a monomial in $\mathbb{K}[\hat{X}(\tau)]$ that is not standard as a linear combination of standard monomials, one can use 
the subduction algorithm (see for example \cite{KM}). The appropriate notation in this context is in the language of Khovanskii basis. 
For details concerning the subduction algorithm in the context of Seshadri stratifications we refer to \cite{CFL}. Recall that the 
quasi-valuation $\mathcal V$ depends on the choice of a total order  $\geq_t$ on $A_{\tau}$ refining the Bruhat order
and preserving the length function. To emphasize this, we write $\mathcal{V}_{\geq_t}$.

\begin{definition}
\begin{enumerate}
\item[{\it i)}] A subset $\mathbb{G}\subseteq \mathbb{K}[\hat{X}(\tau)]$ is called a \emph{Khovanskii basis for the quasi-valuation} $\mathcal{V}_{\geq_t}$, if the image of 
$\mathbb{G}$ in $\mathrm{gr}_{\mathcal{V}_{\geq_t}} \mathbb{K}[\hat{X}(\tau)]$ generates the algebra $\mathrm{gr}_{\mathcal{V}_{\geq_t}}\mathbb{K}[\hat{X}(\tau)]$.
\item[{\it ii)}] A subset $\mathbb{G}\subseteq \mathbb{K}[\hat{X}(\tau)]$ is called a \emph{Khovanskii basis for a Seshadri stratification}, if it is a Khovanskii basis for all possible $\mathcal{V}_{\geq_t}$, where $\geq_t$ is a linear extensions of $\geq$ satisfying: if $\ell(p) > \ell(q)$ then $p>_t q$.
\end{enumerate}
\end{definition} 
Let $\mathbb B(V(\lambda)_\tau)$ be a fixed basis of $V(\lambda)_\tau^*=\mathbb K[\hat X(\tau)]_1$ consisting of path vectors.
Theorem~\ref{some_vanishingA}, Theorem~\ref{LScandidat_FanOfMonoids} and Theorem~\ref{proposition_standard_monomial_basis1} together imply:
\begin{coro}
The set $\mathbb B(V(\lambda)_\tau)$ is a Khovanskii basis 
for the Seshadri stratification defined in Section~\ref{Seshadri_stratification}.
\end{coro}
\subsection{Compatibility with the strata}
We refer to \cite{CFL} for more details, we just quote the results which hold in general for Seshadri stratifications.
Since we consider two Schubert varieties at the same time, we add a $\tau$ or $\kappa$ as index, for example
we write $\mathrm{LS}^+_{\lambda,\tau}$ instead of $\mathrm{LS}^+_{\lambda}$. We fix basis $\mathbb B(V(\lambda)_\tau)$ of $V(\lambda)_\tau^*=\mathbb K[\hat X(\tau)]_1$ 
consisting of path vectors.
\begin{definition}\label{compatibility_subvariety}
Let $\kappa\in A_\tau$. A standard monomial $p_{\underline{a}^1}\cdots p_{\underline{a}^m}\in \mathbb K[\hat X(\tau)]$
on $\hat X(\tau)$ is called \emph{standard on $\hat X(\kappa)$} if  $\max\supp{\underline{a}^1}\leq \kappa$.
\end{definition}

Let $\kappa\in A_\tau$ be such that $\kappa<\tau$. 
By Example~\ref{induction}, we know that the collection of subvarieties $X(\delta)$, $\delta\le \kappa$, 
and the extremal functions $p_\delta$ for $\delta\in A_\kappa$ satisfy the conditions (S1)-(S3), and hence defines a Seshadri 
stratification for $X(\kappa)\hookrightarrow  \mathbb P(V(\lambda)_\tau)$.  Let $\ge_t $ be the total order on $A_\tau$ chosen
in the construction of $\mathcal V$. We keep the notation $\ge_t $ for the induced total order on $A_\kappa$. Let $\mathcal V_\kappa$ be
the associated quasi-valuation on $ \mathbb K[\hat X(\kappa)]$.
In this setup, one gets various natural objects associated to the subset
$A_\kappa\subseteq A_\tau$: 
\begin{itemize}
\item $\mathrm{LS}^+_{\lambda,\kappa}(1):=\{\underline{a}\in \mathrm{LS}^+_{\lambda,\tau}(1)\mid \supp \underline{a}\in A_\kappa\}$, 
\item $\mathbb B(V(\lambda)_\kappa)=\{p_{\underline{a}}\vert_{\hat X(\kappa)} \mid \underline{a}\in \mathrm{LS}^+_{\lambda,\kappa}(1)\}$, and
\item $\mathrm{LS}^+_{\lambda,\kappa}:=\{\underline{a}\in \mathrm{LS}^+_{\lambda,\tau}\mid \supp \underline{a} \subseteq A_\kappa\}$.
\end{itemize}
A natural question arises: What is the connection between these objects and the fan of monoids associated to $\mathcal V_\kappa$, its generating set, etc? 

We consider the vector space $\mathbb Q^{A_\kappa}$ as a subspace of $\mathbb Q^{A_\tau}$. Since every
maximal chain in $A_\kappa$ is a chain in $A_\tau$, by abuse of notation we write $\mathcal V_\kappa(f)\in \mathbb Q^{A_\tau}$ 
for a non-zero function $f\in \mathbb K[\hat{X}(\kappa)]$. For a proof see \cite{CFL}.

\begin{theorem}\label{prop:SMT:for:subvarieties}
Since the Seshadri stratification of $X(\tau)$ is balanced and normal, the following holds:
\begin{enumerate}
    \item[\it i)]\label{balancednormal} for all $\kappa\in A_\tau$, the induced Seshadri stratification on $X(\kappa)$ is balanced and normal;
    \item[\it ii)]\label{fan} the fan of monoids associated to $\mathcal V_\kappa$ is equal to $\mathrm{LS}^+_{\lambda,\kappa}$, $\mathrm{LS}^+_{\lambda,\kappa}(1)$ is 
     its generating set of indecomposables  and $\mathbb B(V(\lambda)_\kappa)$ 
    is a Khovanskii basis for the Seshadri stratification of $X(\kappa)$;
    \item[\it iii)]\label{valuation-compatibility} if $p_{\underline{a}^1}\cdots p_{\underline{a}^m}$ is a standard monomial, standard on $X(\kappa)$, 
    then $\mathcal V_\kappa(p_{\underline{a}^1}\cdots p_{\underline{a}^m}\vert_{\hat X(\kappa)})=\mathcal V(p_{\underline{a}^1}\cdots p_{\underline{a}^m})=
    \sum_{i=1}^m \underline{a}^i$;
    \item[\it iv)]\label{restriction} the restrictions of the standard monomials $p_{\underline{a}^1}\cdots p_{\underline{a}^m}\vert_{\hat X(\kappa)}$, standard on $X(\kappa)$, form a basis of $\mathbb K[\hat{X}(\kappa)]$;
    \item[\it v)]\label{vanishing} a standard monomial $p_{\underline{a}^1}\cdots p_{\underline{a}^m}$ on $\hat X(\tau)$ vanishes on the subvariety $\hat X(\kappa)$ 
    if and only if $p_{\underline{a}^1}\cdots p_{\underline{a}^m}$ is not standard on $\hat X(\kappa)$;
    \item[\it vi)]\label{ideal} the vanishing ideal $\mathcal I(\hat X(\kappa))$ of $\hat{X}(\kappa)$ in $\mathbb K[\hat{X}(\tau)]$ is generated by the elements in 
    $\mathbb B(V(\lambda)_\tau)\setminus \mathbb B(V(\lambda)_\kappa)$, 
    and the ideal has as vector space basis the set of all standard monomials on $\hat X(\tau)$ 
    which are not standard on $\hat X(\kappa)$;
    \item[\it vii)]\label{intersection} for all pairs of elements $\kappa,\delta\in A_\tau$, the scheme theoretic intersection $X(\kappa)\cap X(\delta)$ is reduced. 
    It is the union of those subvarieties $X(\xi)$ such that $\xi\le \kappa$ and $\xi\le \delta$,
    endowed with the induced reduced structure. 
\end{enumerate}
\end{theorem} 
\subsection{Compatibility with the standard monomial theory in \cite{L1}}
Let $\mathbb B(V(\lambda)_\tau)$ be a fixed basis of $V(\lambda)_\tau^*=\mathbb K[\hat X(\tau)]_1$ consisting of path vectors.
All the results stated in Sections~\ref{standard-monomial-theory} and \ref{standard-monomial-theory2} so far assume that 
one fixes at the beginning such a basis consisting of path vectors, no other property is assumed.

In \cite{L1} too, using the quantum Frobenius morphism, a basis of $V(\lambda)_\tau^*=\mathbb K[\hat X(\tau)]_1$ was constructed,
which was the foundation for the standard monomial theory in \textit{ibidem}. By Lemma~\ref{lpathvector}, these vectors 
$p_{\underline a,\ell}$ (for the notation see Appendix III, Section~\ref{The_path_vectors}) are path vectors in the sense of Definition~\ref{define_path_vector}.
It follows that the basis constructed in \textit{ibidem} is a special choice of a basis $\mathbb B(V(\lambda)_\tau)$
of $V(\lambda)_\tau^*=\mathbb K[\hat X(\tau)]_1$ consisting of path vectors.

The notion of a standard monomial used in \cite{L1} is based on the definition of a standard concatenation of LS-paths,
see Definition~\ref{standardconcat}, and the fact that an LS-path of shape $m\lambda$ can always be written as a
standard concatenation of LS-paths of shape $\lambda$, see Proposition~\ref{decompLSpath}. Now the bijection
$\Theta$ between the LS-paths and the fan of monoids described in Appendix I, Section~\ref{The_path_vectors}, { preserves} the notion of standardness.
{ That is to say,} it sends a standard concatenation of LS-paths to a standard sum of elements in the fan of monoids $\textrm{LS}_\lambda^+$,
see Definition~\ref{AppendixI_standardsum} and the text thereafter. It follows that the notion of a standard monomial
in \cite{L1} for special products of the vectors $p_{\underline a,\ell}$ coincides with the notion of a standard monomial in Definition~\ref{def:standard:monomial}. 
Summarizing we have:
\begin{proposition}
The standard monomial theory described in \cite{L1} is the same as in Theorem~\ref{proposition_standard_monomial_basis1},
only the possible choice of path vectors is restricted  in \cite{L1} to those obtained via the quantum Frobenius splitting
(see Appendix III, Section~\ref{The_path_vectors}).
\end{proposition}

See \cite{CFL2} for another approach to this proposition in the framework of LS-algebras.
}
\section{Projective normality}\label{section:projective_normality}
There are meanwhile many different proofs of the normality of Schubert varieties respectively
the projective normality of embedded Schubert varieties. A geometric proof can be found, for example,
in \cite{S4} (finite type case), the proofs in \cite{RR} (finite type case) and \cite{Ma2} (Kac-Moody groups) 
use  (the algebraic geometric) Frobenius splitting, the proof in \cite{L1} uses the standard monomial theory
developed in the same paper.

We give here another proof, which is a direct application of the theory of Seshadri stratifications. 
Let $\mathrm{SR}(A_\tau)$ be the Stanley-Reisner algebra of the poset $A_\tau$. By definition
$$
\mathrm{SR}(A_\tau) := \mathbb{K}[t_\sigma \mid \sigma\in A_\tau]/(t_\sigma t_\kappa \mid \sigma,\kappa \mathrm{\ not\ comparable}).
$$
By \cite{Bj,BjW}, the Stanley-Reisner algebra $\mathrm{SR}(A_\tau)$ is Cohen-Macaulay over any field $\mathbb K$.
By Theorem~\ref{LScandidat_FanOfMonoids} we know that the Seshadri stratification is normal. By Theorem 14.1 in \cite{CFL},
the normality of the stratification and the Cohen-Macaulayness of the Stanley-Reisner algebra $\mathrm{SR}(A_\tau)$
imply:
\begin{theorem}
The embedded Schubert variety $X(\tau)\hookrightarrow \mathbb P(V(\lambda)_\tau)$ is projectively normal.
\end{theorem}
Fix a dominant weight  $\lambda$ and let $Q\supset B$ be the parabolic subgroup of 
$G$ associated to $\lambda$. We identify the dual space 
$V(\lam)^*$ with the space of global sections 
$\mathrm{H}^0(G/Q,\mathcal L_\lambda)$ of the line bundle 
$\mathcal L_\lambda:=G\times_Q \mathbb{K}_{-\lambda}$.
Let $\phi: G/Q\hookrightarrow {\mathbb P}(V(\lam))$ be the 
corresponding embedding.

The projective normality implies: the natural morphisms $\mathbb K[X(\tau)]_m\rightarrow \mathrm{H}^0(X(\tau),\mathcal L_{m\lambda})$ are isomorphisms for all $m\ge 0$. 
So one can reformulate the standard monomial theory in terms of sections. We use the same notation as in the section before.

\begin{theorem}\label{prop:SMT:for:subvarieties2}
\begin{enumerate}
    \item[\it i)] The standard monomials $p_{\underline{a}^1}\cdots p_{\underline{a}^m}$ of degree $m$ with $\underline{a}^1,\ldots,\underline{a}^m\in\mathrm{LS}^+_{\lambda}(1)$ form a basis for $\mathrm{H}^0(X(\tau),\mathcal L_{m\lambda})$.
    \item[\it ii)] A standard monomial $p_{\underline{a}^1}\cdots p_{\underline{a}^m}$ on $X(\tau)$ vanishes on the subvariety $X(\kappa)$ for $\kappa\le \tau$
    if and only if the monomial is not standard on $\hat X(\kappa)$.
    \item[\it iii)] For $\kappa\le \tau$, the restrictions of the standard monomials $p_{\underline{a}^1}\cdots p_{\underline{a}^m}\vert_{\hat X(\kappa)}$ of degree $m$, 
    standard on $X(\kappa)$, form a basis of $\mathrm{H}^0(X(\tau),\mathcal L_{m\lambda})$.
     \item[\it iv)] For $\kappa\le \tau$ and all $m\ge 1$, the restriction map $\mathrm{H}^0(X(\tau),\mathcal L_{m\lambda})\rightarrow \mathrm{H}^0(X(\kappa),\mathcal L_{m\lambda})$
     is surjective.
      \item[\it v)] For all $m\ge 1$, the multiplication map $S^m \mathrm{H}^0(X(\tau),\mathcal L_\lam)\rightarrow \mathrm{H}^0(X(\tau),\mathcal L_{m\lam})$ is surjective.
\end{enumerate}
\end{theorem}
\begin{proof}
Part \textit{i)--iii)} of the theorem are just reformulations of the corresponding results in Theorem~\ref{prop:SMT:for:subvarieties}. Part \textit{v)} is an immediate
consequence of \textit{i)}, part \textit{iv)} is an immediate consequence of \textit{i)} and \textit{ii)}.
\end{proof}
{ 
\begin{remark} Another consequence of standard monomial theory is the vanishing of the higher cohomology: 
for all $m\ge 1$ and $i\ge 1$, $\mathrm{H}^i(X(\tau),\mathcal L_{m\lambda})=0$. 
    For a proof see, for example, \cite{L1}, Theorem 7, or \cite{LLM}, Theorem 6.4, where the proof is given even for unions of Schubert varieties.
\end{remark}
  }  
\section{Newton-Okounkov simplicial complex}\label{section:Newton_Okounkov}
Newton-Okounkov bodies generalize Newton polytopes from toric geometry
to a more general algebro-geometric as well as representation-theoretic setting.
The Newton-Okounkov  body is constructed from an embedded projective variety $X\hookrightarrow \mathbb P(V)$,
its homogeneous coordinate ring $R$, and a valuation $\nu$ on $R$.
Roughly speaking, in the nicest of all cases, the Newton-Okounkov body associated to a projective variety is
a polytope. It is known that in such a case there exists a flat toric degeneration from $X$ to the toric variety
corresponding to this polytope. The Newton-Okounkov body stores a lot of information about the original 
variety, for example the degree of the embedded variety.

Since instead of a valuation, we have only a quasi-valuation, it is not straightforward to define a Newton-Okounkov body.  But one can get an appropriate replacement. We formulate the construction here only for the case of Schubert varieties. See \cite{CFL} for more details in the general case.

Let $\mathcal{C}$ be the set of all maximal chains in $A_\tau$. 
 
\begin{definition}\label{Defn:NOSC}
The \emph{Newton-Okounkov simplicial complex} $\Delta_{\mathcal V}$ associated to the quasi-valuation $\mathcal{V}$ is defined as
$$
\Delta_{\mathcal V}:=\bigcup_{\mathfrak C\in\mathcal{C}}
\overline{ 
\bigcup_{m\ge 1}
\left\{
\frac{1}{m}\,\underline{a}\mid\underline{a}\in \Gamma_{\mathfrak C},\deg \underline{a}=m \right\} 
}\subseteq \mathbb R^{A_\tau}.
$$
\end{definition}
Another object associated to the Seshadri stratification is the order complex $\Delta(A_\tau)$ attached to the poset $(A_\tau,\leq)$.
It is the simplicial complex having the set $A_\tau$ as vertices and all chains $C\subseteq A_\tau$ as faces, i.e. 
$\Delta(A_\tau)=\{C\subseteq A_\tau\mid C\text{ is a chain}\}$.
A geometric realization of $\Delta(A_\tau)$ can be constructed by intersecting the cones $K_C$ (see Section \ref{higher_valued_valuation}) with appropriate hyperplanes:
for a chain $C\subseteq A_\tau$ denote by $\Delta_{C}\subseteq  \mathbb R^{A_\tau}$
the simplex: 
\begin{equation}\label{chainsimplex}
\Delta_{C}:=\textrm{convex hull of\ }\left\{ e_p\mid p\in C\right\}.
\end{equation}
The union of the simplexes 
$$
|\Delta(A_\tau)|:=\bigcup_{C\subseteq A_\tau\textrm{\ chain}}\Delta_{C}\subseteq\mathbb{R}^{A_\tau}$$ 
is the desired \emph{geometric realization} of $\Delta(A_\tau)$. The maximal simplexes are those $\Delta_{\mathfrak{C}}$ 
arising from maximal chains $\mathfrak{C}$ in $A_\tau$. It has been shown in \cite{CFL}, Proposition 13.3: 
\begin{proposition}
The Newton-Okounkov simplicial complex $\Delta_{\mathcal V}$ coincides with the 
geometric realization of $\Delta(A_\tau)$: $\Delta_{\mathcal V}=|\Delta(A_\tau)|$. 
In particular, $\Delta_{\mathcal V}$ is a homogeneous simplicial complex of dimension $r$.
\end{proposition}

In the construction of a toric variety associated to a polytope $P\subseteq \mathbb R^l$, the sequence of 
integral points one gets as the intersections $mP\cap \mathbb Z^l$, $m\ge 1$, plays an important role. Instead of 
multiplying the polytope one could also ``shrink'' the lattice and consider the sets $P\cap \frac{1}{m}\mathbb Z^l$, $m\ge 1$.
This leads to the definition of an integral structure on our simplicial complex $\Delta_{\mathcal V}$.

\begin{definition}
An \textit{integral structure} on $\Delta_{\mathcal V}$ is a collection of subsets $\Delta_{\mathcal V}(n)\subset \Delta_{\mathcal V}$ for
all $n\in\mathbb N$ and an affine embedding $i_{\Delta}: \Delta \rightarrow \mathbb R^{\dim\Delta}$ for each 
simplex $\Delta$ in $\Delta_{\mathcal V}$ such that:
\begin{itemize}
\item  the vertices of $i(\Delta)$ have integral coordinates,
\item if we denote $\Delta(n)=\Delta\cap \Delta_{\mathcal V}(n)$, then
$$
i _\Delta(\Delta(n))=\{v\in i _\Delta(\Delta(n))\mid nv\in \mathbb Z^{\dim\Delta}\}.
$$
\end{itemize}
\end{definition}
The goal is of course to have an integral structure on $\Delta_{\mathcal V}$ such that for all $m\in\mathbb N$:
$$
\Delta_{\mathcal V}(m)= \left\{
\frac{1}{m}\,\underline{a}\mid\underline{a}\in \Gamma_{\mathfrak C},\deg \underline{a}=m \right\},
$$
and to generalize in this way the tools provided by the Newton-Okounkov theory to the setting of Newton-Okounkov
simplicial complexes.

In \cite{CFL} we have constructed such an integral structure in the general case. We will see that in the case of Schubert varieties,
the construction recovers the one given by Dehy in \cite{De}.

Let $\frak C=\{\tau=\tau_r>\ldots>\tau_0=\mathrm{id}\}$ be a maximal chain in $A_\tau$, and for $k=1,\ldots,r$  let $\beta_k$ be the positive 
real root such that $\tau_k=s_{\beta_k}\tau_{k-1}$ and let $b_k=\langle \tau_{k-1}(\lambda),\beta_{k}^\vee\rangle$ be the bond. 

Denote by $\Delta_{\frak C}$ the maximal simplex associated to
the given maximal chain. This means in term of the order complex: $\Delta_{\frak C}\subseteq \mathbb R^{\frak C}\subseteq \mathbb R^{A_\tau}$
is the convex hull of the $e_{\tau_i}$, $i=0,\ldots,r$. So a typical element in $\Delta_{\frak C}$ is of the form $\sum_{i=0}^r a_i e_{\tau_i}$,
where the $a_i\ge 0$ are such that $\sum_{i=0}^r a_i =1$.

Dehy defined the map $i_{\Delta_{\frak C}}:\Delta_{\frak C} \rightarrow \mathbb R^{r}$ as follows:
$$
i_{\Delta_{\frak C}}(e_{\tau_j})=\left\{
\begin{array}{ll}
0,&\textrm{if $j=0$},\\
\sum_{k=1}^j b_k e_k, &\textrm{otherwise,}
\end{array}\right.
$$
and the map is extended to the simplex by: 
$$i_{\Delta_{\frak C}}\left(\sum_{k=0}^r a_je_{\tau_j}\right)=\sum_{k=0}^r a_ji_{\Delta_{\frak C}}(e_{\tau_j}).$$

If $\Delta'$ is not a maximal simplex, then it is contained in some maximal simplex $\Delta_{\frak C}$ and we define
$i_{\Delta'}=i_{\Delta_{\frak C}}\vert_{\Delta'}$. Dehy showed in \cite{De}: these maps glue together and provide a well defined integral 
structure on $\Delta_{\frak C}$ in the sense of Definition~\ref{Defn:NOSC}. 

\begin{lemma}
If $\Delta_{\frak C}$ is the maximal simplex in $\Delta_{\mathcal V}$ associated to the maximal chain $\mathfrak C$, then 
for the integral structure constructed above holds:
$$
\Delta_{\frak C}(m)=\left\{\left.\frac{1}{m}\,{\underline{a}}\,\right\vert\,\underline{a}\in \Gamma_{\mathfrak C},\deg \underline{a}=m \right\}.
$$
\end{lemma}
\begin{proof}
Let $\mathfrak C=\{\tau=\tau_r>\ldots>\tau_0=\mathrm{id}\}$ be the maximal chain, and for $k=1,\ldots,r$  let $\beta_k$ be the corresponding
positive real root and $b_k$ the bond. We have
$$
\Delta_{\frak C}(m)=\left\{\left.\underline{a}=\sum_{i=0}^r a_i e_{\tau_i}\in  \Delta_{\frak C}\right\vert m \big(i_{\Delta_{\frak C}}(\underline{a})\big)\in \mathbb Z^r\right\}.
$$
Now for $\underline{a}=\sum_{i=0}^r a_i e_{\tau_i}\in \Delta_{\frak C}$ we have 
$$
m\big(i_{\Delta}(\sum_{i=0}^r a_i e_{\tau_i})\big)=m\big(\sum_{j=1}^r a_j(\sum_{k=1}^j b_k e_k) \big)=m\big(\sum_{k=1}^rb_k(a_k+\ldots+a_r)e_k\big)\in \mathbb Z^r
$$
if and only if $m\underline{a} \in \mathrm{LS}^+_{\mathfrak C,\lambda}(m)$. Indeed,
recall that $\sum_{i=0}^r a_i=1$ and $a_i\ge 0$ for all $i=0,\ldots,r$. So the degree condition $m(a_0+\ldots+a_r)=m$ is automatically satisfied.
\end{proof}
\begin{coro}
The Newton-Okounkov simplicial complex $\Delta_{\mathcal V}$ admits an integral structure such that $\Delta_{\mathcal V}(m)=\{\frac{1}{m}\underline{a}\mid
\underline{a}\in\Gamma(m)\}$.
\end{coro}

The degree formula in \cite{CFL}, Theorem 13.6, has in this setting the following formulation:

\begin{proposition} 
The degree of the embedded Schubert variety 
$X(\tau)\subseteq \mathbb P(V(\lambda)_\tau)$ is equal to the sum $\sum_{\mathfrak C}\prod_{j=1}^r b_j$ running over all maximal chains in $A_\tau$, and $\prod_{j=1}^r b_j$ is the product of all bonds along the maximal chain $\frak C$.
\end{proposition}

The degree formula can also be found in \cite{Ch}, and in \cite{Kn}  in a
symplectic context.

\section{Appendix I: LS-paths and the LS-lattice}\label{appendix1}
We fix  a dominant weight $\lam\in\Lambda^+$. Let $Q\subseteq G$ be the standard parabolic
subgroup associated to $\lambda$, i.e., $Q$ is generated by the Borel subgroup $B$ and 
the root subgroup $U_{-\alpha}$ for all simple roots $\alpha$ such that $\langle\lambda,\alpha^\vee\rangle=0$.
The Weyl group of $G$ is denoted by $W$, the Weyl group of $Q$ by $W_Q$.

In this section we recall the notion of a Lakshmibai-Seshadri path of shape $\lambda$, or, 
for short an {\it LS-path} \index{LS-path} of shape $\lambda$,
and we explain the connection with the fan of monoids $\mathrm{LS}_\lambda^+$ introduced in Section~\ref{CandidatForMonoid}.

\subsection{Chains of elements in $W/W_Q$}

We start with the notion of a $(d,\lambda)$-chain as it was introduced by Lakshmibai \cite{LS2} 
(see also \cite{SMT,S2}).

A maximal chain joining two elements $\kappa$ and $\sigma\in W/W_Q$, $\kappa>\sigma$, is a pair of sequences
$(\kappa_t,\ldots,\kappa_0;\beta_t,\cdots,\beta_1)$, where
\begin{itemize}
\item $\kappa=\kappa_t>\cdots >\kappa_0=\sigma$ is a chain of elements in $W/W_Q$ and 
\item $\beta_t,\cdots,\beta_1$ are positive real roots such that
\item $s_{\beta_i}\kappa_{i-1}=\kappa_i$ and $\ell(\kappa_{i})=\ell(\kappa_{i-1})+1$, $i=1,\ldots,t$. 
\end{itemize}
\begin{definition}
Given a positive rational number $d$,
the chain $(\kappa_t,\ldots,\kappa_0;\beta_1,\cdots,\beta_r)$ is called a {\it $(d,\lambda)$-chain} joining $\kappa$ and $\sigma$ if in addition
$d \langle \kappa_i(\lam),\beta^\vee_i \rangle\in\Z$ for all $i=1,\ldots,t$.
\end{definition} 
It has been shown in \cite{De} that  if one maximal chain between $\kappa$ and $\sigma$ has the property of being a $(d,\lambda)$-chain,
then all maximal chains joining $\kappa$ and $\sigma$ are $(d,\lambda)$-chains.

In many jointed papers, Seshadri, Lakshmibai, Musili and others  (see \cite{LSII,LSIII, LSIV, LS, LS2,S2})
developed in the framework of the standard monomial theory case by case a combinatorial character formula,
leading Lakshmibai to the notion of what is now called an LS-path.
The name \textit{LS-path} came up in connection with the theory of path models of representations,
which was developed in \cite{L2,L3}.

\begin{definition}\label{LSpath}
An {\it LS-path} $\pi=(\sigma_p,\sigma_{p-1},\ldots,\sigma_1; 0,d_p,\ldots,d_1=1)$ \textit{of shape} $\lambda$ \index{$\pi$, LS-path}
is a pair of sequences of elements in $W/W_Q$ and rational numbers such that
\begin{itemize}
\item $\underline{\sigma}:\sigma_p>\sigma_{p-1}>\ldots>\sigma_1$ 
is a linearly ordered sequence of elements in $W/W_Q$,
\item $\underline{d}:0<d_p<\ldots<d_1=1$ is a sequence of
rational numbers, 
\item for all $i=2,\ldots,p$, there exists a $(d_i,\lambda)$-chain joining $\sigma_{i}$ and 
$\sigma_{i-1}$.
\end{itemize}
The {\it support $\supp\pi$ of an LS-path $\pi$}\index{$\supp\, \pi$; support of a LS-path} is the set $\{\sigma_p,\sigma_{p-1},\ldots,\sigma_1\}$.
\end{definition}
\subsection{Interpretation in terms of a path model}
In \cite{L2}, the definition of LS-paths was transferred from the context of Weyl group combinatorics into the setting of 
combinatorics of piecewise linear paths in $\Lambda_{\mathbb R}=\Lambda\otimes_{\mathbb Z}\mathbb R$.
Let $\Pi$ be the set of all piecewise linear paths $\phi:[0,1]\rightarrow \Lambda_{\mathbb R}$
starting in $0$ and ending in $\phi(1)\in \Lambda$. (Two paths $\pi_1,\pi_2$ are identified
if there exists a nondecreasing, surjective, continuous map
$\psi: [0, 1] \rightarrow [0, 1]$ (for short: a reparameterization) such that $\pi_1 =\pi_2 \circ \psi$).
The endpoint $\phi(1)$ of a path is called the weight of the path.
In \cite{L2,L3}, for all simple roots $\alpha$ folding operators $e_\alpha,f_\alpha$ are defined on the set $\Pi$. 
For a sufficiently nice path $\pi_0$ ending in a dominant integral weight $\lambda=\pi_0(1)$, let $\mathbb B(\pi_0)$\index{$\mathbb B(\pi_0)$, path model} be  
\textit{the smallest set of piecewise linear paths which contains the chosen starting path $\pi_0$ and is stable 
under all folding operators}. This set is called a \textit{path model}\index{path model} for the representation $V(\lambda)$. 

An example for a ``nice'' path is the straight line path $\pi_{\lambda}$ which joins the origin with the dominant weight $\lambda$.
It turns out  that in this special case the set $\mathbb B(\lambda)$ (or rather $\mathbb B(\pi_\lambda)$ to be precise) 
\index{$\mathbb B(\lambda)$, path model given by LS-paths}
is equal to the set of LS-paths of shape $\lambda$ as in Definition~\ref{LSpath}. 
To make this more precise, we have to give an interpretation of a LS-path as a piecewise linear path in $\Lambda_{\mathbb R}$.

Let $\pi=(\sigma_p,\sigma_{p-1},\ldots,\sigma_1; d_{p+1}=0,d_p,\ldots,d_1=1)$ be an LS-path of shape $\lambda$.
We consider $\pi:[0,1]\rightarrow \Lambda_{\mathbb R}$ as the piecewise linear path \cite{L2} defined by :
$$
 t\mapsto  (d_p-d_{p+1})\sigma_p(\lambda)+\ldots+(d_{j}-d_{j+1})\sigma_j(\lambda)+
 (t-d_{j})\sigma_{j-1}(\lambda)\quad\text{for\ }d_{j}\le t\le d_{j-1}.
$$
The endpoint of the path is 
\begin{equation}\label{weightOfpaths}
\pi(1)= \sum_{i=1}^{p} (d_i-d_{i+1})\sigma_i(\lambda).
\end{equation}
The proof of the bijection: 
$$
\textrm{LS-paths of shape $\lambda$ as in Definition~\ref{LSpath}} \quad
\mathop{\longleftrightarrow}^{1:1} 
\quad\textrm{path model $\mathbb B(\lambda)$}
$$
can be found in \cite{L2}. The folding operators play a crucial role in the proof. 

 The following combinatorial character formula was conjectured by Lakshmibai,
proved in many special cases by Seshadri, Lakshmibai, Musili and others in \cite{LSII,LSIII, LSIV, LS, LS2},
and proved in full generality in \cite{L2,L3} in the framework of the path model theory. 
Let $\rho\in \Lambda^+$ be an element such that $\langle\rho,\alpha^\vee\rangle=1$
for all simple roots $\alpha$, and let $D_\alpha$ be the Demazure operator on the group ring
$\mathbb Z[\Lambda]=\mathbb Z[e^\mu\mid \mu\in\Lambda]$ associated to a simple root $\alpha$:
\begin{equation}\label{Demazure_operator}
D_\a (e^\mu):={{e^{\mu+\rho}-e^{s_\a(\mu+\rho)}\over{1-e^{-\a}}}}e^{-\rho}.
\end{equation}
Let \index{$\mathbb{B}_{\tau}$, path model of a Demazure module $V(\lambda)_{\tau}$}
$$
{\mathbb{B}(\lambda)_{\tau}}:=\{\pi=(\sigma_p,\sigma_{p-1},\ldots,\sigma_1; 0,d_p,\ldots,d_1=1)\mid\hbox{\rm $\pi$ LS-path of shape $\lambda$, $\sigma_p\le\tau$}\}
$$ 
be the set of all LS-paths of shape $\lambda$ starting with an element $\sigma_p$ which is smaller or equal to $\tau$.
The following character formula was proved in \cite{L2}:
\begin{theorem}\label{path-character}
If $\tau=s_{i_1}\cdots s_{i_r}$ be a reduced decomposition, then
$$
D_{\a_{i_1}}\circ\cdots\circ D_{\a_{i_r}} e^\lambda=\sum_{\pi\in {\mathbb{B}(\lambda)_{\tau}}} e^{\pi(1)}.
$$
\end{theorem}
It is well known that the projective normality of Schubert varieties is strongly related to Demazure's character formula, see, for example
\cite{S3} and \cite{Jo}.
Or, in other words:
$$
\textrm{Char\,}(\mathrm{H}^0(X(\tau),\mathcal L_\lambda)^*)=D_{\a_{i_1}}\circ\cdots\circ D_{\a_{i_r}} e^\lambda,
$$
holds for all very ample line bundles $\mathcal L_\lambda$ on $G/Q$ if and only if $X(\tau)\subseteq \mathbb P(V(\lambda)_\tau)$ 
is projectively normal for all these $\lambda$. One has to be a bit cautious with the formulation if $G$ is not of finite type.
Over an arbitrary algebraically closed field, the proof of the projective normality {is an essential part in the construction}
of Schubert varieties, see \cite{Ma1,Ma2}.

Assuming the normality of Schubert varieties, this formula shows in connection with Demazure's character formula: 
the set of LS-paths in ${\mathbb B}(\lambda)_{\tau}$ provides a combinatorial character formula for Demazure modules,
and it shows that this set provides an indexing system for a basis of the Demazure module 
$V_\tau(\lambda)$. Or, in its dual version, $\mathbb B(\lambda)_\tau$ provides an indexing system for a basis of 
$\mathrm{H}^0(X(\tau),\mathcal L_\lambda)$, which is of course what Seshadri and Lakshmibai had in mind.

On the set $\Pi$ one has a product, the concatenation of paths. Given $\pi_1,\pi_2\in\Pi$,
by the product $\pi:=\pi_1*\pi_2$\index{$\pi_1*\pi_2$, concatenation of paths} we mean the piecewise linear path obtained by concatenation:
$$ 
\pi(t):=\left\{
\begin{array}{ll}
\pi_1(2t),& \textrm{if\ } 0\le t\le \frac{1}{2};\\
\pi_1(1)+\pi_2(2t-1),& \textrm{if\ } \frac{1}{2}\le t\le 1.\\
\end{array}\right.
$$
Given an LS-path $\pi\in \mathbb B(\lambda)$, then let $i(\pi)$\index{$i(\pi)$, maximal element in the support} 
be the maximal element in $\supp\pi$
and let $e(\pi)$\index{$e(\pi)$, minimal element in the support} be the minimal element in $\supp\pi$. 
We write $\supp\pi_1\ge \supp \pi_2$ if $e(\pi_1)\ge i(\pi_2)$.

\begin{definition}\label {standardconcat}
Given two LS-paths $\pi_1,\pi_2\in \mathbb B(\lambda)$, the concatenation $\pi_1*\pi_2$
is called \textit{standard} if  $\supp\pi_1\ge \supp \pi_2$.\index{standard, concatenation of paths}
\end{definition}

Let $\pi_\lambda$ be again the straight line path joining the origin with the dominant weight $\lambda$.  
The path $\pi_{2\lambda}$ coincides with the concatenation $\pi_\lambda*\pi_\lambda$.
It has been shown in \cite{L2} that the path model $\mathbb B(2\lambda)$ coincides
with the path model $\mathbb B(\pi_\lambda*\pi_\lambda)$ obtained by applying the folding operators to the 
concatenation $\pi_\lambda*\pi_\lambda$. And the elements in $\mathbb B(\pi_\lambda*\pi_\lambda)$ are exactly
the concatenations $\pi:=\pi_1*\pi_2$,  $\pi_1,\pi_2\in  \mathbb B(\lambda)$, which are standard.
Proceeding by induction one can show:
\begin{proposition}\label{decompLSpath}
Every LS-paths of shape $m\lambda$ can be decomposed as a standard product of $m$ LS-paths
of shape $\lambda$. In other words: for any $\pi\in B(m\lambda)$ one can find LS-paths $\pi_1,\ldots,\pi_m\in B(\lambda)$
such that $\supp \pi_1\ge \supp\pi_2\ge \ldots\ge \supp \pi_m$, and $\pi$ is (up to reparameterization)
equal to $\pi_1*\pi_2*\cdots*\pi_m$.
\end{proposition}

\subsection{Interpretation in terms of a fan of monoids}
Recall the situation in Section~\ref{maxchainAndquasivalu}: we have a Schubert variety $X(\tau)\subseteq G/Q$,
embedded in the projective space over a Demazure module $\mathbb P(V(\lambda)_\tau)$, endowed with a Seshadri stratification.
The partially ordered set $A_\tau$ is the set $\{\sigma\in W/W_Q\mid \sigma\le \tau\}$.

Let $\mathcal V$ be the associated quasi-valutation, which depends on a choice of a total order on $A_\tau$, see Section~\ref{higher_valued_valuation}.
We introduced in Section~\ref{CandidatForMonoid}  the fan of monoids 
$\mathrm{LS}_{\lambda}^+=\mathrm{LS}_{\lambda}\cap\mathbb Q_{\ge 0}^{A_\tau}$ as a 
candidate for the fan of monoids $\Gamma$ associated to $\mathcal V$. Note that $\mathrm{LS}_{\lambda}^+$
does not depend of the choice of the total order!

For $\underline{a}\in \mathrm{LS}_{\lambda}$,
we denote by $\supp\underline{a}$ the support $\{\kappa\in A_\tau\mid a_\kappa\not=0\}$ of $\underline{a}$.
Denote by $\mathrm{LS}_{\lambda}^+(m)$ the subset of elements of degree $m$, i.e.:
$$
\mathrm{LS}_{\lambda}^+(m)=\left\{ \underline{a}=\sum_{\kappa\in A_\tau}a_\kappa e_\kappa \in \mathrm{LS}_{\lambda}^+ \mid \sum_{\kappa\in A_\tau} a_\kappa=m \right\}.
$$
We define now a map $\Theta$ between
$$
\mathbb{B}(m\lambda)_{\tau}=
\left\{
\begin{array}{c}
\textrm{LS-paths $\pi$ of}\\
\textrm{shape\ }m\lambda,\ i(\pi)\le \tau
\end{array}\right\}\quad
\mathop{\longrightarrow}^\Theta\quad
\mathrm{LS}_{\lambda}^+(m)=
\left\{\begin{array}{c}
\textrm{elements of degree $m$ in}\\
\textrm{the fan of monoids\ } \mathrm{LS}_{\lambda}^+
\end{array}\right\}.
$$
The map resembles the weight map in \eqref{weightOfpaths}. The difference is that the map remembers the 
Weyl group elements. For  $m\ge 1$, we define the map  by 
\begin{equation}\label{pathLattice}
\pi=(\sigma_p,\ldots,\sigma_1; 0,d_p,\ldots,d_1=1)\mapsto
\Theta(\pi):=\sum_{j=1}^p (d_j-d_{j+1})me_{\sigma_j}.
\end{equation}
where we set $d_{p+1}=0$.

See \cite{Ch}, Proposition 1 and the discussion before Proposition 32 for a different proof of the following result.

\begin{proposition}\label{LSlatticeLSpath}
The map $\Theta$ induces a bijection between the set of LS-paths $\mathbb B(m\lambda)_{\tau}$ of shape $m\lambda$
and  the set $\textrm{\rm LS}_{\lambda}^+(m)$ of elements of degree $m$ in the fan of monoids $\textrm{\rm LS}_{\lambda}^+$. The map respects the support,
i.e. $\supp\pi=\supp\Theta(\pi)$.
\end{proposition}

\begin{proof}
The map $\Theta$ is obviously an injective map from $\mathbb B(m\lambda)_{\tau}$ to $\mathbb Q^{A_\tau}$. 
For a given path $\pi=(\sigma_p,\sigma_{p-1},\ldots,\sigma_1; 0,d_p,\ldots,d_1=1)\in\mathbb B(m\lambda)_{\tau}$, the 
condition $0<d_p<\ldots<d_1=1$ implies $\Theta(\pi)\in \mathbb Q^{A_\tau}_{\ge 0}$, the sum of the coefficients of  
$\Theta(\pi)$ is equal to $m$, and $\supp\pi=\supp\Theta(\pi)$. It remains to show that the image 
is in $\mathrm{LS}_{\lambda}$, and that the map is surjective.

Let $\pi$ be as above. Fix a maximal chain $\mathfrak C=(\kappa_r>\ldots>\kappa_0)$ in 
$A_\tau$ such that $\supp\pi\subseteq \mathfrak C$. 
Since $\pi$ is an LS-path of shape $m\lambda$, for all $i=1,\ldots,p$ there exists a $(d_i,m\lambda)$-chain between 
$\sigma_i$ and $\sigma_{i-1}$. By  \cite{De} we know
that  if one maximal chain between $\sigma_i$ and $\sigma_{i-1}$ is a $(d_i,m\lambda)$-chain,
then all maximal chains joining $\sigma_i$ and $\sigma_{i-1}$ are $(d_i,m\lambda)$-chains. So we can fix the 
$(d_i,m\lambda)$-chains to be subchains of  $\mathfrak C$.

We define an ``extended'' path $\hat \pi=(\kappa_r,\ldots,\kappa_0;\hat d_{r+1}=0,\hat d_r,\ldots,\hat d_0=1)$ as follows:
the sequence of Weyl group elements is the entire maximal chain $\mathfrak C$, and the 
sequence of rational numbers is now only weakly increasing. The rational numbers $\hat d_j$ are given by the following
rules: for $j=1,\ldots,r$ set:

If $\kappa_j>\sigma_p$, then set $\hat d_j=0$. 
If $\sigma_i\ge \kappa_j > \sigma_{i-1}$ for some $i=2,\ldots,p$, then set $\hat d_j=d_i$.
If $\sigma_1\ge \kappa_j$, then set $\hat d_j=1$. 

From the fact that the ``squeezed in'' subchains are $(d_i,m\lambda)$-chains, we conclude that the
condition on $\pi$ to be an LS-path implies for $\hat \pi$:
$$
\forall j=0,\ldots,r:\quad  \hat d_j\langle \kappa_j(m\lambda),\beta^\vee_j\rangle\in\mathbb Z,
$$
where $\beta_j$ is the positive real root such that $s_{\beta_j}\kappa_{j-1}=\kappa_{j}$.

We adapt the definition of $\Theta$ and set:
$\hat\Theta (\hat\pi):=\sum_{j=1}^r (\hat d_j-\hat d_{j+1})me_{\kappa_j}$. By the choice of the $\hat d_j$, $j=1,\ldots,r$, 
we have $\hat\Theta (\hat\pi)=\Theta (\pi)$. Now by comparing the coefficients in the equation:
$$
\Theta (\pi)=\hat\Theta (\hat\pi)=\sum_{j=1}^r (\hat d_j-\hat d_{j+1})me_{\kappa_j}=\sum_{j=0}^r a_j e_{\kappa_j}
$$
we get: 
$$
\begin{array}{lcl}
m\hat d_r&=&a_r\\ 
m\hat d_{r-1}&=&m\hat d_{r}+a_{r-1}= a_r+a_{r-1}\\
m\hat d_{r-2}&=&m\hat d_{r-1}+a_{r-2}= a_r+a_{r-1}+a_{r-2}\\
 \vdots     &\vdots&\vdots\\
m\hat d_{0}&=&m\hat d_{1}+a_0=\sum_{i=0}^r a_i\\
\end{array}.
$$ 
Since $b_j=\langle\kappa_{j-1}(\lambda),\beta^\vee_j\rangle=\vert \langle\kappa_j(\lambda),\beta^\vee_j\rangle\vert$ is the 
bond between $\kappa_{j-1}$ and $\kappa_{j}$, we see that the
LS-path conditions are equivalent to the lattice conditions for $\mathrm{LS}_{\mathfrak C,\lambda}$. We have
for all $i=1,\ldots,r$:
\begin{equation}\label{equivLatticePath}
\hat d_j\langle \kappa_j(m\lambda),\beta^\vee_j\rangle=m\hat d_j\langle \kappa_j(\lambda),\beta^\vee_j\rangle\in\mathbb Z
\Leftrightarrow (a_r+\ldots+a_j)b_j\in \mathbb Z.
\end{equation}
It follows that $\Theta$ is indeed a well defined injective map from $\mathbb B(m\lambda)_\tau$ to 
$\mathrm{LS}_{\lambda}^+(m)\subseteq \mathbb Q^{A_\tau}$. 

The arguments used above work also vice versa: starting with an element in $\underline{a}\in \mathrm{LS}_{\lambda}^+(m)$,
fix a maximal chain such that $\underline{a}\in \mathrm{LS}_{\mathfrak C,\lambda}^+(m)$.
One attaches to $\underline{a}$ an ``extended'' LS-path $\hat\pi$ having as sequence of Weyl group elements
the elements of the maximal chain, and the rational numbers $\hat d_i$ are given by $\hat d_i=(a_r+\ldots+a_i)/m$, $i=0,\ldots,r$.
In addition, set $\hat d_{r+1}=0$.
Let $\pi$ be the pair of sequences obtained by omitting those $\kappa_{j+1}$ and $\hat d_{j+1}$ such that 
$\hat d_{j+1}=\hat d_{j}$, $j=0,\ldots,r$. The equivalence in \eqref{equivLatticePath} shows: 
this is a LS-path $\pi$ of shape $m\lambda$ contained in $\mathbb B(m\lambda)_\tau$. 

It follows: the map $\Theta$ induces a bijection between the LS-paths $\pi$ in $\mathbb B(m\lambda)_{\tau}$ 
and the elements in $\mathrm{LS}^+_{\lambda}(m)$.
\end{proof}
The set $\mathrm{LS}^+_{\lambda}$ is in general a fan of monoids, but a single monoid. The sum of two elements
$\underline a^1,\underline a^2 \in \mathrm{LS}^+_{\lambda}$ makes sense (meaning is again an 
element in $\mathrm{LS}^+_{\lambda}$) only if there exists a maximal chain $\mathfrak C$ in $A_\tau$
such that $\supp \underline a^1,\supp \underline a^2\subseteq \mathfrak C$. This leads to the following
definition:
\begin{definition}\label{AppendixI_standardsum}
A \textit{sum} $\underline a=\underline a^1+\ldots+\underline a^m$ of elements $\underline a^1,\ldots,\underline a^m\in 
\mathrm{LS}^+_{\lambda}$ is called \textit{standard} if $\supp \underline a^1\ge \supp \underline a^2\ge \ldots\ge \supp \underline a^m$.
In particular, a standard sum is again an element in $\mathrm{LS}^+_{\lambda}$.
\end{definition}
The map $\Theta$ respects the support, so it turns a concatenation $\pi_1*\cdots*\pi_m$, which is standard, 
into a sum which is standard: $\Theta(\pi_1*\cdots*\pi_m)=\sum_{i=1}^m \Theta(\pi_i)$.
An important consequence: by Proposition~\ref{decompLSpath}, 
an LS-path $\pi$ of shape $m\lambda$ can always be written in a unique way as a concatenation
of LS-paths of shape $\lambda$: $\pi=\pi_1*\cdots*\pi_m$ such that $\supp \pi_1\ge \supp\pi_2\ge \ldots\ge \supp \pi_m$.
For an element in $\mathrm{LS}^+_{\lambda}$ this implies: 

\begin{lemma}\label{latticedecomp}
Every element $\underline a\in  \mathrm{LS}^+_{\lambda}(m)$
has a unique decomposition $\underline a=\underline a^1+\ldots+\underline a^m$ into $m$ elements
$\underline a^i\in \mathrm{LS}^+_{\lambda}(1)$, $i=1,\ldots,m$, such that 
$\supp \underline a^1\ge \supp \underline a^2\ge \ldots\ge \supp \underline a^m$.
\end{lemma}
See also \cite{Ch} for a different proof of Lemma~\ref{latticedecomp}.

\section{Appendix II: A filtration on 
$V_{\mathbb Z}(s\lambda)$}\label{A_filtration_on_V}
\subsection{Some comments on the Appendix II and III}{  The reason to add this and the following section as { appendices} is the slightly different approach compared to the one in \cite{L1}. The main goal in \textit{ibidem} was to show that a standard monomial theory in the sense of Lakshmibai, Musili and Seshadri exists for all Schubert varieties in a generalized flag variety, and this in the setting of a symmetrizable Kac-Moody group. So it was sufficient to construct \textit{one} standard monomial theory. For this the path vectors $p_\pi$ in $V(\lambda)^*$ have been constructed, and to prove their linear independence, a basis $\{v_\pi\mid\textrm{$\pi$ LS path}\}$ of $V_{\mathbb Z}(\lambda)^*$
was constructed. Both bases depend on certain choices, which, at that time, was not relevant to achieve the main goal.

To show that the construction of a standard monomial theory via a Seshadri stratification is compatible with the construction in \cite{L1}, one has to take care of these possible choices. For this reason we introduce the filtration of $V_{\mathbb Z}(\lambda)_\tau$
in \eqref{appendix_two_filtration_one} respectively 
\eqref{appendix_two_filtration_two}, and its dual filtration 
in \eqref{appendix_drei_filtration_one} respectively 
\eqref{appendix_drei_filtration_two}. Note that the definition of a path vector in Definition~\ref{define_path_vector2} (respectively Definition~\ref{define_path_vector}) is more general than the corresponding definition in \cite{L1}. Indeed, Corollary~\ref{differentEll} shows how these two definitions are related. Theorem~\ref{some_vanishingA} shows that the path vectors are representatives of the filtration associated to the Seshadri stratification, { independently} of the choice of a compatible total order. The main ingredient in the proof is Corollary~\ref{factoring}, which enables us to proceed with the proof by induction. 

Many of the tools and some of the rather technichal statements 
in this and the following section can be, at least implicitly, found in \cite{L1}. But the more general notion 
of a path vector and the idea to use filtrations instead of a fixed ordered basis made it necessary to adapt some of the proofs
and definitions.
For example, the preorder defined in Definition~\ref{preorderOnTensorProductOfWeight} is slightly different from the one used in \cite{L1}.

{  In order to} present a text written in a uniform style 
(and not always refer to \cite{L1} with some explanation 
why the results still hold in this slightly different setting)
we have decided to rewrite part of \cite{L1}. We have declared this part as { appendices} because apart from minor
changes and the adaption to the new setting, { none} of the results is new except for the following statements: part \textit{iii)} of Theorem~\ref{AuxBasisTheorem}, as well as 
Proposition~\ref{product_path_vector} and the statements after. They are new because they are now formulated for path vectors in the sense of Definition~\ref{define_path_vector2}.

To conclude, one last observation. The filtration defined in 
\eqref{section_six_filt_one} and \eqref{section_six_filt_two}
(see also \eqref{appendix_two_filtration_one} and \eqref{appendix_two_filtration_two})
is by construction of combinatorial nature, and so
is its dual version defined in \eqref{section_six_dual_filtr_one}
and \eqref{section_six_dual_filtr_two}. 
{ 
These filtrations are strongly related to the definition of the path vectors. The role of the quantum Frobenius
splitting is now reduced to that of a tool proving the linear independence of the vectors $v_{\underline{a},\underline{\sigma}}$.
So it would be interesting to have a
purely representation theoretic interpretation of the filtration
defined in \eqref{section_six_filt_one} and \eqref{section_six_filt_two}.}
}
\subsection{Notation} We start with the complex version of the symmetrizable Kac-Moody algebra $\mathfrak g_{\mathbb C}$.
Let $U(\mathfrak g)_{\mathbb Z}$ be the Kostant-$\mathbb Z$-form of the enveloping algebra $U(\mathfrak g)$, and 
for a ring $\hat R$ set $U(\mathfrak g)_{\hat R}=U(\mathfrak g)_{\mathbb Z}\otimes_{\mathbb Z}\hat R $. 
\index{$U(\mathfrak g)_{\hat R}$, form over a ring $\hat R$}

For a simple root $\alpha_i$, $1\le i\le n$, denote by $X_i$ the Chevalley generator of weight $\alpha_i$ and let 
$X_{-i}$ be the Chevalley generator of weight $-\alpha_i$. We denote by $U(\mathfrak g)^+_{\mathbb Z}$ respectively 
$U(\mathfrak g)^-_{\mathbb Z}$ the part of the Kostant-$\mathbb Z$-form generated by the divided powers of the $X_i$
respectively the $X_{-i}$. For a ring $\hat R$ we set $U(\mathfrak g)^+_{\hat R}=U(\mathfrak g)^+_{\mathbb Z}\otimes_{\mathbb Z}\hat R$,
$U(\mathfrak g)^-_{\hat R}$ is defined similarly.\index{$U(\mathfrak g)^\pm_{\hat R}$, form over a ring $\hat R$}

We fix a dominant integral weight $\lambda$ and choose a highest weight vector $v_\lambda$ in the complex irreducible 
$\mathfrak g_{\mathbb C}$-module $V_{\mathbb C}(\lambda)$. Let $V_{\mathbb Z}(\lam) = U(\mathfrak g)_{\mathbb Z}.v_{\lambda}$
be the corresponding Kostant-$\mathbb Z$-form of the module. We associate to $\tau\in W/W_Q$ in a canonical way
an extremal weight vector $v_{\tau}$ of weight $\tau(\lambda)$ as follows: let $\tau=s_{\alpha_{i_1}}\cdots s_{\alpha_{i_r}}$
be a reduced decomposition and set $n_r=\langle \lambda,\alpha_{i_r}^\vee\rangle$,  $n_{r-1}=\langle s_{\alpha_{i_r}}(\lambda),\alpha_{i_{r-1}}^\vee\rangle$,
$\ldots$, and $n_{1}=\langle (s_{\alpha_{i_2}}\cdots s_{\alpha_{i_r}})(\lambda),\alpha_{i_1}^\vee\rangle$. The numbers are all positive,
and we set
\begin{equation}\label{extremalweightvector}
v_\tau= X_{-\alpha_{i_1}}^{(n_1)}\cdots X_{-\alpha_{i_{r-1}}}^{(n_{r-1})}X_{-\alpha_{i_r}}^{(n_r)} v_\lambda.
\end{equation}
Here $X^{(n)}$ is the usual abbreviation for the divided power $\frac{X^n}{n!}$. The fact that $v_\tau$ is independent of the chosen
reduced decomposition follows from the Verma identities, see \cite{Verma} respectively \cite{Lu}, Section 39.3: if $\tau=s_{\alpha_{i_1}}\cdots s_{\alpha_{i_r}}=
s_{\alpha_{j_1}}\cdots s_{\alpha_{j_r}}$
are two reduced decompositions and $(n_1,\ldots,n_r)$ respectively $(m_1,\ldots,m_r)$ are the associated tupels of integers, then one has
\begin{equation}\label{vermaId}
X_{-\alpha_{i_1}}^{(n_1)}\cdots X_{-\alpha_{i_{r-1}}}^{(n_{r-1})}X_{-\alpha_{i_r}}^{(n_r)}=
X_{-\alpha_{j_1}}^{(m_1)}\cdots X_{-\alpha_{j_{r-1}}}^{(m_{r-1})}X_{-\alpha_{j_r}}^{(m_r)}
\end{equation}
in the enveloping algebra of $\mathfrak g$.
The weight vector  $v_{\tau}$ generates the weight space $V_{\mathbb Z}(\lam)_{\tau(\lambda)}$, which is free of rank one.

We denote by $V_{\mathbb Z}(\lam)_{\tau}$ the corresponding Demazure module 
$V_{\mathbb Z}(\lam)_{\tau}=U(\mathfrak g)^+_{\mathbb Z}v_\tau$.
For a ring $\hat R$ let $V_{\hat R}(\lambda)=V_{\mathbb Z}(\lambda)\otimes_{\mathbb Z} \hat R$ be
the corresponding module for $U(\mathfrak g)_{\hat R}$ over the ring $\hat R$. Similarly, for a ring $\hat R$ we write 
$V_{\hat R}(\lambda)_\tau=V_{\mathbb Z}(\lambda)_\tau\otimes_{\mathbb Z} \hat R$ for the Demazure module over the ring $\hat R$.
\index{$V_{\hat R}(\lambda)_\tau$, form of the Demazure module over a ring $\hat R$}
\subsection{Some combinatorics}\label{algorithm}
A first step in the construction of the filtration of $V_{\mathbb Z}(s\lambda)_\tau$
is a combinatorial procedure which associates  to a given element $\underline{a}\in\mathrm{LS}_{\lambda}^+(s)$
a sequence of pairs of simple roots and nonnegative integers.

Let $\supp\underline{a}=\{\tau_q,\tau_{q-1},\ldots,\tau_1\}$ be the support of $\underline{a}$ with $\tau_q>\tau_{q-1}>\ldots>\tau_1$. We fix a reduced decomposition $\tau_q=s_{i_1}\cdots s_{i_t}$.  
We have $s_{i_1}\tau_q<\tau_q$. Let now $0\le j<q$ be minimal such that $s_{i_1}\tau_i\le \tau_i$ for all $i > j$.
Or, in other words, $j$ is maximal such that $s_{i_1}\tau_j > \tau_j$, and one knows hence: $s_{i_1}\tau_{j+1}\ge \tau_{j}$.

Let $\underline{a}'\in\mathbb Q^{A_\tau}$ be the element obtained from $\underline{a}$ as follows:
\begin{equation}\label{aprim}
\underline{a}'=\left\{\begin{array}{ll}
\sum_{h=j+1}^q a_{\tau_h}e_{s_{i_1}\tau_{h}} +\sum_{h=1}^{j} a_{\tau_h}e_{\tau_h},&\textrm{if\ } s_{i_1}\tau_{j+1}\not=\tau_j;\\
\sum_{h=j+2}^q a_{\tau_h}e_{s_{i_1}\tau_{h}} + (a_{\tau_{j+1}}+a_{\tau_{j}})e_{\tau_{j}} +
\sum_{h=1}^{j-1} a_{\tau_h}e_{\tau_h},& \textrm{if\ } s_{i_1}\tau_{j+1}=\tau_j.
\end{array}\right.
\end{equation}
\begin{lemma}\label{chainForAprime}
There exists a maximal chain $\mathfrak C'$ such that $\underline{a}'\in \mathrm{LS}_{\lambda,\mathfrak C'}^+(s)$.
\end{lemma}
\begin{proof}
Let $\alpha$ be a simple root and let $\sigma\in W/W_Q$ be such that $s_\alpha\sigma < \sigma$.
If $\kappa<\sigma$ is of length one less than $\sigma$, then either $s_\alpha\sigma =\kappa$,
or $s_\alpha\kappa < \kappa$. And, in the latter case, one has $s_\alpha\kappa<s_\alpha\sigma$,
 and the bond between $\kappa$ and $\sigma$
is the same as the bond between $s_\alpha\kappa$ and $s_\alpha\sigma$. 
Starting with a maximal chain $\mathfrak C$ such that $\underline{a}\in \mathrm{LS}_{\lambda,\mathfrak C}^+(s)$,
one uses the procedure just described to construct stepwise out of $\mathfrak C$  a maximal chain $\mathfrak C'$ such that
 $\underline{a'}\in \mathrm{LS}_{\lambda,\mathfrak C'}^+(s)$.
\end{proof}
\begin{remark}
A proof of Lemma~\ref{chainForAprime} in the language of paths can be found in \cite{L1}.
\end{remark}
By the bijection in Proposition~\ref{LSlatticeLSpath} between LS-paths of shape a multiple of $\lambda$
and elements of the fan of monoids $\mathrm{LS}_\lambda^+$, we can associate to
$\underline{a}\in \mathrm{LS}_{\lambda}^+(s)$ an element in the weight lattice. We set:
$$
\textrm{\rm weight}(\underline{a})=\sum_{i=1}^q a_i\tau_i(\lambda)\in\Lambda.
$$
It is the weight of the LS-path $\Theta^{-1}(\underline{a})$ (see \eqref{weightOfpaths}) of shape $s\lambda$ associated to 
$\underline{a}$ via the bijection in  \eqref{pathLattice}.
In the same way we can associate an element in the weight lattice to $\underline{a}'$ (see \eqref{aprim}), 
and, by construction, we have 
\begin{equation}\label{weightdifference}
\textrm{\rm weight}(\underline{a'})-\textrm{\rm weight}(\underline{a})=n_1\alpha_{i_1}=
(a_q\vert\langle\tau_q(\lambda),\alpha_{i_1}^\vee\rangle\vert+\ldots+
a_{j+1}\vert\langle\tau_{j+1}(\lambda),\alpha_{i_1}^\vee\rangle\vert) \alpha_{i_1}.
\end{equation}
So by the procedure above we associate to $\underline{a}\in\mathrm{LS}_{\lambda}^+(s)$ 
the pair $((n_1,\alpha_{i_1}), \underline{a}')$, where $n_1\in\mathbb N$, $\alpha_{i_1}$ is a simple root
and $\underline{a}'\in\mathrm{LS}_{\lambda}^+(s)$.

Note that $s_{i_1}\tau_q=s_{i_2}\cdots s_{i_t}$ is a reduced decomposition of $s_{i_1}\tau_q$, the maximal 
element in $\supp\underline{a}'$. So by repeating
the procedure with $\underline{a}'$, we associate to $\underline{a}$ the tuple 
$((n_1,\alpha_{i_1}),(n_2,\alpha_{i_2}), \underline{a}'')$, where $n_1,n_2\in\mathbb N$, $\alpha_{i_1},\alpha_{i_2}$ are simple roots
and $\underline{a}''\in\mathrm{LS}_{\lambda}^+(s)$.

By repeating the procedure again and again, we arrive at a sequence of integers and simple roots:
$((n_1,\alpha_{i_1}),(n_2,\alpha_{i_2}), \ldots, (n_t,\alpha_{i_t}),\underline{a}''')$, where $\underline{a}''' \in \mathrm{LS}_{\lambda}^+(s)$
is equal to $se_{{\rm id}}$. The simple roots in the sequence are determined by the reduced decomposition $\tau_q=s_{i_1}\cdots s_{i_t}$ 
and the integers are determined by $\underline{a}$ through the procedure explained above.

The procedure associates hence to $\underline{a}$ and the chosen reduced decomposition $\underline{\tau}_q$ of $\tau_q$ a 
sequence of pairs of simple roots and nonnegative integers: 
\begin{equation}\label{sequenceofrootsandnumbers}
S(\underline{a},\underline{\tau}_q):=\left((n_1,\alpha_{i_1}),(n_2,\alpha_{i_2}),\ldots,(n_t,\alpha_{i_t})\right).
\index{$S(\underline{a},\underline{\sigma})$, sequence of simple roots and numbers}
\end{equation} 
\subsection{Some vectors and a filtration}\label{filtrationspaces}
Given $\underline{a}\in\mathrm{LS}_{\lambda}^+(s)$, let 
$\sigma\in\supp\underline{a}$ be the maximal element.
We fix a reduced decomposition $\underline{\sigma}$, i.e. $\sigma=s_{i_1}\cdots s_{i_t}$,
and let $S(\underline{a},\underline{\sigma})=\left((n_1,\alpha_{i_1}),(n_2,\alpha_{i_2}),\ldots,(n_t,\alpha_{i_t})\right)$ 
be the associated sequence of simple roots and integers 
as in \eqref{sequenceofrootsandnumbers}.
\begin{definition}\label{definitionvatau}
The vector $v_{\underline{a},\underline{\sigma}}\in V_{\mathbb Z}(s\lambda)_\tau$ 
\textit{associated to $\underline{a}\in \mathrm{LS}_{\lambda}^+(s)$ 
and the reduced decomposition $\underline{\sigma}$} is defined by:
$$
v_{\underline{a},\underline{\sigma}}=X_{-i_1}^{(n_1)}\cdots X_{-i_t}^{(n_t)} v_{s\lambda}.
\index{$v_{\underline{a},\underline{\sigma}}$, vector associated to $\underline{a}$}
$$
\begin{remark}
One may view this procedure as a generalization of the construction of extremal weight vectors
in \eqref{extremalweightvector}. Indeed, it is easy to check that if $\underline{a}\in \mathrm{LS}_{\lambda}^+(1)$ is of weight
$\sigma(\lambda)$ for some $\sigma\in W/W_{}$, then $v_{\underline{a},\underline{\sigma}}=v_\sigma$, independent
of the chosen reduced decomposition. But in general, note that $v_{\underline{a},\underline{\sigma}}$ depends
on the choice of the reduced decomposition $\underline{\sigma}$ of $\sigma$, 
see Remark~\ref{dependenceReducedDecomposition}. 
The procedure to get the sequence  $S(\underline{a},\underline{\sigma})$ was inspired by \cite{RS}.
K. N. Raghavan and P. Sankaran define in this article a procedure to attach such a 
sequence to the Young tableaux defined by Lakshmibai, Musili and Seshadri in \cite{LSIV},\cite{LS}.
These tableaux exist, for example, for all representations of the simple groups of classical type {\tt A,B,C,D}.
\end{remark}
\begin{remark}
Let $\tau=s_{j_1}\cdots s_{j_r}$ be a reduced decomposition.
It is well known that $V_{\mathbb Z}(s\lambda)_\tau$ is
spanned by the vectors of the form
$X_{-j_1}^{(p_1)}\cdots X_{-j_t}^{(p_r)} v_{s\lambda}$. By choosing for $\sigma\le \tau$ a compatible
reduced decomposition, one sees easily that $V_{\mathbb Z}(s\lambda)_\sigma \subseteq V_{\mathbb Z}(s\lambda)_\tau$.
It follows hence: $v_{\underline{a},\underline{\sigma}}\in V_{\mathbb Z}(s\lambda)_\tau$.
\end{remark}

We denote by $V_{\mathbb Z}(s\lambda)_{\tau,\underline{a}}$ the $\mathbb Z$-submodule:\index{$V_{\mathbb Z}(s\lambda)_{\tau,\underline{a}}$,
submodule in the filtration}
\begin{equation}\label{appendix_two_filtration_one}
V_{\mathbb Z}(s\lambda)_{\tau,\unlhd\underline{a}}=\left\langle 
v_{\underline{a}',\underline{\sigma}'}\left\vert 
\begin{array}{l}\underline{a}\unrhd \underline{a}',\underline{a}'\in \mathrm{LS}_{\lambda}^+(s),
\ \underline{\sigma}' \textit{\ reduced decomposition}\\
\textit{\ of the maximal element in\ }\supp\underline{a}'
\end{array}\right.\right\rangle,
\end{equation}

and by $V_{\mathbb Z}(s\lambda)_{\tau,\lhd\underline{a}}$ the $\mathbb Z$-submodule:
\begin{equation}\label{appendix_two_filtration_two}
V_{\mathbb Z}(s\lambda)_{\tau,\lhd\underline{a}}=\left\langle 
v_{\underline{a}',\underline{\sigma}'}\left\vert 
\begin{array}{l}\underline{a}\rhd \underline{a}',\underline{a}'\in \mathrm{LS}_{\lambda}^+(s),
\ \underline{\sigma}' \textit{\ reduced decomposition}\\
\textit{\ of the maximal element in\ }\supp\underline{a}'
\end{array}\right.\right\rangle.
%
\end{equation}
\end{definition}
We refer to  Definition~\ref{newpartialorder} for the definition of ``$\unrhd$''.

\begin{remark}\label{dependenceReducedDecomposition}
In general, $v_{\underline{a},\underline{\sigma}}$ depends on the chosen reduced decomposition 
$\underline{\sigma}$ of $\sigma$. So in $V_{\mathbb Z}(s\lambda)_{\tau,\underline{a}}$ one may find several
vectors $v_{\underline{a},\underline{\sigma}}$ and $v_{\underline{a},\underline{\sigma}'}$ indexed by the same
$\underline{a}\in \mathrm{LS}_{\lambda}^+(s)$, but with different reduced decompositions $\underline{\sigma}$ and 
$\underline{\sigma}'$ of the maximal element $\sigma\in\supp \underline{a}$. The next aim is to analyze the difference
$v_{\underline{a},\underline{\sigma}}-v_{\underline{a},\underline{\sigma}'}$.
\end{remark}
\subsection{The Frobenius splitting trick}\label{FrobeniusTrick}
To have a closer look at the $\mathbb Z$-submodules $V_{\mathbb Z}(s\lambda)_{\tau,\underline{a}}$ defined above,
we use representation theory of quantum groups at a root of unity. To keep the notation simple, we assume 
$\mathfrak g$ to be of type {\tt A}, {\tt D} or {\tt E}, so the associated
Cartan matrix is symmetric. The construction works  in the same way
for an arbitrary symmetrizable Kac-Moody algebra, only the details are more technical. In particular,
the Frobenius splitting involves the Langlands dual Kac-Moody algebra. The details can be found in \cite{L1}.

Let $\ell$ be an even positive integer. Let $A:=\Z[q,q^{-1}]$ be the ring of Laurent polynomials
and denote by $R$ the ring $A/I$, where $I$ is the ideal
generated by the $2\ell$-th cyclotomic polynomial.  Let ${\tilde R}$ be the ring obtained by adjoining
all roots of unity to $\Z$. We fix an embedding $R\hookrightarrow{\tilde R}$.
If $\mathbb K$ is an algebraically closed field and $\charc\, \mathbb K=0$, then we consider $\mathbb K$
as an $R$- respectively $\tR$-module by the inclusions $R\subseteq {\tilde R}\subseteq \mathbb K$.
If $\charc\, \mathbb K=p>0$, then we consider $\mathbb K$ as an $R$- respectively $\tR$-module by extending the
canonical map $\mathbb Z\rightarrow \mathbb K$ to a map $R\hookrightarrow {\tilde R}\rightarrow \mathbb K$ (where
the first map is given by the projection $\mathbb Z\rightarrow\mathbb Z/p\mathbb Z$ and the
inclusion $\mathbb Z/p\mathbb Z\subset \mathbb K$). 

{  
Similar to the classical case (see \cite{KL1, KL2, L1, Lu} for more details), let
$U_{q}(\mathfrak g)$ be the quantum group associated to { the symmetrizable Kac-Moody algebra $\mathfrak{g}$} over the field $\mathbb C(q)$. Denote by 
$U_{q, A}(\mathfrak{g})$
the corresponding Kostant-Lusztig form of the quantum group \cite{Lu}
and let $U_{\xi}(\mathfrak g)_\tR:= { U_{q,A}(\mathfrak{g})\otimes \tR}$ be the 
corresponding quantum group over $\tR$ at a $2\ell$-th root of unity $\xi$.

For a simple root $\alpha_i$ denote by $E_{\pm i}$ the generators of the quantum group { $U_{\xi}(\mathfrak g)_{\tR}$} 
corresponding to the root $\alpha_i$ respectively $-\alpha_i$.
It has been shown by Lusztig (see \cite{Lu}) that the assignment 
$$
X^{(k)}_i\mapsto E_{i}^{(k\ell)},
$$
extends to an algebra homomorphism $\mathrm{Fr}':U(\mathfrak g)^+_\tR\rightarrow U_{\xi}(\mathfrak g)^+_\tR$, 
called the splitting of the Frobenius map. Here $U_{\xi}(\mathfrak g)^\pm_\tR$ 
denotes the part of the Lusztig form generated by the divided powers of the 
$E_{\pm i}$. There is a corresponding Frobenius splitting map $\mathrm{Fr}':U(\mathfrak g)^-_\tR\rightarrow U_{\xi}(\mathfrak g)^-_\tR$
sending $X^{(k)}_{-i}\mapsto E_{-i}^{(k\ell)}$.
\begin{remark}
The reason to restrict to quantum groups at a $2\ell$-th root of unity for an even $\ell$ is just for convenience.
A splitting also exists { for} other roots of unity, but one has to be careful because there are restrictions for certain
types. For our considerations the case $2\ell$-th root of unity for an even $\ell$  is sufficient,
and the assumptions made here makes it possible to present a uniform approach. 
\end{remark}
For a dominant weight $\lam\in \Lambda^+$ let $M(\lambda)$ 
be the simple $U_q(\mathfrak{g})$-module of highest weight $\lambda$
over the field $\mathbb C(q)$. 
We fix an $A$-lattice in $M(\lambda)$ by choosing a highest 
weight vector $m_{\lambda}$, and we set: 
$M_{A}(\lambda)=U_{q, A}(\mathfrak{g})m_\lambda$.
Set $M_{R}(\lambda):= M_A(\lambda)\otimes_A \tR$, 
then $M_{R}(\lambda)$ is a $U_{\xi}(\mathfrak g)_\tR$-module 
such that its character is given by the Weyl character formula. 
It is called the Weyl module of highest weight $\lambda$ for  
$U_\xi({\mathfrak g})_\tR$, see \cite{APW,KL1, KL2} for details and further properties.}

The extremal weight spaces are free of rank one.
Fix a highest weight vector $m_\lambda$ and, using the quantum Verma identities \cite{Lu}, one can associate 
to $\tau\in W/W_Q$ a corresponding extremal weight vector $m_{\tau}\in M_\tR(\lam)$
of weight $\tau(\lambda)$, it generates the corresponding rank one $\tR$-module.
Denote by $M_\tR(\lambda)_\tau:=U_\xi({\mathfrak g})_\tR^+\cdot m_{\tau}$
the Demazure submodule associated to $\tau$. Via the Frobenius splitting $\mathrm{Fr}'$, the  $U_\xi({\mathfrak g})_\tR^+$-module
$M_\tR(\lambda)$ and its tensor powers become $U({\mathfrak g})^+_\tR$- respectively  $U({\mathfrak g})^-_\tR$-modules.
In general, one can not glue these actions together to get a global $U({\mathfrak g})_\tR$-representation.

The module $M_\tR(\lambda)_{}^{\otimes\ell}$ has a weight space decompostion $\bigoplus_{\mu\in \Lambda} \big(M_\tR(\lambda)_{}^{\otimes\ell}\big)_{\mu}$.
Denote by $(M_\tR(\lambda)_{}^{\otimes\ell})^{1/\ell}$ the $\tR$- submodule $\bigoplus_{\mu\in \ell\Lambda} (M_\tR(\lambda)_{}^{\otimes\ell})_\mu$.
The $\tR$-submodule $(M_\tR(\lambda)_{}^{\otimes\ell})^{1/\ell}$ 
is not a  $U_\xi({\mathfrak g})_\tR$-submodule, but it is, via the Frobenius splitting map,  a $U({\mathfrak g})^+_\tR$- respectively  
$U({\mathfrak g})^-_\tR$-submodule. It has been shown in \cite{L1}: on this submodule one can glue these two module structures to get a 
$U({\mathfrak g})_\tR$-module structure.

In particular, we have a natural inclusion of $U({\mathfrak g})^+_\tR$- respectively of  $U({\mathfrak g})^-_\tR$-modules:
\begin{equation}\label{quantuminj1}
V_\tR(\lambda)_{}\hookrightarrow (M_\tR(\lambda)_{}^{\otimes\ell})^{1/\ell}\subseteq M_\tR(\lambda)_{}^{\otimes\ell},
\end{equation}
which sends the highest weight vector $v_{\lambda}$ to the highest weight vector $m_\lambda^{\otimes\ell}$,
and, more generally, an extremal weight vector $v_{\kappa}\in V_\tR(\lambda)_{\tau}$ to $m_{\kappa}^{\otimes\ell}$.
More generally, we consider natural inclusions of $U({\mathfrak g})^+_\tR$- respectively of  $U({\mathfrak g})^-_\tR$-modules:
$$
V_{\tR}(s\lambda)\hookrightarrow  ( (M_\tR(s\lambda)^{\otimes \ell})^{\frac{1}{\ell}})\subseteq M_\tR(s\lambda)^{\otimes \ell}
\hookrightarrow M_\tR(\lambda)^{\otimes s\ell}.
$$
In the comultiplication formula $\Delta(E_{\pm j}^{(m)})$ for the quantum group \cite{Lu}, the generators $K_j^{\pm}$ associated to the simple
roots $\alpha_j$ play an important role:
\begin{equation}\label{comultiplication_rule}
\begin{array}{rcl}
\Delta(E_{j}^{(p)})&=& \sum_{p'+p''=p}
q^{p' p''}E_{j}^{(p')}K_{j}^{p''}\otimes E_{j}^{(p'')}\\
\Delta(E_{-j}^{(p)})&=& \sum_{p'+p''=p}
q^{-p' p''}E_{-j}^{(p')}\otimes K_{-j}^{p'} E_{-j}^{(p'')}.
\end{array}
\end{equation}
Recall that we consider quantum groups at a root of unity, so $q=\xi$ is a $2\ell$-th root of unity. For the representations we consider,
the generators $K_{\pm j}$ operate on weight spaces by multiplication with some root of unity.
As long as the precise value is not important, we use the notation 
$\v^{\bullet}$ \index{$\v^{\bullet}$, some root of unity}  for ``some root of unity''.

\subsection{The vectors $v_{\underline{a}, \underline{\sigma}}$ and their image in $M_{\tR}(\lambda)^{\otimes s\ell}$}
Given $\underline{a}\in\mathrm{LS}_{\lambda}^+(s)$, let 
$\sigma\in\supp\underline{a}$ be the maximal element.
We fix a reduced decomposition $\underline{\sigma}$, i.e. $\sigma=s_{i_1}\cdots s_{i_t}$.
As a next step we consider the image of $v_{\underline{a}, \underline{\sigma}}$ in $M_{\tR}(\lambda)^{\otimes s\ell}$
with respect to the embedding (for an appropriate $\ell$):
$$
V_{\tR}(s\lambda)\hookrightarrow  ( M_\tR(s\lambda)^{\otimes \ell})^{\frac{1}{\ell}}\hookrightarrow
 M_\tR(s\lambda)^{\otimes \ell} \hookrightarrow M_\tR(\lambda)^{\otimes s\ell}.
$$
These are embeddings of $U(\mathfrak g)^\pm$-modules, acting on the quantum group modules via $\mathrm{Fr}'$. We will use the notation $\circ$ to denote these module structures. For example, for $m\in M_{\tilde{R}}(\lambda)^{\otimes s\ell}$ and $X\in U(\mathfrak g)^\pm$, $X\circ m:=\mathrm{Fr}'(X)(m)$.

We will see that the image of $v_{\underline{a},\underline{\sigma}}$ (see Definition~\ref{definitionvatau}) has a natural ``maximal term'' which is independent of the  choice
of the reduced decomposition $\underline{\sigma}$. To make this statement more precise, we define a preorder on tuples
of weight vectors in $M_{\tR}(\lambda)$.

Let ``$\succ$'' be the usual partial order on the set of weights: $\nu\succ\chi$ if the difference $\nu-\chi$ is a non-negative
sum of positive roots. In the following we define a preorder on the set of ordered tuples of non-zero weight vectors.
Since rescaling does not affect the preorder, we write the tuples as tensor products. 
\begin{definition}\label{preorderOnTensorProductOfWeight}
Let $m_{\nu_1},\ldots, m_{\nu_t}$
and $m_{\chi_1},\ldots, m_{\chi_t}$ be weight vectors  in $M_{\tR}(\lambda)$ of weight $\nu_1,\ldots, \nu_t$ respectively $\chi_1,\ldots, \chi_t$.
We say $m_{\nu_1}\otimes \cdots\otimes m_{\nu_t}$ is \textit{smaller} than $m_{\chi_1}\otimes \cdots\otimes m_{\chi_t}$, in symbols:
$$
m_{\nu_1}\otimes \cdots\otimes m_{\nu_t} \LHD\index{$\LHD$, preorder, defined on tensor product of weight vectors} m_{\chi_1}\otimes \cdots\otimes m_{\chi_t}, 
$$
\begin{itemize}
\item[-]if $\nu_1\succ\chi_1$ and $m_{\nu_1}\in U_{\xi}(\mathfrak g)^+_\tR\cdot m_{\chi_1}$, 
\item[-]or $\nu_1=\chi_1$,  $m_{\nu_1}, m_{\chi_1}$ are linearly dependent, $\nu_2\succ\chi_2$ and $m_{\nu_2}\in U_{\xi}(\mathfrak g)^+_\tR\cdot m_{\chi_2}$, 
\item[-]or $\nu_1=\chi_1$, $\nu_2=\chi_2$,  $m_{\nu_1}, m_{\chi_1}$ are linearly dependent,  $m_{\nu_2}, m_{\chi_2}$ are linearly dependent,
$\nu_3\succ\chi_3$ and $m_{\nu_3}\in U_{\xi}(\mathfrak g)^+_\tR\cdot m_{\chi_3}$, 
\par\noindent
\item[-]and so on...
\end{itemize}
\end{definition}

Note that $\succ$ turns into $\LHD$. It is easy to see that for $\underline{a},\underline{a}'\in \mathrm{LS}_{\lambda}^+(s)$: 
$m^{\underline{a}}\RHD m^{\underline{a}'}$ if and only if $\underline{a}'\rhd \underline{a}$.
We use sometimes the following variation of the preorder. 
\begin{definition}
Given weight vectors
$m_{\nu_1},\ldots, m_{\nu_t}$, $m_{\chi_1},\ldots, m_{\chi_t}$ in $M_{\tR}(\lambda)$ as above, 
and a permutation $\mathfrak s\in \mathfrak S_t$, we write
$$
m_{\nu_1}\otimes \cdots\otimes m_{\nu_t} \LHD_{\mathfrak s} m_{\chi_1}\otimes \cdots\otimes m_{\chi_t}. 
$$
\begin{itemize}
\item[-]if $\nu_{\mathfrak s(1)}\succ \chi_{\mathfrak s(1)}$ and $m_{\nu_{\mathfrak s(1)}}\in U_{\xi}(\mathfrak g)^+_\tR\cdot m_{\chi_{\mathfrak s(1)}}$, 
\item[-]or $\nu_{\mathfrak s(1)}=\chi_{\mathfrak s(1)}$, $m_{\nu_{\mathfrak s(1)}}$ and $m_{\chi_{\mathfrak s(1)}}$ are linearly dependent, and
$\nu_{\mathfrak s(2)}\succ\chi_{\mathfrak s(2)}$ and the weight vector $m_{\nu_{\mathfrak s(2)}}\in U_{\xi}(\mathfrak g)^+_\tR\cdot m_{\chi_{\mathfrak s(2)}}$,
\par\noindent
\item[-]and so on...
\end{itemize}
\end{definition}

Let $\mathfrak C=(\tau_r,\ldots,\tau_0)$ be a maximal chain in $A_\tau$ and let
$\underline{a}\in \mathrm{LS}_{\mathfrak C,\lambda}^+(s)$. Let $\ell$ be an even number such that 
$\ell a_{\tau_j}\in\mathbb N$ for all $j=0,\ldots,r$. Let $\sigma=\tau_q\in\supp{\underline{a}}$ be the unique maximal element.
We associate to $\underline{a}$ the vector
$$
m^{\underline{a}}:=
\underbrace{ m_{\tau_q}\otimes\ldots\otimes m_{\tau_q} }_{\ell a_q}
\otimes
\underbrace{
m_{\tau_{q-1}}\otimes\ldots\otimes 
m_{\tau_{q-1}}
}_{\ell a_{q-1}}
\otimes \ldots\otimes
\underbrace{
m_{\tau_0}\otimes\ldots\otimes m_{\tau_0}}_{\ell a_{0}}\in
(M_{\tR}(\lam))^{\otimes{s\ell}}.
$$
Sometimes it is convenient to permute the tensor factors. If $\mathfrak s\in\mathfrak S_{s\ell}$
is a permutation, then let 
$$
m^{\underline{a},\mathfrak s}:=m_{\kappa_1}\otimes\cdots\otimes m_{\kappa_{s\ell}}\in (M_{\tR}(\lambda))^{\otimes{s\ell}}
$$
be the tensor product of the same factors as in $m^{\underline{a}}$, with the first factor $m_{\kappa_1}$ be the same
as the $\mathfrak s(1)$-th factor in $m^{\underline{a}}$, and the second factor $m_{\kappa_2}$ be the same
as the $\mathfrak s(2)$-th factor in $m^{\underline{a}}$ and so on. 
\begin{remark}\label{importantRemark}
The following proposition provides an important combinatorial tool, we will refer to it later quite often.
For this reason it is important to keep the following in mind: the weight space decomposition of $M_{\tR}(\lam)$ 
makes the module $M_{\tR}(\lam)^{\otimes{s\ell }}$ being  multigraded, so an element can be written as a linear 
combination of its multihomogeneous components.

The summands presented in \eqref{summandsOfva} are multihomogeneous, but they do not correspond
necessarily to a multihomogeneous decomposition. Indeed, several different terms might be of the same 
multihomogeneous degree. The summands listed in \eqref{summandsOfva} are coming from the successive
application of the root vectors $X_{-i}$ respectively the corresponding elements $E_{-i}$ in the quantum group.

There is one case we will encouter later, where this distinction is not necessary. Let $\kappa_1,\ldots,\kappa_{s\ell}$
be elements in $W/W_Q$. The multihomogeneous weight space of weight $(\kappa_1(\lambda),\ldots,\kappa_{s\ell}(\lambda))$
has rank one and is spanned by $m_{\kappa_1}\otimes \cdots\otimes m_{\kappa_{s\ell}}$. So if in some 
multihomogeneous component of $v_{\underline a,\underline{\sigma}}=X_{-{i_1}}^{(n_1)}\cdots X_{-{i_r}}^{(n_r)}
\circ m_\lam^{\otimes s\ell}$, viewed as an element in $M_{\tR}(\lam)^{\otimes{s\ell }}$, a summand of the form 
$m_{\kappa_1}\otimes \cdots\otimes m_{\kappa_{s\ell}}$ shows up, then also one of the terms in 
\eqref{summandsOfva} is, up to a scalar multiple, equal to $m_{\kappa_1}\otimes \cdots\otimes m_{\kappa_{s\ell}}$.
Recall that elements $m^{\underline{a}}$ for $\underline{a}\in \mathrm{LS}_{\lambda}^+(s)$ are of this form.
\end{remark}

\begin{proposition}\label{leadingtermauxbasis}
Let $\underline{a}\in \mathrm{LS}_{\lambda,\mathfrak C}^+(s)$, let $\underline{\sigma}=s_{i_1}\cdots s_{i_r}$ be a reduced decomposition
of the maximal element in $\supp\underline{a}$, and let $\ell$ be an even number such that 
$\ell a_\sigma\in \mathbb N$ for all $\sigma\in \supp\underline{a}$.
\begin{itemize}
\item[\textit{i)}] In the expression for 
$v_{\underline a,\underline{\sigma}}=X_{-{i_1}}^{(n_1)}\cdots X_{-{i_r}}^{(n_r)}
\circ m_\lam^{\otimes s\ell}$ as an element in $M_{\tR}(\lam)^{\otimes{s\ell }}$ we get:
$v_{\underline a,\underline{\sigma}}=
 m^{\underline{a}}+\textrm{summands of the form}$
\begin{equation}\label{summandsOfva}
(E_{-i_1}^{(h_1)}\ldots E_{-i_r}^{(h_r)}m_{\lam})\otimes\ldots\otimes
(E_{-i_1}^{(p_1)}\ldots  E_{-i_r}^{(p_r)}m_{\lam})
\end{equation}
which either vanish or are strictly smaller than $m^{\underline{a}}$:
$$
m^{\underline{a}}\RHD
(E_{-i_1}^{(h_1)}\ldots E_{-i_r}^{(h_r)}m_{\lam})\otimes\ldots\otimes
(E_{-i_1}^{(p_1)}\ldots E_{-i_r}^{(p_r)}m_{\lam}).$$
In particular, one of these summands is a non-zero multiple of $m^{\underline{a}''}$ for
$\underline{a}''\in  \mathrm{LS}_{\lambda}^+(s)$
only if $\underline{a}\rhd\underline{a}''$.
\item[\textit{ii)}] Let $\mathfrak s\in\mathfrak S_{s\ell}$ be a permutation. 
In the expression for  $v_{\underline a,\underline{\sigma}}$ as an element in $M_{\tR}(\lam)^{\otimes{s\ell }}$ we get:
 $v_{\underline a,\underline{\sigma}}=  \v^{\bullet} m^{\underline{a},\mathfrak s}+$  summands strictly smaller
 than $m^{\underline{a},\mathfrak s}$ with respect to $\RHD_{\mathfrak s}$. In particular, one of 
 these summands is a non-zero multiple of $m^{\underline{a}'',\mathfrak s}$ for $\underline{a}''\in  \mathrm{LS}_{\lambda}^+(s)$
only if $\underline{a}\rhd\underline{a}''$.
\end{itemize}
\end{proposition}

Before we come to the proof, let us point out one important consequence.
The summands occurring in \eqref{summandsOfva} may become linearly dependent, cancel each
other or add up to something. So recall that $\RHD$ is a preorder on the
associated tuples in \eqref{summandsOfva} and not on the elements of $M_{\tR}(\lam)^{\otimes{s\ell }}$. 

Let $\underline{a}, \underline{a}''\in  \mathrm{LS}_{\lambda}^+(s)$, $\underline{a}\not= \underline{a}''$,
and suppose  we have written $v_{\underline a,\underline{\sigma}}$ as a linear combination
of multihomogeneous tensors. As pointed out in Remark~\ref{importantRemark},
a non-zero multiple of  $m^{\underline{a}''}$ occurs in such a sum with a non-zero coefficient only if 
at least one of the terms in \eqref{summandsOfva} is a non-zero multiple of $m^{\underline{a}''}$.
Summarizing we have:

\begin{coro}\label{MultigradedWeightSpace}
In an expression of $v_{\underline a,\underline{\sigma}}$ as a linear combination
of multihomogeneous elements in $M_{\tR}(\lam)^{\otimes s\ell}$, the coefficient in front of 
$m^{\underline{a}''}$ is different from zero only if $\underline{a}\rhd\underline{a}''$.
\end{coro}

\begin{proof}(of Proposition \ref{leadingtermauxbasis})
The claim in \textit{i)} holds obviously if the maximal element in $\supp{\underline{a}}$ is
of length zero. We proceed by induction on the length of the maximal element in
$\supp{\underline{a}}$.

Let $\underline{a}'$ be obtained from $\underline{a}$ as in \eqref{aprim}.
By induction, we know that $v_{\underline{a}',\underline{s_{i_1}\sigma}}$ has, viewed as an element in $M_{\tR}(\lambda)^{\otimes{s\ell }}$,
a unique maximal summand which is equal to
$m^{\underline{a}'}$. This element has the form
\begin{equation}\label{mastrich}
m^{\underline{a}'}:=
\underbrace{ m_{s_{i_1}\tau_q}\otimes\ldots\otimes m_{s_{i_1}\tau_q} }_{\ell a_q}
\otimes\ldots
\otimes \left(
\begin{array}{c}
\textrm{\rm remaining part} \\
\textrm{\rm same as in\ }m^{\underline{a}}
\end{array}
\right).
\end{equation}
\begin{enumerate}
\item[(i)] If we apply $X_{-i_1}^{(n_1)}$ to the term $m^{\underline{a}'}$, then we get
a sum that runs over all possible $(s\cdot\ell)$-tuples:
\begin{equation}\label{sumXi1}
X_{-i_1}^{(n_1)}\circ m^{\underline{a}'}=\sum \v^{\bullet} (E_{-i_1}^{(k_1)}m_{s_{i_1}\tau_q})\otimes\ldots\otimes
(E_{-i_1}^{(k_{s\ell})}m_{\tau_{0}}),
\end{equation}
such that $n_1\ell=k_1+\ldots+k_{s\ell}$.  Since all tensor factors are weight vectors and the $K_{-i}$ acts by multiplication
with some root of unity, we omit the $K_{-i}$ in the tensor product formula above and summarize the action of the powers of the $K_{-i}$
by the symbol $\v^{\bullet}$ for ``some root of unity''.

For an extremal weight vector  $m_{\kappa}$ in $M_{\tR}(\lambda)$ we know: $E_{-i}^{(k)} m_{\kappa}=0$
if $k>\vert\langle\kappa(\lambda),\alpha_i^\vee\rangle\vert$.
So if one of the $k$'s is too big, then such a term is zero.
With respect to the chosen ordering, for a maximal term in the sum \eqref{sumXi1} the $k_1$ will be necessarily maximal,
and hence $k_1=\langle s_{i_1}\tau_q(\lam),\alpha^\vee_{i_1}\rangle$. A summand in \eqref{sumXi1}
will be a maximal term with respect to ``$\LHD$'' only if the first tensor factor is (up to rescaling) equal to $m_{\tau_q}$. 

If a summand in  \eqref{sumXi1}  is such that $k_1$ is strictly smaller than $\langle s_{i_1}\tau_q(\lam),\alpha^\vee_{i_1}\rangle$, 
then we get a weight vector in $M_\tR(\lam)_{\tau_q}=U_{\xi}(\mathfrak g)^+_\tR\cdot m_{\tau_q}$, 
say of weight $\mu$, such that  $\mu\succ\tau_q(\lam)$. 
It follows that, independent of the form of the remaining part of the tensor product, 
this summand is either zero or strictly smaller than
$m^{\underline{a}}$ with respect to the preorder ``$\RHD$''.

The same arguments apply to the exponents $k_2$, $k_3$ etc. up to $k_{max}$,
where $k_{max}=\ell(a_q+\ldots+a_{j+1})$. Recall that
$n_1=a_q\vert\langle\tau_q(\lambda),\alpha_{i_1}^\vee\rangle\vert+\ldots+
a_{j+1}\vert\langle\tau_{j+1}(\lambda),\alpha_{i_1}^\vee\rangle\vert)$ by \eqref{weightdifference}.
So applying $X_{-i_1}^{(n_1)}$ to the term $m^{\underline{a}'}$ gives (up to a root of unity) the desired maximal
term $m^{\underline{a}}$ (and this one shows up exactly once), plus terms  which are strictly smaller than
$m^{\underline{a}}$ with respect to the preorder ``$\LHD$''.

It remains to show that the coefficient in front of $m^{\underline{a}}$ is equal to $1$ and not just some root of unity. 
We write $m^{\underline{a}'}$ as $\mathfrak m'_1\otimes \mathfrak m'_2$, where the first part has
$\ell(a_q+\ldots+a_{j+1})$-tensor factors and is of the form $m_{s_{i_1}\tau_q}\otimes \cdots\otimes m_{s_{i_1}\tau_{j+1}}$.
For the tensor product we get (in the  comultiplication rule \eqref{comultiplication_rule} we have $p''=0$ for the first term):
$$
X_{-i_1}^{(n_1)}\circ m^{\underline{a}'}=E_{-i_1}^{(\ell n_1)}\cdot (\mathfrak m'_1\otimes \mathfrak m'_2)
=(E_{-i_1}^{(\ell n_1)}\cdot \mathfrak m'_1)\otimes (K_{-i_1}^{\ell n_1}\cdot \mathfrak m'_2) +\textrm{smaller terms}.
$$
Since the tensor factors in $ \mathfrak m'_1$ are all extremal weight vectors and the divided power of $E_{-i_1}$ is maximal,
rewriting $(E_{-i_1}^{(\ell n_1)}\cdot \mathfrak m'_1)$ as a sum as in $\eqref{sumXi1}$, only one
summand is non-zero, and this is $\v^{\bullet} m_{\tau_q}\otimes \cdots\otimes m_{\tau_{j+1}}$.
But now a simple induction procedure shows that in this case the root of unities in front of the term in \eqref{comultiplication_rule}
get cancelled by the powers of the roots of unities coming from the action of the powers of $K_{-i_1}$ on the tensor factors.
So we get $\v^{\bullet} =1$ for this part of the tensor product. 

It remains to consider the second part:  $(K_{-i_1}^{\ell n_1}\cdot \mathfrak m'_2)$. We know the weight $\mu=\textrm{weight}({\underline{a}})$
is an element in $\Lambda$. In particular, $\langle\mu,\alpha_{i_1}^\vee\rangle\in\mathbb Z$. Let $\mu_1$ be the weight
$a_q\tau_q(\lambda)+\ldots+a_{j+1}\tau_{j+1}(\lambda)$ and set $\mu_2=\mu-\mu_1$. These are in general just rational weights,
but we know: $\langle\tau_q(\lambda),\alpha_{i_1}^\vee\rangle,\ldots,\langle\tau_{j+1}(\lambda),\alpha_{i_1}^\vee\rangle \le 0$
and $\langle\tau_{j}(\lambda),\alpha_{i_1}^\vee\rangle > 0$. By the local integrality property (see \cite{L2}), this implies
$\langle\mu_1,\alpha_{i_1}^\vee\rangle\in\mathbb Z$ and hence $\langle\mu_2,\alpha_{i_1}^\vee\rangle\in\mathbb Z$.
It follows: $(K_{-i_1}^{\ell n_1}\cdot \mathfrak m'_2)=\xi^{\ell n_1\langle\ell\mu_2,\alpha_{i_1}^\vee\rangle}\mathfrak m'_2=\mathfrak m'_2$
because: $\langle \mu_2,\alpha_{i_1}^\vee\rangle\in\mathbb Z$, $\ell$ is even and hence $\ell^2$ is divisible by $2\ell$, and $\xi$ is a $2\ell$-th
root of unity.

\item[(ii)] Next consider a summand $\mathfrak m$ in the expression of $v_{\underline{a}',\underline{s_{i_1}\sigma}}$ such that 
$\mathfrak m\not = m^{\underline{a}'}$, but the first $\ell(a_q+\ldots +a_{j+1})$ tensor factors 
of $\mathfrak m$ and $m^{\underline{a}'}$ coincide. In other words, up to some non-zero constant,
$$
\mathfrak m= \underbrace{m_{s_{i_1}\tau_q}\otimes \cdots\otimes m_{s_{i_1}\tau_{j+1}}}_{\ell(a_q+\ldots +a_{j+1})\textrm{-tensor factors}}\otimes
\mathfrak m'.
$$
The same procedure as in 1) shows: applying $X_{-i_1}^{(n_1)}$ to the term $\mathfrak m$ gives 
a sum of tensor products of weight vectors which has a unique maximal element with respect to ``$\LHD$'',
and this element is 
$$
\widehat{\mathfrak m}=
\underbrace{m_{\tau_q}\otimes \cdots\otimes m_{\tau_{j+1}}}_{\ell(a_q+\ldots +a_{j+1})\textrm{-tensor factors}}\otimes \mathfrak m',
$$
so the first $\ell(a_q+\ldots +a_{j+1})$ tensor factors 
of $\widehat{\mathfrak m}$ and $m^{\underline{a}}$ coincide. Now the assumption $m^{\underline{a}'}\RHD \mathfrak m$
implies $m^{\underline{a}}\RHD \widehat{\mathfrak m}$, and hence $m^{\underline{a}}$ is larger with respect to ``$\RHD$''
than all summands in this sum.

\item[(iii)] It remains to consider the following case: $\mathfrak m$ is a summand in the expression of $v_{\underline{a}',\underline{s_{i_1}\sigma}}$ and 
there exists a $1\le d \leq \ell(a_q+\ldots +a_{j+1})$, such that (up to non-zero constants) the first $(d-1)$ tensor factors of $\mathfrak m$ and $m^{\underline{a}'}$
coincide, but the $d$-th factor does not agree. The $d$-th tensor factor of $m^{\underline{a}'}$ is an extremal weight vector,
say $m_{s_{i_1} \tau_i}$ for some $q\ge i\ge j+1$. So $m^{\underline{a}'}\RHD \mathfrak m$ implies for 
the $d$-th tensor factor $m_{\nu_{d}}$ of $\mathfrak m$: for the weight $\nu_{d}$ we have $s_{i_1} \tau_i(\lambda)\prec \nu_{d}$,
and for the weight vector $m_{\nu_{d}}$ we have $m_{\nu_{d}}\in M_\tR(\lam)_{s_{i_1}\tau_i}
=U_{\xi}(\mathfrak g)^+_\tR\cdot m_{s_{i_1}\tau_i}$.

We apply now $X_{-i_1}^{(n_1)}$ to the term $\mathfrak m$. We consider first what happens to the first 
$(d-1)$ tensor factors. The same procedure as in 1) shows: if one does not apply the maximal  possible divided power of 
$E_{-i_1}$ to the first, second, $\ldots$, $(d-1)$-th tensor factor (maximal means applying gives still a non-zero result), 
then the resulting tensor product is a weight vector which is automatically strictly smaller than  $m^{\underline{a}}$ with respect to ``$\RHD$''. 

So let us assume we have a summand $\widehat{\mathfrak m}$ in $X_{-i_1}^{(n_1)}\circ\mathfrak m$ where one has applied the maximal possible 
divided power of $E_{-i_1}$ to the first $(d-1)$ tensor factors
of $\mathfrak m$. But whatever divided power of $E_{-i_1}$ one applies to $m_{\nu_{d}}$ ($=$ the $d$-th tensor
factor): $m_{\nu_{d}}\in M_\tR(\lam)_{s_{i_1}\tau_i}$ implies $E_{-i_1}^{(k)}\cdot m_{\nu_{d}}\in M_\tR(\lam)_{\tau_i}$ for all $k\ge 0$.
And, if the result for $k>0$ is not equal to zero,  then the $\alpha_{i_1}$-string through $\nu_{d}$ will never meet $\tau_i(\lambda)$,
which implies $m^{\underline{a}}\RHD \widehat{\mathfrak m}$.

Now for $\underline{a}''\in  \mathrm{LS}_{\lambda}^+(s)$ the inequality
$m^{\underline{a}}\RHD m^{\underline{a}''}$ implies automatically $\underline{a}\rhd\underline{a}''$,
which finishes the proof of part one of the proposition.
\end{enumerate}

It remains to prove the second part of the proposition. As we have seen in \eqref{comultiplication_rule}, up to 
multiplication by a root of unity (coming from the action of some $K_j$), the summands in \eqref{summandsOfva}
showing up in the expression of $v_{\underline{a}, \underline{\sigma}}$ are symmetric with respect to permutation of the factors.

Since we argued always only using the summands as in \eqref{summandsOfva} (remember, we look at them rather as tuples
and not as elements in a $\tR$-module), we can apply the arguments as before also to prove part two. 

The only difference is the change of the order in which we look at the tensor factors. Instead of
using the lexicographic indexing of the tensor factors: first the left most factor, then the second most left factor etc., 
we take the enumeration given by the permutation $\mathfrak s$: the ``first'' factor is the $\mathfrak s(1)$-th, 
the ``second'' factor is the $\mathfrak s(2)$-th etc. Now we 
use the preorder $\RHD_{\mathfrak s}$,  and the same arguments as before to show that in the expression for 
$v_{\underline{a},\underline{\sigma}}$, with respect to the preorder $\RHD_{\mathfrak s}$, there is a unique maximal element 
$m^{\underline{a},\mathfrak s}$ showing up as a summand with a root of unity as a coefficient. And all
the other terms are strictly smaller than $m^{\underline{a},\mathfrak s}$ with respect to $\RHD_{\mathfrak s}$.

If one of the terms in \eqref{summandsOfva} is, up to rescaling, equal to $m^{\underline{a}'',\mathfrak s}$,
then, by the symmetry, up to rescaling, $m^{\underline{a}''}$ appears as one of the terms in \eqref{summandsOfva}. Now $m^{\underline{a},\mathfrak s}\RHD_{\mathfrak s} m^{\underline{a}'',\mathfrak s}$ implies $m^{\underline{a}}\RHD m^{\underline{a}''}$ and hence $\underline{a}\rhd \underline{a}''$.
\end{proof}

\subsection{A $U_{\mathbb Z}(\mathfrak g)^+$-stable filtration with one dimensional leaves}\label{Zfiltration}
Let $\mu\in\Lambda$ be an integral weight. Denote by $V_{\mathbb Z}(s\lambda)_\mu$ the corresponding weight space
in $V_{\mathbb Z}(s\lambda)$. Given $\tau\in W/W_Q$, let $\mathrm{LS}_{\lambda}^+(s)$ be the set of 
elements of degree $s$ in the fan of monoids $\mathrm{LS}_{\lambda}^+$ (associated to the partially 
ordered set with bonds $A_\tau$ in Section~\ref{CandidatForMonoid}).

The vectors $v_{\underline{a},\underline\sigma}$
for $\underline{a}\in\mathrm{LS}_{\lambda}^+(s)$ are weight vectors of weight $\textrm{weight}(\underline{a})$ in $V_{\mathbb Z}(s\lambda)_\tau$.
Recall that $\underline{\sigma}$ is a reduced decomposition of the 
maximal element in $\supp\underline{a}$.

If $V_{\mathbb Z}(s\lambda)_\mu\not=\{0\}$, then there exists a large enough element $\tau\in W/W_Q$, such that 
the weight space $V_{\mathbb Z}(s\lambda)_\mu$ is completely contained in the Demazure module $V_{\mathbb Z}(s\lambda)_\tau$. 

In the following we fix $\mu$ such that $V_{\mathbb Z}(s\lambda)_\mu\not=\{0\}$, and 
$\tau\in W/W_Q$ such that $V_{\mathbb Z}(s\lambda)_\mu$ is completely contained in the Demazure module $V_{\mathbb Z}(s\lambda)_\tau$.

For every $\underline{a}\in\mathrm{LS}_{\lambda}^+(s)$ of weight $\textrm{weight}(\underline{a})=\mu$ 
fix a reduced decomposition $\underline\sigma$ of the maximal element $\sigma$ in $\supp \underline{a}$. 
We denote 
$$\mathbb B(s\lambda)_\mu:=\{v_{\underline{a}, \underline{\sigma}}\mid \underline{a}\in\mathrm{LS}^+_\lambda(s), \textrm{weight}(\underline{a})=\mu \}.$$

\begin{lemma}\label{basisweightspace}
The set $\mathbb B(s\lambda)_\mu$ form a basis for $V_{\mathbb Z}(s\lambda)_\mu$.
\end{lemma}

\begin{proof}
The set $\mathrm{LS}_{\lambda}^+(s)$ is finite, so we can fix an even number $\ell$ such that $\ell a_\kappa\in\mathbb N$ for 
all $\underline{a}\in  \mathrm{LS}_{\lambda}^+(s)$. We consider the embedding $V_{\tR}(s\lambda)_\mu\hookrightarrow
V_{\tR}(s\lambda)_\tau\hookrightarrow
M_\tR(\lam)_{\tau}^{\otimes{s\ell }}$. By Proposition~\ref{leadingtermauxbasis}, an element 
$v_{\underline{a}, \underline{\sigma}}\in \mathbb B(s\lambda)_\mu$ has as maximal term $m^{\underline{a}}$ in its expression in 
$M_\tR(\lam)_\tau^{\otimes{s\ell }}$, independent of the choice of the reduced decomposition  $\underline{\sigma}$.
It follows immediately that the elements in $\mathbb B(s\lambda)_\mu$ are linearly independent over $\tR$. The coefficient
of the leading term is $1$, so they remain linearly independent after base change with any algebraically closed field $\mathbb K$.
Now the Kac-Weyl character formula for $V_{\mathbb C}(s\lambda)$ and the character formula for the path model
(\cite{L3}, Theorem 9.1) implies that the vectors $v_{\underline{a}, \underline{\sigma}}$ of weight $\mu$ are not only linearly independent but form 
a basis for the weight space $V_{\mathbb C}(s\lambda)_\mu$. By construction, the dimensions of the weight spaces
are independent of the choice of the algebraically closed field $\mathbb K$, so the (images of the) vectors $v_{\underline{a}, \underline{\sigma}}$
form a basis for  $V_{\mathbb K}(s\lambda)_\mu$ for any algebraically closed field $\mathbb K$.

As a consequence we see that the $\mathbb Z$-lattice $L$ spanned by the $v_{\underline{a}, \underline{\sigma}}$ of weight $\mu$ in
$V_{\mathbb Z}(s\lambda)_\mu$ has the same rank as $V_{\mathbb Z}(s\lambda)_\mu$. It follows that 
the quotient $V_{\mathbb Z}(s\lambda)_\mu/L$ is a torsion module. Let $\bar v\in V_{\mathbb Z}(s\lambda)_\mu/L\setminus\{0\}$
and fix a representative $v\in V_{\mathbb Z}(s\lambda)_\mu$. Let $d>1$ be the minimal positive integer such that $dv\in L$, 
say $dv=\sum_{\underline{a}} d_{\underline{a}} v_{\underline{a}, \underline{\sigma}}$, where the sum is running over all $\underline{a}$ of weight $\mu$. 

Since $d$ is supposed to be minimal, there exists a  prime $p$ dividing $d$ and an index $\underline{a}'$ such that 
$d_{\underline{a}'}\not=0$ and $p$ does not divide $d_{\underline{a}'}$. If $\mathbb K$ is an algebraically closed field
of characteristic $p$, then the image of $dv$ in $V_{\mathbb K}(s\lambda)_\mu$ is zero, which implies that
the linear combination $\sum_{\underline{a}} \bar d_{\underline{a}} v_{\underline{a},\underline{\sigma}}$ in $V_{\mathbb K}(s\lambda)_\mu$
is equal to zero. Since the images of the $v_{\underline{a}, \underline{\sigma}}$ in $V_{\mathbb K}(s\lambda)_\mu$ are linearly independent,
this implies that the coefficients $\bar d_{\underline{a}}$ are all equal to zero, in contradiction to the assumption that 
$d_{\underline{a}'}$ is not divisible by $p$. 
It follows: $V_{\mathbb Z}(s\lambda)_\mu/L=\{0\}$, and hence $\mathbb B(s\lambda)_\mu$ spans $V_{\mathbb Z}(s\lambda)_{\mu}$.
\end{proof}

Given $\underline{a}\in \mathrm{LS}_{\lambda}^+(s)$ and the vectors  $v_{\underline{a},\underline{\sigma}}\in V_{\mathbb Z}(s\lambda)_\tau$, 
it is natural to inspect the dependence  on the choice of the reduced decompositions $\underline{\sigma}$. 
Let $\underline{\sigma}'$ be a different reduced decomposition of the maximal element in $\supp\underline{a}$, and let 
$v_{\underline{a},\underline{\sigma}'}\in V_{\mathbb Z}(s\lambda)_\tau$ be the corresponding vector.
\begin{lemma}\label{uppertriangualrdecomp}
We have:
$$
v_{\underline{a},\underline{\sigma}'}=v_{\underline{a},\underline{\sigma}}+
\sum_{\substack{\underline{a}''\in \mathrm{LS}_{\lambda}^+(s)\\ 
\mathrm{weight}(\underline{a}'')=\mu;\,\underline{a}\rhd \underline{a}''}}
b_{\underline{a}''}v_{\underline{a}'',\underline{\sigma}''}.
$$
\end{lemma}
\begin{proof}
Let $\mu =\textrm{weight}(\underline{a})$. We assume first $\tau'\ge \tau$ is large enough so that 
 $V_{\mathbb Z}(s\lambda)_\mu$ is completely contained in  $V_{\mathbb Z}(s\lambda)_{\tau'}$.
 We fix a basis $\mathbb B(s\lambda)_\mu$ of $V_{\mathbb Z}(s\lambda)_\mu$ as in Lemma~\ref{basisweightspace}.
 Without loss of generality we can fix the basis so that $v_{\underline{a},\underline{\sigma}}\in \mathbb B(s\lambda)_\mu$.
 Lemma~\ref{basisweightspace} implies then: we can write $v_{\underline{a},\underline{\sigma}'}$ as a linear combination
 of the elements in $\mathbb B(s\lambda)_\mu$.

To be more precise, we use the Frobenius splitting trick.
 The set $\mathrm{LS}_{\lambda,\tau'}^+(s)$ is finite, so we can fix an even number $\ell$ such that $\ell a_\kappa\in\mathbb N$ for
all $\underline{a}\in  \mathrm{LS}_{\lambda, \tau'}^+(s)$. We consider the embedding $V_{\tR}(s\lambda)_{\mu}\hookrightarrow
V_{\tR}(s\lambda)_{\tau'}\hookrightarrow M_\tR(\lam)_{\tau'}^{\otimes{s\ell }}$. As an immediate consequence of Proposition~\ref{leadingtermauxbasis}
and Corollary~\ref{MultigradedWeightSpace} we get
$$
v_{\underline{a},\underline{\sigma}'}=v_{\underline{a},\underline{\sigma}}+
\sum_{\substack{v_{\underline{a}'',\underline{\sigma}''}\in \mathbb B(s\lambda)_\mu\\ 
\underline{a}\rhd \underline{a}''}}
b_{\underline{a}''}v_{\underline{a}'',\underline{\sigma}''}.
$$
It remains to replace $\tau'$ by $\tau$ again. Now $\underline{a}\rhd \underline{a}''$ implies in particular that $\sigma$,  the 
largest element in $\supp\underline{a}$, is greater or equal in the Bruhat order than $\sigma''$, the largest element in $\supp\underline{a}''$.
It follows that  $v_{\underline{a}'',\underline{\sigma}''}\in \mathbb B(s\lambda)_\mu$ and $\underline{a}\rhd \underline{a}''$
implies automatically $\underline{a}\in \mathrm{LS}_{\lambda}^+(s)$, $\textrm{weight}(\underline{a})=\mu$ and 
$\underline{a}\rhd \underline{a}''$, which finishes the proof.
\end{proof}
Because of Lemma~\ref{uppertriangualrdecomp} it makes sense to write just $v_{\underline{a}}$ instead
of $v_{\underline{a},\underline{\sigma}}$ if no confusion is possible. It means we fix for $\underline{a}\in 
\mathrm{LS}_{\lambda}^+(s)$ a reduced decomposition $\underline{\sigma}$ for the maximal element in $\supp{\underline{a}}$, 
and the difference  $v_{\underline{a},\underline{\sigma}}-v_{\underline{a},\underline{\sigma}'}$ we get by choosing a different
decomposition $\underline{\sigma}'$ can be neglected in the corresponding context.

\begin{lemma}\label{BasisDividedPower}
Let $\underline{a}\in \mathrm{LS}_{\lambda}^+(s)$ be of weight $\mu$ and $1\le i\le n$.
If $k\ge 1$, then $X_{i}^{(k)} v_{\underline{a}}=\sum_{\underline{a}\rhd \underline{a}'} d_{\underline{a},\underline{a}'} v_{\underline{a}'}$
for some $\underline{a'}\in \mathrm{LS}_{\lambda}^+(s)$.
\end{lemma}

\begin{proof}\rm
If necessary, we choose first an element $\tau'\ge \tau$ which is large enough so that for  $\mu'=\mu+k\alpha_i$ with $k\ge 0$ we have:
$V_{\mathbb Z}(s\lambda)_{\mu'} \subseteq V_{\mathbb Z}(s\lambda)_{\tau'}$.

The set $\mathrm{LS}_{\lambda,\tau'}^+(s)$ associated to $A_{\tau'}$  is finite, so we can fix an even integer $\ell$ such that 
$\ell a_\kappa\in\mathbb N$ for all $\underline{a}\in  \mathrm{LS}_{\lambda,\tau'}^+(s)$ (we add a $\tau'$ as index to emphasize
that we consider the set $A_{\tau'}$). 

We consider the embedding $V_{\mathbb Z}(s\lambda)_\mu
\hookrightarrow V_{\tR}(s\lambda)_{\tau'}\hookrightarrow
M_\tR(\lam)_{\tau'}^{\otimes{s\ell }}$.
By Proposition~\ref{leadingtermauxbasis} we can write 
$v_{\underline{a}}$  in $M_\tR(\lam)_{\tau'}^{\otimes{s\ell }}$  as 
$m^{\underline{a}}$ plus multihomogeneous terms which are strictly smaller than $m^{\underline{a}}$ with respect to 
$\RHD$. Applying $X_{i}^{(k)}$ to one of the latter terms gives again a sum of tensor products of weight vectors,
all strictly smaller than $m^{\underline{a}}$.

It remains to consider $X_{i}^{(k)}\circ m^{\underline{a}}$. By the definition
of $\RHD$, the result is a sum of tensor products of weight vectors, all strictly smaller than $m^{\underline{a}}$.

The vector $X_{i}^{(k)} v_{\underline{a}}$ is of weight $\mu'=\textrm{weight}(\underline{a})+k\alpha_i$.
The weight space $V_{\mathbb Z}(s\lambda)_{\mu'}$ has the set $\mathbb B(s\lambda)_{\mu'}$ as basis (Lemma~\ref{basisweightspace}), 
it consists of vectors of the form $v_{\underline{a}'}$, $\textrm{weight}(\underline{a}')=\mu'$. 

Now Proposition~\ref{leadingtermauxbasis} together with Corollary~\ref{MultigradedWeightSpace} implies: 
in an expression of $X_{i}^{(k)} v_{\underline{a}}$
as a linear combination of the elements in $\mathbb B(s\lambda)_{\mu'}$ only those $v_{\underline{a}'}$ can occur
such that $\underline{a}\rhd \underline{a}'$.

Now we come back to the assumption at the beginning: $\tau'$ is large enough so that 
$V_{\mathbb Z}(s\lambda)_{\mu'} \subseteq V_{\mathbb Z}(s\lambda)_{\tau'}$. But $\underline{a}\rhd \underline{a}'$ 
implies in particular that the maximal element in $\supp\underline{a}$ is greater or equal to the maximal element in $\supp\underline{a'}$.
But this implies automatically: if $\tau$ is greater or equal  to the maximal elements in $\supp\underline{a}$, then $v_{\underline{a}}$
and $v_{\underline{a}'}$ are already elements in $V_{\mathbb Z}(s\lambda)_\tau$. So the expression
$X_{i}^{(k)} v_{\underline{a}}=\sum_{\underline{a}\rhd \underline{a}'} d_{\underline{a},\underline{a}'} v_{\underline{a}'}$
holds already in $V_{\mathbb Z}(s\lambda)_\tau$.
\end{proof}
{  We summarize the results of this section in the following (already known) statements and in a theorem:
\begin{itemize}
\item[\it i)] The Demazure module $V_{\mathbb Z}(s\lambda)_{\tau}$ is a direct summand of the $\mathbb Z$-module $V_{\mathbb Z}(s\lambda)$,
having as basis the set $\{v_{\underline{a}'}\mid \underline{a}'\in \mathrm{LS}_{\lambda}^+(s)\}$.
\item[\it ii)] For $\underline{a}\in \mathrm{LS}_{\lambda}^+(s)$, the $\mathbb Z$-submodule $V_{\mathbb Z}(s\lambda)_{\tau,\underline{a}}$ of $V_{\mathbb Z}(s\lambda)_\tau$ is a $U(\mathfrak g)_{\mathbb Z}^+$-stable 
direct summand of $V_{\mathbb Z}(s\lambda)_{\tau}$, and it has as basis the set 
$\{v_{\underline{a}'}\mid \underline{a}'\in \mathrm{LS}_{\lambda}^+(s),
\underline{a}\unrhd \underline{a}'\}$.
\end{itemize}
\begin{theorem}\label{AuxBasisTheorem}
The $\mathbb Z$-submodules $V_{\mathbb Z}(s\lambda)_{\tau,\underline{a}}$, 
$ \underline{a}\in \mathrm{LS}_{\lambda}^+(s)$, define a 
$U(\mathfrak g)_{\mathbb Z}^+$-stable collection of  subspaces of
$V_{\mathbb Z}(s\lambda)_{\tau}$ with leaves 
$V_{\mathbb Z}(s\lambda)_{\tau,\underline{a}}/V_{\mathbb Z}(s\lambda)_{\tau,\lhd\underline{a}}$ free of rank 1. 
\end{theorem}}
\begin{proof}
Given a weight $\mu$, we can always choose an element $\tau'$ such that the weight space $V_{\mathbb Z}(s\lambda)_\mu$
is contained in the Demazure module $V_{\mathbb Z}(s\lambda)_{\tau'}$.  By Lemma~\ref{basisweightspace}, we have hence a basis 
$\mathbb B(s\lambda)_\mu$ for each weight space $V_{\mathbb Z}(s\lambda)_\mu$,
so the union $\mathbb B$ of the $\mathbb B(s\lambda)_\mu$ provides a basis for the representation space $V_{\mathbb Z}(s\lambda)$.

We come now back to the fixed element $\tau$. The subset $\mathbb B_\tau(s)=\{v_{\underline{a}'}\mid \underline{a}'\in \mathrm{LS}_{\lambda}^+(s)\}\subset \mathbb B$ 
consists of basis vectors contained in $V_{\mathbb Z}(s\lambda)_{\tau}$. By Lemma~\ref{BasisDividedPower},
the $\mathbb Z$-span of $\mathbb B_\tau(s)$ is stable under the action of $U(\mathfrak g)_{\mathbb Z}^+$. The extremal weight vector
$v_\tau$ is an element in $\mathbb B_\tau(s)$, which implies that $\mathbb B_\tau(s)$  is a basis for $V_{\mathbb Z}(s\lambda)_{\tau}$,
and hence the latter is a direct summand of $V_{\mathbb Z}(s\lambda)$.

By Lemma~\ref{uppertriangualrdecomp}, the subset  $\mathbb B_{\tau,\underline{a}}(s)=\{v_{\underline{a}'}\in \mathbb B_{\tau}(s)\mid 
\underline{a}'\unlhd \underline{a} \}$  is a generating set for 
$V_{\mathbb Z}(s\lambda)_{\tau,\underline{a}}$ and, because of the linear independency, it is a basis. This implies
also that $V_{\mathbb Z}(s\lambda)_{\tau,\underline{a}}$ is a direct summand of $V_{\mathbb Z}(s\lambda)_{\tau}$.
By Lemma~\ref{BasisDividedPower}, $V_{\mathbb Z}(s\lambda)_{\tau,\underline{a}}$ is stable under 
the action of $U(\mathfrak g)_{\mathbb Z}^+$. The subquotient 
$V_{\mathbb Z}(s\lambda)_{\tau,\underline{a}}/V_{\mathbb Z}(s\lambda)_{\tau,\lhd\underline{a}}$
is spanned by the class $\bar v_{\underline{a}}$ of $v_{\underline{a}}$, and hence is free of rank $1$.
\end{proof}


\section{Appendix III: The path vectors and some relations}\label{The_path_vectors}

We keep the same notation as in Section~\ref{A_filtration_on_V}. 

In the Appendix II we defined a filtration of  $V_{\mathbb Z}(s\lambda)_{\tau}$
given by the subspaces $V_{\mathbb Z}(s\lambda)_{\tau,\underline{a}}$, $\underline{a}\in \mathrm{LS}_{\lambda}^+(s)$. The vectors
$v_{\underline a,\underline \sigma}\in V_{\mathbb Z}(s\lambda)_{\tau,\underline{a}}$ defined in Definition~\ref{definitionvatau}
are representatives of the leaves associated to this filtration on $V_{\mathbb Z}(s\lambda)_{\tau}$.

We will define a kind of dual basis which again induces a filtration, but this time on the dual space. We recall Definition~\ref{define_path_vector}:

\begin{definition}\label{define_path_vector2}
Let $\underline{a}\in \mathrm{LS}_{\lambda}^+(s)$.  A \textit{path vector associated to $\underline a$} is 
a linear function $p_{\underline a}\in  V(s\lambda)^*_{\tau}$ which is a 
$T$-eigenvector of weight $(-\textrm{weight}(\underline a))$, and such that 
\begin{enumerate}
\item[\it i)]
{  there exists a reduced decomposition $\underline{\sigma}$ of the maximal element $\sigma$ in $\mathrm{supp}\,\underline{a}$ such that $p_{\underline a}(v_{\underline a,\underline \sigma})=1$};
\item[\it ii)] for $\underline{a'}\in\mathrm{LS}_\lambda^+(s)$ and for some reduced decomposition $\underline{\sigma'}$ of the maximal element $\sigma'$ in $\mathrm{supp}\,\underline{a'}$,  $p_{\underline a}(v_{\underline a',\underline \sigma'})\not=0$
implies $\underline a'\rhd \underline a$.
\end{enumerate}
\end{definition}

%
\begin{remark} 
This definition is independent of the choice of the reduced decomposition $\underline{\sigma}$ of $\sigma$ (resp. $\underline{\sigma'}$ of $\sigma'$),
see Lemma~\ref{indipendentproofpathvector}. 
\end{remark}

More generally, for an appropriate ring $R$ (resp. an algebraically closed field $\mathbb K$), we say 
 $p_{\underline a}\in  V_{R}(s\lambda)^*_{\tau}$ (resp. $V_{\mathbb K}(s\lambda)^*_{\tau}$) is a \textit{path vector}
 if it satisfies the same conditions with respect to $v_{\underline a',\underline \sigma'}$ and $v_{\underline a,\underline \sigma}\in V_{R}(s\lambda)_{\tau}$
 (resp. $V_{\mathbb K}(s\lambda)_{\tau})$.

We will use the Frobenius splitting trick to construct examples of such functions in $V_{\tR}(s\lambda)^*_{\tau}$.
We fix a dominant weight $\lam\in \Lambda^+$
and an element $\tau\in W/W_Q$. Given $\underline{a}\in {\mathrm{LS}^+_{\lambda}}(s)$, we fix
an even number $\ell$ such that $\ell a_\kappa\in \mathbb N$ for all $\kappa\in\supp\underline a$.
We have seen (see Section~\ref{FrobeniusTrick}) that we have a natural inclusion of $U({\mathfrak g})^+_\tR$-modules:
\begin{equation}\label{quantuminj2}
V_\tR(s\lambda)_{\tau}\hookrightarrow M_\tR(\lambda)_{\tau}^{\otimes s\ell},
\end{equation}
The inclusion induces a restriction map
for the corresponding dual modules:
\begin{equation}\label{quantumres}
\textrm{res}:(M_\tR(\lambda)^*_{\tau})^{\otimes s\ell}\rightarrow V_\tR(s\lambda)_{\tau}^*.
\end{equation}
For an extremal weight vector $m_{\kappa}\in M_\tR(\lambda)_\tau$ denote by $x_\kappa\in M_\tR(\lambda)_\tau^*$
the dual vector. Let $\supp\underline{a}= \{\tau_q,\ldots, \tau_1\}$ with $\tau_q>\ldots> \tau_1$.
As a first step we associate to $\underline{a}$ the tensor product of the following
extremal weight vectors:
\begin{equation}\label{xu}
x_{\underline{a}}:=
\underbrace{x_{\tau_q}\otimes\ldots\otimes x_{\tau_q}}_{{\ell} a_{q}} 
\otimes 
\underbrace{x_{\tau_{q-1}}\otimes\ldots\otimes x_{\tau_{q-1}}}_{{\ell}a_{q-1}} 
\otimes\ldots\otimes
\underbrace{x_{\tau_1}\otimes\ldots\otimes 
x_{\tau_1}}_{{\ell}a_1}\in
(M_\tR(\lam)_{\tau}^*)^{\otimes{s\ell}}.
\end{equation}

\begin{lemma}\label{lpathvector}
The image $p_{\underline{a},\ell}:=\mathrm{res}(x_{\underline{a}})$ in $V_\tR(s\lambda)_{\tau}^*$ is a \textit{path vector}
in $V_\tR(s\lambda)_{\tau}^*$.
\end{lemma}

\begin{proof}
This is a direct consequence of Proposition~\ref{leadingtermauxbasis} and Corollary~\ref{MultigradedWeightSpace}, which imply:
\begin{equation}\label{KindOfDual}
p_{\underline{a},\ell}(v_{\underline{a},\underline\sigma})=1\quad\textrm{and}\quad p_{\underline{a},\ell}(v_{\underline{a}',\underline\sigma'})\not=0
\quad\textrm{only if\ } \underline{a} \unlhd\underline{a}'.
\end{equation}
\end{proof}

\begin{definition}\label{pathvectors}
For $\underline{a}\in\mathrm{LS}_\lambda^+(s)$ and $\ell$ chosen as above, we call $p_{\underline{a},\ell}\in V_\tR(s\lambda)_{\tau}^*$ the \textit{path vector associated to 
$\underline{a}$ and $\ell$}. 
\end{definition}

The set ${\mathrm{LS}^+_{\lambda}}(s)$ is finite, so we can choose an even number $\ell$ such that $\ell a_\kappa\in \mathbb N$ for all 
$\kappa\in\supp\underline a$ and all $\underline a\in {\mathrm{LS}^+_{\lambda}}(s)$. For such a choice of $\ell$ we get:

\begin{lemma}\label{BasisPathVector}
The collection of path vectors $\{p_{\underline{a},\ell}\mid \underline a\in {\mathrm{LS}^+_{\lambda}}(s)\}$ is a basis for 
$V_\tR(s\lambda)_{\tau}^*$.
\end{lemma}

\begin{proof}
By fixing a linearization of the partial order $\unlhd$ on $\mathrm{LS}^+_{\lambda}(s)$, \eqref{KindOfDual} implies that the base change matrix between $\{p_{\underline{a},\ell}\mid \underline{a}\in\mathrm{LS}^+_{\lambda}(s)\}$ and the dual basis of $\{v_{\underline{b}}\mid \underline{b}\in\mathrm{LS}^+_{\lambda}(s)\}$ in Theorem~\ref{AuxBasisTheorem} is upper-triangular unipotent.
\end{proof}

\begin{coro}\label{differentEll}
Let $\ell$ be as above. For $\underline a\in {\mathrm{LS}^+_{\lambda}}(s)$ let $p_{\underline a}\in V_\tR(s\lambda)_{\tau}^*$ be a path vector. Then
$p_{\underline{a}}=p_{\underline{a},\ell} +\sum_{\underline{a}'\rhd \underline{a}} d_{\underline{a}',\ell}p_{\underline{a}',\ell}$.
\end{coro}

\begin{proof}
Fix the same linearization as in the proof of Lemma \ref{BasisPathVector}. From the definition of path vectors, the base change matrix between $\{p_{\underline{a}}\mid\underline{a}\in\mathrm{LS}^+_{\lambda}(s)\}$ and the dual basis of $\{v_{\underline{b}}\mid \underline{b}\in\mathrm{LS}^+_{\lambda}(s)\}$ in Theorem~\ref{AuxBasisTheorem} is also upper-triangular unipotent.
\end{proof}

Therefore it makes sense to define the following  $\tR$-submodules of $V_\tR(s\lambda)_{\tau,\underline a}^*$
for $\underline{a}\in {\mathrm{LS}^+_{\lambda}}(s)$:
\begin{equation}\label{appendix_drei_filtration_one}
V_\tR(s\lambda)_{\tau,\unrhd\underline a}^*
=\left\langle 
p_{\underline{a}'} \bigg\vert p_{\underline{a}'}\textrm{\ a path vector,\ } \underline{a'}\in {\mathrm{LS}^+_{\lambda}}(s),\underline{a}'\unrhd\underline{a}
\right\rangle_\tR;
\end{equation}
and set 
\begin{equation}\label{appendix_drei_filtration_two}
V_\tR(s\lambda)_{\tau,\rhd\underline a}^*=\left\langle 
p_{\underline{a}'} \bigg\vert p_{\underline{a}'}\textrm{\ a path vector,\ }\underline{a}'\in {\mathrm{LS}^+_{\lambda}}(s), \underline{a}'\rhd\underline{a} 
\right\rangle_\tR.
\end{equation}

\begin{coro}
These subspaces define a filtration of $\tR$-modules on $V_\tR(s\lambda)_{\tau}^*$ with leaves $V_\tR(s\lambda)_{\tau,\underline a}^*/V_\tR(s\lambda)_{\tau,\rhd\underline a}^*$ free of rank $1$. A path vector $p_{\underline{a}}$ is a representative of such a leaf 
associated to ${\underline{a}}$.
\end{coro}

\begin{remark}\label{ell}
We  fix an even number $\ell$ such that for all $\underline a$ in ${\mathrm{LS}^+_{\lambda}}(1)$ and  all  $\kappa\in\supp\underline a$
one has: $\ell a_\kappa\in \mathbb N$.
If $\underline a\in{\mathrm{LS}^+_{\lambda}}(s)$ for some $s>1$, then $\underline{a}$ has a unique decomposition
$\underline a=\underline{a}^1+\ldots +\underline{a}^s$, where $\underline{a}^1,\ldots,\underline{a}^s\in {\mathrm{LS}^+_{\lambda}}(1)$ (Lemma \ref{latticedecomp}).
So the number $\ell$ chosen above has the property: for all $m\ge 1$ and all $\underline a$ in ${\mathrm{LS}^+_{\lambda}}(m)$ and all  $\kappa\in\supp\underline a$
one has $\ell a_\kappa\in \mathbb N$.
\end{remark}

\begin{lemma}\label{pathvectorspecialrelation}
Let $\mathfrak C$ be a maximal chain, $\underline{a}\in \mathrm{LS}_{\mathfrak C,\lam}^+(s_1)$
and $\underline{b}\in \mathrm{LS}_{\mathfrak C,\lam}^+(s_2)$. Then $p_{\underline{a},\ell}p_{\underline{b},\ell}
=cp_{\underline{a}+\underline{b},\ell} +\sum_{\underline{a}'\rhd \underline{a}+\underline{b}} d_{\underline{a}'}p_{\underline{a}'}$,
where $c$ is a root of unity, $p_{\underline{a}+\underline{b},\ell}$ is the path vector associated to 
$\ell$ and $\underline{a}+\underline{b}\in \mathrm{LS}_{\mathfrak C,\lam}^+(s_1+s_2)$, and the $p_{\underline{a}'}$
are path vectors  for the leafs associated to $\underline{a'}\in  \mathrm{LS}_{\lam}^+(s_1+s_2)$.
\end{lemma}
\begin{proof}
This is an immediate consequence of part \textit{ii)} of Proposition~\ref{leadingtermauxbasis}, subject to the following 
consideration: We view $p_{\underline{a},\ell}p_{\underline{b},\ell}$ as a function on $V_\tR((s_1+s_2)\lambda)_{\tau}$
via the restriction of $p_{\underline{a},\ell}\otimes p_{\underline{b},\ell}\in V_\tR(s_1\lambda)^*_{\tau}\otimes V_\tR(s_2\lambda)_{\tau}^*$
to $ V_\tR((s_1+s_2)\lambda)_{\tau}\subseteq  V_\tR(s_1\lambda)_{\tau}\otimes V_\tR(s_2\lambda)_{\tau}$.

For the definition of  $p_{\underline{a},\ell}$ (resp. $p_{\underline{b},\ell}$) we have viewed $V_\tR(s_1\lambda)_{\tau}$ (resp. $V_\tR(s_2\lambda)_{\tau}$) as being embedded in $M_\tR(\lambda)_\tau^{\otimes s_1\ell}$ (resp. $M_\tR(\lambda)_\tau^{\otimes s_2\ell}$). So we view $V_\tR((s_1+s_2)\lambda)_{\tau}$ and the tensor product $V_\tR(s_1\lambda)_{\tau}\otimes V_\tR(s_2\lambda)_{\tau}$ as being 
embedded in $M_\tR(\lambda)_\tau^{\otimes (s_1+s_2)\ell}$.

Now one has to be careful with the coproduct. One can view $M_\tR(\lambda)^{\otimes (s_1+s_2)\ell}$ as a $U(\frak g)_\tR^\pm$ module
by viewing $M_\tR(\lambda)^{\otimes (s_1+s_2)\ell}$ as a tensor product of $U_\xi(\frak g)_\tR^\pm$-modules, and one takes
then $\mathrm{Fr'}$ to turn it into a $U(\frak g)_\tR^\pm$-module. This is the way we view the module in the context of Proposition~\ref{leadingtermauxbasis}.

The way we look here at the module is different. We see $M_\tR(\lambda)^{\otimes s_1\ell}$ and $M_\tR(\lambda)^{\otimes s_2\ell}$ 
as tensor products of $U_\xi(\frak g)_\tR^\pm$-modules, and we make them via $\mathrm{Fr'}$ into $U(\frak g)_\tR^\pm$-modules, and one
takes then the tensor product of $U(\frak g)_\tR^\pm$-modules. 

So we have on $M_\tR(\lambda)^{\otimes (s_1+s_2)\ell}$ two different $U(\frak g)_\tR^\pm$-module structures, and this makes a difference
for the image of a vector $v_{\underline a,\underline{\sigma}}$, $\underline a\in  {\mathrm{LS}^+_{\lambda}}(s_1+s_2)$, in $M_\tR(\lambda)^{\otimes (s_1+s_2)\ell}$. But it is now easy to see that the only difference is: to get 
$$
v_{\underline a,\underline{\sigma}}=X_{-{i_1}}^{(n_1)}\cdots X_{-{i_r}}^{(n_r)}\circ (m_\lam^{\otimes s_1\ell}\otimes m_\lam^{\otimes s_2\ell}),
$$ 
with respect to the second $U(\frak g)_\tR^\pm$-structure, one has first to apply the coproduct
for the Lie algebra:
\begin{equation}\label{less_summands}
\sum_{j_1=0}^{n_1}\cdots\sum_{j_r=0}^{n_r}(X_{-{i_1}}^{(j_1)}\cdots X_{-{i_r}}^{(j_r)}\circ m_\lam^{\otimes s_1\ell})
\otimes (X_{-{i_1}}^{(n_1-j_1)}\cdots X_{-{i_r}}^{(n_r-j_r)}\circ m_\lam^{\otimes s_2\ell}),
\end{equation}
and then apply the Frobenius morphism. One checks now as before: the maximal term is still $m^{\underline a}$,
but the coefficient in front of it might be a root of unity. Since we use first the coproduct for the Lie algebra and 
then the Frobenius morphism, there are less summands than in the formula in \eqref{summandsOfva}.

But, up to multiplying by roots of unities, all terms $(E_{-i_1}^{(h_1)}\ldots E_{-i_r}^{(h_r)}m_{\lam})\otimes\ldots\otimes
(E_{-i_1}^{(p_1)}\ldots  E_{-i_r}^{(p_r)}m_{\lam})$ which show up in \eqref{less_summands} after applying the Frobenius
morphism, also occur also in \eqref{summandsOfva}. So we still can apply Proposition~\ref{leadingtermauxbasis}.
\end{proof}

\begin{proposition}\label{product_path_vector}
Let $\mathfrak C$ be a maximal chain, $\underline{a}\in \mathrm{LS}_{\mathfrak C,\lam}^+(s_1)$
and $\underline{b}\in \mathrm{LS}_{\mathfrak C,\lam}^+(s_2)$. If $p_{\underline{a}}$ and $p_{\underline{b}}$
are path vectors associated to the leaves $\underline{a}$ respectively $\underline{b}$, then $p_{\underline{a}}p_{\underline{b}}$ is a path vector associated to $\underline{a}+\underline{b}\in \mathrm{LS}_{\mathfrak C,\lam}^+(s_1+s_2)$, up to multiplying by a root of unity.
\end{proposition}

\begin{proof}
Let $\ell$ be an even number such that $\ell a_\kappa\in \mathbb N$ for all 
$\kappa\in\supp\underline a$ and all $\underline a\in {\mathrm{LS}^+_{\lambda}}(1)$. (The same property holds
then for all $\underline a\in {\mathrm{LS}^+_{\lambda}}(s)$ and all $s\ge 1$, see Remark~\ref{ell}).

By Corollary~\ref{differentEll}, we know for $\underline a\in {\mathrm{LS}^+_{\lambda}}(s_1)$:
$p_{\underline{a}}=p_{\underline{a},\ell} +\sum_{\underline{a}'\unrhd \underline{a}} d_{\underline{a}',\ell}p_{\underline{a}',\ell}$,
and we have a similar expression for $p_{\underline{b}}$.

Lemma~\ref{pathvectorspecialrelation} implies:  $p_{\underline{a},\ell}p_{\underline{b},\ell}$ is a path vector associated to $\underline{a}+\underline{b}\in \mathrm{LS}_{\mathfrak C,\lam}^+(s_1+s_2)$. It remains
to show: $\underline{a}'\unrhd \underline{a}$ and $\underline{b}'\unrhd \underline{b}$ (and strict inequality for at least one term) 
implies $p_{\underline{a}',\ell}p_{\underline{b}',\ell}$ is a sum of path vectors $p_{\underline c'}$ such that 
$\underline{c}'\rhd \underline{a}+\underline{b}$. This holds by Lemma~\ref{pathvectorspecialrelation}
if there exists a maximal chain $\frak C$ containing both $\supp \underline{a}'$ and $\supp \underline{b}'$. 

If this is not the case, then we use the same arguments as in the proof of Lemma~\ref{pathvectorspecialrelation}.
We view $p_{\underline{a}',\ell}p_{\underline{b}',\ell}$ as a function on $V_\tR((s_1+s_2)\lambda)_{\tau}$
via the restriction of $p_{\underline{a}',\ell}\otimes p_{\underline{b}',\ell}\in V_\tR(s_1\lambda)^*_{\tau}\otimes V_\tR(s_2\lambda)_{\tau}^*$
to $ V_\tR((s_1+s_2)\lambda)_{\tau}\subseteq  V_\tR(s_1\lambda)_{\tau}\otimes V_\tR(s_2\lambda)_{\tau}$.
Now with the same argument as in the proof of Lemma~\ref{pathvectorspecialrelation}, we can apply again 
Proposition~\ref{leadingtermauxbasis}. We conclude that $p_{\underline{a}',\ell}p_{\underline{b}',\ell}(v_{\underline c,\underline{\sigma}})\not=0$
for some $\underline{c}$ in $\mathrm{LS}_{\mathfrak C,\lam}^+(s_1+s_2)$ only if $m^{\underline{a}'+\underline{b}'}:=
\bigotimes_{\kappa\in A_\tau} m_\kappa^{\otimes(a'_\kappa+b'_\kappa)}$ occurs
in the expression of $v_{\underline c,\underline{\sigma}}$ in Proposition~\ref{leadingtermauxbasis}. 

The notation for $m^{\underline{a}'+\underline{b}'}$ looks ambiguous at a first glance. The sum $\underline{a}'+\underline{b}'$ makes sense in $\mathbb Q_{\geq 0}^{A_\tau}$,
but the support is now not anymore contained in a maximal chain. We fix an ordering of the tensor factors of $m^{\underline{a}'+\underline{b}'}$. Recall (compare to the comments before Proposition~\ref{leadingtermauxbasis}) that $M_{\tR}(\lam)^{\otimes{(s_1+s_2)\ell }}$
is multigraded, the homogeneous component of 
$M_{\tR}(\lam)^{\otimes{(s_1+s_2)\ell }}$ having the same multidegree as $m^{\underline{a}'+\underline{b}'}$ is just the span
of $m^{\underline{a}'+\underline{b}'}$. So if $m^{\underline{a}'+\underline{b}'}$ shows up in a presentation 
$v_{\underline c,\underline{\sigma}}$, then, up to multiplication by a unit, it must appear among one of the terms
in \eqref{summandsOfva}. 

But the tensor product is symmetric up to multiplying by roots of unities. Assume that $m^{\underline{a}'+\underline{b}'}$ shows up in \eqref{summandsOfva} with respect to the chosen ordering. Then for any other ordering to write down tensor factors in $m^{\underline{a}'+\underline{b}'}$, this pure tensor appears as one of the terms in \eqref{summandsOfva}. This implies: if $m^{\underline c}\RHD m^{\underline{a}'+\underline{b}'}$ with respect to one ordering of the tensor factors of $m^{\underline{a}'+\underline{b}'}$,
then $m^{\underline c}\RHD m^{\underline{a}'+\underline{b}'}$ with respect  any ordering of the
tensor factors of $m^{\underline{a}'+\underline{b}'}$.  Note that equality (up to rescaling) is not possible by assumption.

As a next step we want to show that $\underline{a}'\unrhd \underline{a}$ and $\underline{b}'\unrhd \underline{b}$ (and strict inequality for at least one term)
imply: there exists an ordering of the tensor factors such that $m^{\underline{a}'+\underline{b}'}\RHD m^{\underline{a}+\underline{b}}$.

Let $q$ be minimal such that there exist elements $\zeta_{q+1}<\ldots<\zeta_r=\tau$ in $A_\tau$ with the properties:
$$
\begin{array}{cc}
\{\xi\in\supp \underline{a}'\mid \ell(\xi)>q\},& \{\xi\in\supp \underline{a}\mid \ell(\xi)>q\},\\ 
\{\xi\in\supp \underline{b}'\mid \ell(\xi)>q\},&\{\xi\in\supp \underline{b}\mid \ell(\xi)>q\}
\end{array}
$$
are all subsets of $\{\zeta_{q+1},\ldots,\zeta_r=\tau\}$, and $a_{\zeta_j}=a'_{\zeta_j}$, $b_{\zeta_j}=b'_{\zeta_j}$ for $j\ge q+1$.

Since $\underline{a}'\unrhd \underline{a}$ and $\underline{b}'\unrhd \underline{b}$, the minimality of $q$ implies that we can find 
in $\supp \underline{a}'\cup\supp \underline{b}'$ an element $\zeta_q$ of length $q$. Without loss of generality we
assume $\zeta_q\in \supp \underline{a}'$ and $a'_{\zeta_q}>a_{\zeta_q}\ge 0$. If $b_{\zeta_q}=0$, then set
$$
m^{\underline{a}'+\underline{b}'}=m_{\zeta_r}^{\otimes(a_{\zeta_r}+b_{\zeta_r})}\otimes\cdots\otimes
m_{\zeta_{q+1}}^{\otimes(a_{\zeta_{q+1}}+b_{\zeta_{q+1}})}\otimes m_{\zeta_{q}}^{\otimes a'_{\zeta_{q}}}\otimes\cdots
$$
and compare the tensor product with:
$$
m^{\underline{a}+\underline{b}}=m_{\zeta_r}^{\otimes(a_{\zeta_r}+b_{\zeta_r})}\otimes\cdots\otimes
m_{\zeta_{q+1}}^{\otimes(a_{\zeta_{q+1}}+b_{\zeta_{q+1}})}\otimes m_{\zeta_{q}}^{\otimes a_{\zeta_{q}}}\otimes\cdots.
$$
Since $\zeta_{q}$ is greater than any of the elements in $\supp\underline{a}\cup \supp\underline{b}$ of length smaller than
$q$, $a'_{\zeta_{q}}>a_{\zeta_{q}}$ implies for the tensor products: $m^{\underline{a}'+\underline{b}'}\RHD m^{\underline{a}+\underline{b}}$.

If $b_{\zeta_q}>0$, then necessarily  $b'_{\zeta_q}\ge b_{\zeta_q}$, and the same arguments imply for 
$$
m^{\underline{a}'+\underline{b}'}=m_{\zeta_r}^{\otimes(a_{\zeta_r}+b_{\zeta_r})}\otimes\cdots\otimes
m_{\zeta_{q+1}}^{\otimes(a_{\zeta_{q+1}}+b_{\zeta_{q+1}})}\otimes m_{\zeta_{q}}^{\otimes a'_{\zeta_{q}}+ b'_{\zeta_{q}}}\otimes\cdots
$$
and 
$$
m^{\underline{a}+\underline{b}}=m_{\zeta_r}^{\otimes(a_{\zeta_r}+b_{\zeta_r})}\otimes\cdots\otimes
m_{\zeta_{q+1}}^{\otimes(a_{\zeta_{q+1}}+b_{\zeta_{q+1}})}\otimes m_{\zeta_{q}}^{\otimes a_{\zeta_{q}} + b_{\zeta_{q}}}\otimes\cdots
$$
that $m^{\underline{a}'+\underline{b}'}\RHD m^{\underline{a}+\underline{b}}$. By transitivity of $\RHD$, it follows 
$m^{\underline{c}} \RHD m^{\underline{a}+\underline{b}}$, which in turn implies $\underline{c} \rhd \underline{a}+\underline{b}$.
\end{proof}
By induction we get as an immediate consequence:

\begin{coro}\label{powerpathvector}
If $p_{\underline{a}}$, $\underline{a}\in \mathrm{LS}_{\mathfrak C,\lam}^+(s)$, is a path vector, then 
$p^m_{\underline{a}}$, $m\ge 1$, is a path vector associated to $m\underline{a}\in \mathrm{LS}_{\mathfrak C,\lam}^+(ms)$, up to multiplying by a root of unity.
\end{coro}

In the following we work over an algebraically closed field $\mathbb K$.

\begin{proposition}\label{seperating_path_vector}
Let $\underline{a}\in \mathrm{LS}_{\mathfrak C,\lam}^+(s)$ be such that $a_\tau\ge m$ for some integer $m>0$. If $p_{\underline{a}}$
is a path vector associated to $\underline{a}$, then  $p_{\underline{a}}=c p_\tau^m p_{\underline b}$, where $c$ is some root of unity, $\underline b=\underline{a}-m e_\tau\in \mathrm{LS}_{\mathfrak C,\lam}^+(s-m)$, $p_{\underline b}$ is a path vector
associated to $\underline b$ and $p_\tau$ is the dual of the extremal weight vector $v_\tau\in V(\lambda)_\tau$.
\end{proposition}

\begin{proof}
The proof is by decreasing induction along $\rhd$. First assume $\underline a=s e_\tau$ is the unique maximal element in $\mathrm{LS}_{\mathfrak C,\lam}^+(s)$.
By weight considerations we know that, up to multiplication by a non-zero scalar, $p_{\underline a}$ is equal to $p_\tau^s$, and
Proposition~\ref{product_path_vector} implies that the scalar is a root of unity.

Suppose now $\underline{a}\in \mathrm{LS}_{\mathfrak C,\lam}^+(s)$ is such that $a_\tau\ge m>0$. 
After having chosen an appropriate $\ell$ as in Corollary~\ref{differentEll},
we know $p_{\underline{a}}=p_{\underline{a},\ell} +\sum_{\underline{a}'\rhd \underline{a}} d_{\underline{a}',\ell}p_{\underline{a}',\ell}$.
Now $\underline{a}'\rhd \underline{a}$ implies $a'_\tau\ge a_\tau\ge m$, so by induction, we can write (up to rescaling by a root of unity) $p_{\underline{a}',\ell}$
as a product $p_\tau^m p_{\underline{b}'}$, where $\underline b'=\underline{a}'-m e_\tau\rhd \underline b$.

It remains to consider $p_{\underline{a},\ell}$. We can write  $\underline{a}=m e_\tau+\underline{b}$ where $m e_\tau\in \mathrm{LS}_{\mathfrak C,\lam}^+(m)$
and $\underline{b}\in \mathrm{LS}_{\mathfrak C,\lam}^+(s-m)$. By Lemma~\ref{pathvectorspecialrelation}, $p^m_{\tau}p_{\underline{b},\ell}$
is, up to rescaling by a root of unity, equal to 
$$
p_{\underline{a},\ell} +\sum_{\underline{a}''\rhd \underline{a}} d_{\underline{a}''}p_{\underline{a}''}
$$
Again by induction, we can write the $p_{\underline{a}''}$ as a product $p_\tau^m p_{\underline{b}''}$ where
$\underline b''=\underline{a}''-m e_\tau\rhd \underline b$. It follows: we can write (up to multiplication by a root of unity)
$$
p_{\underline{a}}=p_\tau^m\left(p_{\underline{b},\ell} +\sum_{\underline{b}'\rhd \underline{b}} c_{\underline{b}'}p_{\underline{b}'}\right),
$$
which proves the claim.
\end{proof}

By combining Corollary~\ref{powerpathvector} and Proposition~\ref{seperating_path_vector} we get:

\begin{coro}\label{factoring}
Let $\underline{a}\in \mathrm{LS}_{\mathfrak C,\lam}^+(s)$ be such that $a_\tau>0$ and let $N\in \mathbb N$ be a positive integer
such that $Na_\tau\in \mathbb N$. If $p_{\underline{a}}$ is a path vector associated to $\underline{a}$, then  
$p^N_{\underline{a}}$ is, up to multiplying by a root of unity, equal to $p_\tau^{Na_\tau} p_{\underline b}$, where
$\underline b=N\underline{a}-Na_\tau e_\tau\in \mathrm{LS}_{\mathfrak C,\lam}^+(Ns-Na_\tau)$ and $p_{\underline b}$ is a path vector
associated to $\underline b$.
\end{coro}

\printindex

\end{document}